\documentclass[11pt]{article}
\usepackage[utf8]{inputenc}
\usepackage{authblk}
\usepackage[margin=1in]{geometry}
\usepackage[titletoc,title]{appendix}
\usepackage[dvipsnames,table,xcdraw]{xcolor}
\usepackage{color}
\usepackage{setspace}
\usepackage{amsmath}
\usepackage{bbm}
\usepackage{amsmath,amsfonts,amssymb,amsthm,mathtools,bm}
\usepackage{algorithm2e}
\SetKwComment{Comment}{/* }{ */}
\SetKwInput{KwInput}{Input}
\SetKwFor{When}{when}{do}{}

\usepackage{graphicx,float}
\usepackage{subcaption}  
\usepackage{booktabs}
\usepackage{booktabs,multirow}
\usepackage{tikz-cd}
\usepackage{comment}

\SetKwProg{Init}{Initialize}{}{}
\RestyleAlgo{ruled}

\usepackage[colorlinks=true, pdfborder={0 0 0}, linkcolor=red]{hyperref}

\usepackage[style=authoryear,date=year,maxcitenames=2,maxbibnames=20,giveninits,url=false,backend=biber,uniquelist=false,uniquename=false,dashed=false,natbib=true,safeinputenc]{biblatex}
\bibliography{references}
\usepackage{wasysym,url}

\allowdisplaybreaks
\newtheorem{assumption}{Assumption}

\newtheorem{theorem}{Theorem}
\newtheorem{proposition}{Proposition}
\newtheorem{corollary}{Corollary}
\newtheorem{lemma}{Lemma}
\newtheorem{remark}{Remark}

\usepackage{amsmath,amsfonts,bm}

\def\eqref#1{equation~\ref{#1}}

\def\1{\bm{1}}

\DeclareMathAlphabet{\mathsfit}{\encodingdefault}{\sfdefault}{m}{sl}
\SetMathAlphabet{\mathsfit}{bold}{\encodingdefault}{\sfdefault}{bx}{n}

\newcommand{\E}{\mathbb{E}}

\newcommand{\R}{\mathbb{R}}

\usepackage{nicefrac}       
\newcommand{\rev}[1]{#1}

\makeatletter
\renewcommand{\eqref}[1]{\textup{\tagform@{\ref{#1}}}}
\makeatother
\usepackage{enumitem}

\newcommand{\Pp}{\mathbb{P}}
\newcommand{\F}{\mathcal{F}}
\newcommand{\one}{\mathbf{1}}
\newcommand{\norm}[1]{\left\|#1\right\|}
\newcommand{\abs}[1]{\left|#1\right|}
\newcommand{\ip}[2]{\left\langle #1,#2\right\rangle}
\newcommand{\op}{\mathrm{op}}

\newcommand{\dd}{\mathrm{d}}

\hypersetup{pdftitle={Exact Tail Laws for Nonlinear SGD with Infinite-Variance Gradients},
pdfauthor={Xiaoli Li, Qianqian Lei, Ziyang Wei, Wei Biao Wu}}

\begin{document}

\title{\textbf{Exact Tail Laws for Nonlinear SGD with Infinite-Variance Gradients}}
\author[1]{Xiaoli Li}
\author[1]{Qianqian Lei}
\author[1]{Ziyang Wei}
\author[1]{Wei Biao Wu}

\affil[1]{Department of Statistics, University of Chicago}

\maketitle

\begin{abstract}
With infinite-variance gradient noise, SGD can have infinite expected excess objective even when most iterates are close to the optimum. We characterize its accuracy through exact first-order tail probabilities for the last iterate of nonlinear, constant-stepsize SGD. For regularly varying innovations of index $t\in(1,2)$, the leading measure transports one large gradient shock through the local Hessian flow. This is the same measure that governs the central stable Ornstein--Uhlenbeck limit; we prove its relative accuracy at shrinking tolerances above the central scale, with state-dependent noise and a growing observation horizon.
The admissible window is nonempty for every available gradient moment $p>1$. The main technical step is a moment bound at the central fluctuation scale and a lower-moment control of the nonlinear error. The result yields explicit projection and excess-objective tails, extreme quantiles, and a factor $b^{1-t}$ for fixed batch size $b$. For affine population gradients, a sharper comparison gives the same law at fixed tolerances, including random sample Hessians. Coupled simulations in dimensions one to one hundred examine relative calibration, objective quantiles, batch effects, and the error due to random curvature.
\end{abstract}


\section{Introduction}\label{sec:intro}

The expected optimization error can be a poor description of an accurate iterate.
With additive infinite-variance gradient noise, even a strongly restoring SGD recursion has infinite expected excess objective after one update, while most iterates can remain close to the optimum.
Tail probabilities answer a more useful question in this setting: how often does the last iterate exceed a prescribed error tolerance, and how does that probability change with the stepsize, curvature, or batch size?
We derive exact first-order asymptotics for these probabilities, including their leading constants.

Consider iid samples $\varepsilon_n$, a differentiable objective $M:\mathbb R^d\to\mathbb R$, and a jointly measurable unbiased stochastic gradient $g$:
\begin{equation}\label{eq:population-problem}
 m(x)=\nabla M(x)=\E g(x,\varepsilon),\qquad
 x_n=x_{n-1}-\gamma g(x_{n-1},\varepsilon_n).
\end{equation}
The stepsize is constant within each run.
Let $x^\star$ be the target, $X_n=x_n-x^\star$, $H=\nabla m(x^\star)$, and $\eta=-g(x^\star,\varepsilon)$.
When $\eta$ is regularly varying with index $t\in(1,2)$ and radial-normalized tail measure $\mu_\eta$, our projected tail law is
\begin{equation}\label{eq:intro-main-tail}
 \Pp(\ell^\top X_{n_\gamma}>u_\gamma)
 \sim C_\ell\frac{\Pp(|\eta|>u_\gamma/\gamma)}{\gamma},\qquad
 C_\ell=\int_0^\infty\mu_\eta\{z:\ell^\top e^{-sH}z>1\}\,\dd s.
\end{equation}
Here $\ell\ne0$ is fixed, $C_\ell>0$, and $\sim$ means that the ratio tends to one.
A single large gradient shock has direction governed by $\mu_\eta$; the local flow $e^{-sH}$ transports a shock of age $s$, and $1/\gamma$ counts effective shock opportunities.
The vector result also specifies the direction of a large error and the tail of the excess objective.


The main technical difficulty is that the one-big-jump principle is most direct for linear recursions and weighted sums, whereas SGD is nonlinear
\rev{, and in \eqref{eq:intro-main-tail} the weights, the memory length, the horizon and the threshold all depend on $\gamma$}. Extreme value theory and regular variation provide precise tail tools for heavy-tailed sums, linear processes, and stochastic recursions \citep{BinghamGoldieTeugels1987,EmbrechtsKluppelbergMikosch1997,HultSamorodnitsky2008,MikoschWintenberger2013}.  Existing heavy-tailed SGD theory has established convergence guarantees, stable-limit behavior, and high-probability bounds under weak moment assumptions or clipping schemes \citep{WangEtAl2021,GorbunovDanilovaGasnikov2020,CutkoskyMehta2021,NguyenEtAl2023,BlanchetMijatovicYang2024}. 
\citet{Rhee_2023} develop sample-path large deviations and metastability for heavy-tailed dynamical systems with a factorized noise term $\sigma(x)Z$. Their results hold on fixed time horizons under global Lipschitz conditions. The untruncated case further assumes bounded drift and noise coefficients. Additional estimates are needed beyond a direct application of these fixed-window results. In contrast, our result identifies the exact shrinking-threshold tail probability of the nonlinear constant-stepsize iterate, including the sharp constant. 

\subsection{Contributions}
\begin{itemize}
\item \textbf{Relative tail asymptotics for nonlinear SGD.}
Under centered dissipativity and the stated regularity
assumptions, we establish a vector tail law for nonlinear
constant-stepsize SGD.
The result identifies the leading measure and gives vanishing
relative error at shrinking tolerances over a diverging
effective observation horizon.

\item \textbf{Moment bounds and nonlinear comparison.}
A central-scale moment bound and a lower-moment coupling
control the nonlinear error at the rare-event probability scale.
For a quadratic local remainder, these estimates yield
$(t-p)/t<a<(t-1)/t$, a nonempty threshold window for every
admissible gradient moment $p>1$.
This includes models whose noise has fewer
moments than the innovation at the target.
We also obtain corresponding windows for superlinear
H\"older remainders.

\item \textbf{Objective extremes and fixed-tolerance extensions.}
The tail measure yields projection and directional laws,
extreme objective quantiles, and fixed-batch scaling.
For affine population gradients, a sharper comparison extends
the law to fixed tolerances under an additional moment condition,
including models with random sample Hessians.
\end{itemize}

\subsection{Related literature}

\textbf{SGD asymptotics and uncertainty quantification.}
Classical stochastic approximation traces back to \citet{RobbinsMonro1951}, with asymptotic efficiency of averaged recursions developed by \citet{PolyakJuditsky1992} and in the stochastic approximation monograph of \citet{KushnerYin2003}.  Non-asymptotic convergence theory for stochastic approximation and SGD is now extensive; see, for example, \citet{MoulinesBach2011}.  More closely related to the inferential motivation of this paper, \citet{ChenLeeTongZhang2020} develop statistical inference for model parameters computed by SGD, \citet{FangXuYang2018} propose online bootstrap confidence intervals, \citet{YuEtAl2021} analyze constant-stepsize SGD inference in nonconvex regimes, and \citet{SheshukovaEtAl2025} establish Gaussian approximation and multiplier bootstrap validity for averaged SGD, while \citet{wei2025gaussian} develop Gaussian approximation and concentration results for constant-learning-rate SGD. These works are primarily finite-variance and Gaussian in nature. Most distributional analyses of SGD rely, either explicitly or implicitly, on a finite-variance picture~\citep{Zhu2023,merad_convergence_2023}. 

\textbf{Heavy-tailed SGD and stochastic optimization.}
Heavy-tailed stochastic gradient noise has been observed and modeled in modern learning problems \citep{SimsekliSagunGurbuzbalaban2019,GurbuzbalabanSimsekliZhu2021,HodgkinsonMahoney2021,DamekMentemeier2026}.  Theoretical work has studied convergence rates and limiting distributions of SGD under infinite-variance noise \citep{WangEtAl2021,BlanchetMijatovicYang2024}, as well as robust or clipped methods that achieve high-probability guarantees under weak moment assumptions \citep{GorbunovDanilovaGasnikov2020,CutkoskyMehta2021,NguyenEtAl2023, LiuZhou2025}.  These results show that optimization remains possible under heavy-tailed noise and provide convergence or high-probability control.  

\textbf{Extreme value theory and regular variation.}
The mathematical mechanism behind our result is the one-big-jump principle for regularly varying random variables.  Regular variation and Potter-type bounds \citep{BinghamGoldieTeugels1987} are classical tools in extreme value theory \citep{BinghamGoldieTeugels1987,EmbrechtsKluppelbergMikosch1997}. \citet{DenisovDiekerShneer2008} characterize uniform big-jump domains for random walks. Tail asymptotics for heavy-tailed linear processes, infinite weighted series, and stochastic recurrence equations have been developed in \citet{HultSamorodnitsky2008} and \citet{MikoschWintenberger2013}. \rev{\citet{Rhee_2023} develop sample-path large deviations for recursions with state-dependent multiplicative noise under global Lipschitz conditions; related work studies exit times and multi-basin dynamics \citep{Rhee_2025, Rhee_2026}.} Our linearized recursion is related to these weighted-sum settings, but the learning-rate-dependent weights and the nonlinear SGD transfer require additional analysis.  In particular, the sharp constant in \eqref{eq:intro-main-tail} arises from a Riemann-sum limit of the exponentially stable impulse response of the linearized dynamics.

\textbf{Nonconvex and dissipative constant-stepsize dynamics.}
Constant-stepsize SGD can be viewed as a Markov chain, and this viewpoint has been used to study asymptotic expansions, bias, and stationary behavior \citep{DieuleveutDurmusBach2020,YuEtAl2021}.  Dissipativity assumptions are also standard in nonconvex sampling and learning dynamics, including analyses of stochastic gradient Langevin dynamics \citep{RaginskyRakhlinTelgarsky2017,ZhangAkyildizDamoulasSabanis2023}.  Our assumptions follow this stability perspective but target a different phenomenon: rather than Gaussian fluctuations or stationary bias, we quantify the exact heavy-tailed rare-event probability of the last iterate.

\section{Setting and main tail results}\label{sec:nonlinear-transfer}

\paragraph{Notation.}
The dimension $d$ is fixed, $|\cdot|$ is the Euclidean norm, and $\|\cdot\|_{\op}$ is the matrix operator norm.
For a random vector $Z$, let $\|Z\|_r=(\E|Z|^r)^{1/r}$.
We write $I_d$ for the identity and $S^{d-1}=\{z:|z|=1\}$.
The initial point $x_0$ is fixed and deterministic; $\mathcal F_n=\sigma(\varepsilon_1,\ldots,\varepsilon_n)$ is the natural filtration.
Constants may depend on the fixed model, dimension, and exponents, but not on $n$ or $\gamma$.
The exponent $p\in(1,2)$ denotes an available gradient moment, and $t\in(1,2)$ denotes the tail index.

\begin{assumption}[Centered dissipativity]\label{ass:centered-diss}
The continuously differentiable objective $M$ has gradient $m$, $m(x^\star)=0$, and, for some $\alpha>0$,
\begin{equation}\label{eq:centered-dissipativity}
 \langle x-x^\star,m(x)\rangle\ge\alpha|x-x^\star|^2,\qquad x\in\mathbb R^d.
\end{equation}
\end{assumption}

Assumption~\ref{ass:centered-diss} is a one-point monotonicity condition centered at the target $x^\star$. 
It is weaker than global strong convexity, as it does not require monotonicity between arbitrary pairs of points, but it is stronger than a purely tail dissipativity condition since \eqref{eq:centered-dissipativity} is imposed globally. This type of stability condition is common in the analysis of stochastic approximation and dissipative dynamical systems, and it is compatible with many nonconvex statistical objectives~\citep{YuEtAl2021}. 
The nonconvex models in Section~\ref{sec:experiments} satisfy this condition explicitly.

\begin{assumption}[$L^p$ gradient regularity]\label{ass:noise-regularity}
For some $1<p<2$, $L_p,L_g<\infty$, and $\kappa>0$,
\begin{align}
 \|g(x,\varepsilon)\|_p&\le L_p(1+|x-x^\star|),\quad x\in\mathbb R^d,
 \label{eq:global-lp-growth}\\
 \|g(x,\varepsilon)-g(x^\star,\varepsilon)\|_p&\le L_g|x-x^\star|,
 \quad |x-x^\star|\le\kappa.
 \label{eq:local-centered-lip}
\end{align}
\end{assumption}
The gradient has global $L^p$ growth and local $L^p$ increments; sample gradients may have random Lipschitz constants.

\begin{assumption}[Local quadratic remainder]\label{ass:local-quadratic}
The gradient is differentiable at $x^\star$, with $H=\nabla m(x^\star)$, and for some $C_R<\infty$,
\begin{equation}\label{eq:def-H-general}
 |m(x)-H(x-x^\star)|\le C_R|x-x^\star|^2,\qquad |x-x^\star|\le\kappa.
\end{equation}
\end{assumption}
Here $H$ is symmetric and positive definite, with $H\succeq\alpha I_d$.
A locally Lipschitz Hessian is sufficient.
Remark~\ref{remark:holder} discusses any superlinear remainder with a corresponding threshold window.

\begin{assumption}[Regularly varying innovation]\label{ass:vector-rv}
The centered innovation $\eta=-g(x^\star,\varepsilon)$ has index $t\in(1,2)$ and a nonzero Radon tail measure $\mu_\eta$ on $\mathbb R^d\setminus\{0\}$ satisfying
\begin{equation}\label{eq:mrv}
 \frac{\Pp(r^{-1}\eta\in A)}{\Pp(|\eta|>r)}
 \longrightarrow\mu_\eta(A),\qquad r\to\infty,
\end{equation}
for every Borel set $A$ bounded away from zero with $\mu_\eta(\partial A)=0$.
The measure satisfies $\mu_\eta(cA)=c^{-t}\mu_\eta(A)$ for $c>0$ and $\mu_\eta\{|z|>1\}=1$.
\end{assumption}


Spherical Student-$t$ noise with $t$ degrees of freedom is an example that has $\Pp(|\eta|>r)\sim c_\eta r^{-t}$ for a constant $c_\eta>0$ and a rotation-invariant $\mu_\eta$. The two objects in \eqref{eq:intro-main-tail} are
\begin{equation}\label{eq:Lgamma}
 L_\gamma(u):=\frac{\Pp(|\eta|>u/\gamma)}{\gamma},\qquad
 \nu_H(A):=\int_0^\infty\mu_\eta\{z:e^{-sH}z\in A\}\,\dd s.
\end{equation}
One unit of effective time contains $1/\gamma$ iterations, and each noise step exceeds $u$ with probability $\Pp(|\eta|>u/\gamma)$, so $L_\gamma(u)$ is the rate of large shocks. The measure $\nu_H$ is finite on sets bounded away from zero and satisfies $\nu_H(cA)=c^{-t}\nu_H(A)$.

\begin{theorem}[Nonlinear vector tail]\label{thm:nonlinear-exact-tail}
Assume Assumptions~\ref{ass:centered-diss}--\ref{ass:vector-rv}, with
\begin{equation}\label{eq:nonlinear-parameter-window}
 1<p\le t,\qquad
 \frac{t-p}{t}<a<\frac{t-1}{t},\qquad u_\gamma=c\gamma^a,\quad c>0.
\end{equation}
As $\gamma\downarrow0$, let $n_\gamma\in\mathbb N$ satisfy $n_\gamma\gamma\to\infty$.
If $x_0\ne x^\star$, also require $n_\gamma\gamma/\log(1/\gamma)\to\infty$.
For every fixed Borel set $A$ satisfying
\begin{equation}\label{eq:A-continuity}
 \operatorname{dist}(A,0)>0,\qquad \nu_H(\partial A)=0,
\end{equation}
we have
\begin{equation}\label{eq:nonlinear-vector-tail}
 \Pp(X_{n_\gamma}/u_\gamma\in A)
 =L_\gamma(u_\gamma)\{\nu_H(A)+o(1)\}.
\end{equation}
If $\nu_H(A)>0$, this is a relative asymptotic equivalence.
\end{theorem}

For a half-space $A_\ell=\{z:\ell^\top z>1\}$, its boundary has zero tail measure.
Thus~\eqref{eq:nonlinear-vector-tail} gives~\eqref{eq:intro-main-tail}, with
\begin{equation}\label{eq:linear-tail-constant}
 C_\ell=\nu_H(A_\ell).
\end{equation}
The theorem also applies to norm exceedances and to directional sets
$\{z:|z|>1,\ z/|z|\in C\}$ when their boundaries have zero $\nu_H$ mass.
These sets describe both the size and direction of extreme deviations.
More precisely, writing $D_1=\{z:|z|>1\}$,
\begin{equation}\label{eq:conditional-extremes}
 \mathcal L(X_{n_\gamma}/u_\gamma\mid |X_{n_\gamma}|>u_\gamma)
 \Longrightarrow\frac{\nu_H(\,\cdot\cap D_1)}{\nu_H(D_1)}.
\end{equation}
The symbol $\mathcal L$ denotes a probability law and $\Longrightarrow$ weak convergence.
Mapping $z$ to $z/|z|$ gives the conditional direction of a large error (Appendix~\ref{app:conditional-law}).

\begin{corollary}[Eigenvector projections]\label{cor:nonlinear-left-eigen-tail}
Under Theorem~\ref{thm:nonlinear-exact-tail}, let $|v|=1$, $H^\top v=\lambda_vv$, and
$\mu_\eta\{z:v^\top z>1\}>0$.
Then
\begin{equation}\label{eq:nonlinear-eigen-tail}
 \Pp(v^\top X_{n_\gamma}>u_\gamma)
 \sim\frac{\Pp(v^\top\eta>u_\gamma/\gamma)}{t\lambda_v\gamma}.
\end{equation}
\end{corollary}
For a general direction, the propagated direction $e^{-sH^\top}\ell$ varies with $s$, so the full measure in~\eqref{eq:linear-tail-constant} matters.
The eigenvector formula isolates the inverse dependence on local curvature.
Appendix~\ref{app:transfer-proof} also gives a projected-only result, with its own sufficient moment and threshold conditions.

\begin{corollary}[Excess objective value]\label{cor:risk-tail}
Under Theorem~\ref{thm:nonlinear-exact-tail}, define
$A_M=\{z:z^\top Hz/2>1\}$ and $C_M=\nu_H(A_M)$.
Then $0<C_M<\infty$ and
\begin{equation}\label{eq:risk-tail}
 \Pp\{M(x_{n_\gamma})-M(x^\star)>u_\gamma^2\}
 \sim C_M L_\gamma(u_\gamma).
\end{equation}
\end{corollary}
This is a tail statement for the optimization objective itself.
For a pure power innovation tail, the objective threshold has exponent $t/2$ in this joint limit.
Tail probabilities remain informative even when expected excess objective value is infinite, as in additive infinite-variance quadratic models.
The resulting quantiles describe increasingly rare objective deviations.

\begin{remark}[Reading the constants]
    Write the radial-normalized tail measure in polar form as
$\mu_\eta(\dd r,\dd\omega)=t r^{-t-1}\dd r\,S_\eta(\dd\omega)$,
where $S_\eta$ is a probability measure on $S^{d-1}$.
For the special case $H=\lambda I_d$, direct radial integration gives
\begin{equation}\label{eq:isotropic-constants}
 C_\ell=\frac{\int(\ell^\top\omega)_+^t S_\eta(\dd\omega)}{t\lambda},
 \qquad
 C_M=\frac{2^{-t/2}}{t}\lambda^{t/2-1},
 \qquad r_+=\max(r,0).
\end{equation}
Indeed, the flow multiplies a shock of age $s$ by $e^{-\lambda s}$, and the objective event requires its propagated radius to exceed $\sqrt{2/\lambda}$.
The projected constant depends on the angular distribution, whereas this isotropic objective constant does not.
Holding the innovation tail law fixed, increasing local curvature suppresses projected tails as $\lambda^{-1}$ and objective-value tails as $\lambda^{t/2-1}$.
The difference comes from the objective's curvature weighting.
For a nonscalar $H$, directions can change through unequal contraction rates, and~\eqref{eq:Lgamma} retains the full interaction between curvature and tail geometry.
\end{remark}

\begin{remark}[Three scales]
    For Student innovations, typical fluctuations have scale
proportional to $\gamma^{1-1/t}$.
The upper bound on $a$ places $u_\gamma$ above this scale,
while $a>0$ localizes the event near $x^\star$.
The lower bound on $a$ makes the nonlinear comparison error
negligible relative to the tail probability.
For general regularly varying innovations, $L_\gamma$
retains the slowly varying factor.

The propagated measure $\nu_H$ is also the L\'evy measure
of the stationary distribution of a stable
Ornstein--Uhlenbeck process with drift $z\mapsto-Hz$
and driving L\'evy measure $\mu_\eta$
\citep{GlowackiKulikSoulier2026}.
Theorem~\ref{thm:nonlinear-exact-tail} shows that this measure
governs shrinking-threshold tail probabilities of nonlinear
SGD with vanishing relative error.
\end{remark}

\section{Why the linearization is relatively accurate}\label{sec:linear-subquadratic}

We separate the argument into uniform moment control, a linear single-jump limit, and a nonlinear comparison.
This order exposes where each assumption enters.
All statements are uniform over $n$ for sufficiently small $\gamma$ when a finite-time bound is asserted.

\subsection{Central-scale moments and the coupled recursion}\label{sec:lp-moment}
Choose the innovation quantile and central scale
\begin{equation}\label{eq:central-scale}
 b_\gamma=\inf\{r\ge1:\Pp(|\eta|>r)\le\gamma\},\qquad
 s_\gamma=\gamma b_\gamma.
\end{equation}
Regular variation gives $\Pp(|\eta|>b_\gamma)\sim\gamma$ and makes $s_\gamma$ regularly varying at zero with index $1-1/t$.
For Student innovations, $s_\gamma$ is proportional to $\gamma^{1-1/t}$.

\begin{theorem}[Central-scale moment bound]\label{thm:central-moment}
Under Assumptions~\ref{ass:centered-diss}, \ref{ass:noise-regularity}, and~\ref{ass:vector-rv}, with $1<p<t$, there exist $C,c,\gamma_0>0$ such that, for all $0<\gamma<\gamma_0$ and integers $n\ge0$,
\begin{equation}\label{eq:central-moment}
 \E|X_n|^p\le Ce^{-cn\gamma}|X_0|^p+Cs_\gamma^p.
\end{equation}
\end{theorem}
The regularized Lyapunov function $V_\gamma(x)=(s_\gamma^2+|x|^2)^{p/2}$ resolves the singular curvature of $|x|^p$ near zero.
Its Taylor remainder uses truncated second moments for small innovations and $p$th moments for large ones.
Both contributions are $O(\gamma s_\gamma^p)$.
The state-dependent noise term is absorbed into the dissipative drift.
Appendix~\ref{app:central-proof} gives the proof.
Without regular variation, the weaker $O(\gamma^{p-1})$ bound remains available under global moment growth alone (Appendix~\ref{app:moment-proof}).

For linear drift $m(x)=\lambda x$, zero initialization, and symmetric $t$-stable innovation with $\E e^{i\theta\eta}=e^{-|\theta|^t}$,
\[
 X_n\overset d=\gamma\left\{\sum_{j<n}(1-\lambda\gamma)^{tj}\right\}^{1/t}\eta.
\]
Hence $\E|X_n|^p\sim\E|\eta|^p(t\lambda)^{-p/t}\gamma^{p(1-1/t)}$ when $n\gamma\to\infty$.
The moment theorem attains this polynomial order, improving the coarse exponent by $1-p/t>0$.

Using the same innovations, define
\begin{equation}\label{eq:eta-def}
 Y_n=B_\gamma Y_{n-1}+\gamma\eta_n,\quad
 B_\gamma=I_d-\gamma H,\quad Y_0=X_0,\quad \eta_n=-g(x^\star,\varepsilon_n).
\end{equation}
Then $\E\eta_n=0$ and
\begin{equation}\label{eq:vector-linear-iteration}
 Y_n=B_\gamma^nX_0+\gamma T_{\gamma,n},\qquad
 T_{\gamma,n}=\sum_{j=0}^{n-1}B_\gamma^j\eta_{n-j}.
\end{equation}
Write $y_n=x^\star+Y_n$ and $\Delta_n=X_n-Y_n=x_n-y_n$.

\subsection{Linear single-jump limit}
\begin{assumption}[Linear rare-event regime]\label{ass:vector-threshold}
The positive sequences $u_\gamma$ and integers $n_\gamma\to\infty$ satisfy
\begin{equation}\label{eq:vector-thresholds}
 n_\gamma\gamma\to\infty,\quad
 r_\gamma=u_\gamma/\gamma\to\infty,\quad
 L_\gamma(u_\gamma)\to0,\quad
 |B_\gamma^{n_\gamma}X_0|/u_\gamma\to0.
\end{equation}
\end{assumption}
\begin{theorem}[Linear vector tail]\label{thm:vector-linear-tail}
Assume Assumption~\ref{ass:centered-diss}, differentiability of $m$ at $x^\star$, and
Assumptions~\ref{ass:vector-rv} and~\ref{ass:vector-threshold}.
For $A$ satisfying~\eqref{eq:A-continuity},
\begin{equation}\label{eq:vector-T-tail-asymp}
 \Pp(Y_{n_\gamma}/u_\gamma\in A)
 =L_\gamma(u_\gamma)\{\nu_H(A)+o(1)\}.
\end{equation}
If $\nu_H(A)>0$, this probability is also asymptotic to
$\sum_{j<n_\gamma}\Pp(\gamma B_\gamma^j\eta\in u_\gamma A)$.
\end{theorem}

\begin{corollary}[Linear-functional constant]\label{cor:linear-functional-tail}
Under Theorem~\ref{thm:vector-linear-tail}, let $\ell\in\mathbb R^d\setminus\{0\}$ be fixed with $C_\ell>0$. Then
\begin{equation}\label{eq:linear-functional-Y-tail-asymp}
 \Pp(\ell^\top Y_{n_\gamma}>u_\gamma)
 \sim \Lambda_\gamma(\ell,u_\gamma)
 \sim C_\ell L_\gamma(u_\gamma),\quad
 \Lambda_\gamma(\ell,u)=\sum_{j<n_\gamma}\Pp(\gamma\ell^\top B_\gamma^j\eta>u).
\end{equation}
\end{corollary}
The finite sum uses the exact innovation distribution when available.
The weights converge to the Hessian flow, turning the sum into the integral for $\nu_H$; truncated moments control the small jumps uniformly over the horizon.
Appendix~\ref{app:linear-proof} gives the vector and projected-tail proofs.
The linear theorem also allows fixed positive tolerances.

\subsection{A lower-moment nonlinear transfer}\label{subsec:transfer-mechanism}
\begin{theorem}[Central-scale linearization error]\label{thm:central-delta}
Under Assumptions~\ref{ass:centered-diss}--\ref{ass:vector-rv}, with $1<p<t$, for every $1\le q\le p$ there are $C_q,c_q,\gamma_q>0$ such that
\begin{equation}\label{eq:central-delta}
 \E|\Delta_n|^q\le C_qs_\gamma^p
 +C_qe^{-c_qn\gamma}(|X_0|^p+|X_0|^q),
 \qquad n\ge1,\quad 0<\gamma<\gamma_q.
\end{equation}
\end{theorem}

The first-moment case $q=1$, combined with the linear tail law, yields Theorem~\ref{thm:nonlinear-exact-tail}; the proof is given in Appendix~\ref{app:transfer-proof}. When the gradient growth and local increment conditions hold for every $1<p<t$, choosing a fixed $p$ sufficiently close to $t$ covers every fixed threshold $u_\gamma=c\gamma^a$ with $0<a<(t-1)/t$.


\begin{remark}[Superlinear Remainders]\label{remark:holder}
The local quadratic remainder Assumption~\ref{ass:local-quadratic} can be extended to the H\"older local remainders, i.e. assuming 
$|m(x)-H(x-x^\star)|\le C|x-x^\star|^\beta$ for $|x-x^\star|\le\kappa$ instead,
where $1<\beta\le2$. This further extends the scope; see Appendix~\ref{app:extensions}. 
\end{remark}
\section{Consequences and extensions}\label{sec:consequences}

The tail law quantifies extreme objective values and the effects of stepsize, batching, and observation time.
We also describe when it extends to fixed tolerances and how local curvature emerges from nonlinear dynamics.

\textbf{How large are extreme objective values?}
Quantiles describe the upper tail even when the expected objective is infinite.
Suppose $\Pp(|\eta|>r)\sim c_\eta r^{-t}$ with $c_\eta>0$, and define
$Q_\gamma(1-\delta)=\inf\{q:\Pp(M(x_{n_\gamma})-M(x^\star)>q)\le\delta\}$.
For $1<p<t$, Corollary~\ref{cor:risk-tail} gives
\begin{equation}\label{eq:main-extreme-quantile}
 Q_\gamma(1-\delta_\gamma)
 \sim\left(\frac{C_Mc_\eta\gamma^{t-1}}{\delta_\gamma}\right)^{2/t},
 \qquad
 \delta_\gamma=\gamma^\zeta,\quad 0<\zeta<p-1.
\end{equation}
This formula relates an extreme loss level to its exceedance probability.
The restriction on $\zeta$ places the corresponding iterate threshold within the nonlinear theorem's admissible window.
Appendix~\ref{app:extensions} gives the inversion argument and extensions to slowly varying tails and H\"older remainders.

\textbf{What do smaller steps and larger batches change?}
Averaging a fixed number $b$ of independent gradients multiplies the leading innovation tail probability by $b^{1-t}$.
For pure power tails, the dependence on stepsize and batch size is therefore $(\gamma/b)^{t-1}$ at a common admissible tolerance.
For example, when $t=1.5$, batches of four and sixteen reduce the leading tail probability by factors of $1/2$ and $1/4$.
Comparing $(\gamma,b)$ with $(\gamma/b,1)$ also requires appropriate observation horizons; Appendix~\ref{app:extensions} states the equal-sample-budget condition.

\textbf{When can the tolerance remain fixed?}
The main nonlinear theorem uses shrinking tolerances.
For an affine population gradient, the following extension also covers any fixed $u>0$.

\begin{corollary}[Affine gradients at fixed tolerances]\label{cor:affine-fixed}
Assume Assumptions~\ref{ass:centered-diss}, \ref{ass:noise-regularity}, and~\ref{ass:vector-rv}, with $m(x)=H(x-x^\star)$ globally and $t^2/(2t-1)<p<t$.
Under the initialization and horizon conditions of Theorem~\ref{thm:nonlinear-exact-tail}, for every fixed $u>0$ and $A$ satisfying~\eqref{eq:A-continuity},
\begin{equation}\label{eq:affine-fixed-tail}
 \Pp(X_{n_\gamma}\in uA)
 =L_\gamma(u)\{\nu_H(A)+o(1)\}.
\end{equation}
The conclusion also holds when the complete gradient conditions hold at the finite endpoint $p=t$.
\end{corollary}

Sample gradients may still have random curvature, as in least squares with bounded random designs and heavy-tailed labels.
At $t=1.5$, the sufficient moment condition is $p>9/8$.
Appendix~\ref{app:affine-fixed} gives the proof.

\textbf{How does observation time affect the tail?}
Starting from $X_0=0$, a finite horizon $n_\gamma\gamma\to T\in(0,\infty)$ replaces the integral in $\nu_H$ by $\int_0^T$.
For a unit eigenvector $Hv=\lambda_vv$, the resulting constant is
$C_v(T)=C_v(1-e^{-t\lambda_vT})$, where
$C_v=\nu_H(\{z:v^\top z>1\})$.
This factor describes how the tail approaches its long-time level.

\textbf{Why does local curvature describe nonlinear dynamics?}
For a locally Lipschitz nonlinear drift, let $\varphi_s(z)$ solve
$\dot\varphi_s=-m(x^\star+\varphi_s)$ with $\varphi_0(z)=z$.
The flow transports a large shock through the nonlinear dynamics, giving the deterministic quantity
$C_A^{\rm flow}(u)=\int_0^\infty\mu_\eta\{z:\varphi_s(z)\in uA\}\,\dd s$.
Under its stated assumptions, Proposition~\ref{prop:flow-localization} proves
$u^tC_A^{\rm flow}(u)\to\nu_H(A)$ as $u\downarrow0$.
Thus the nonlinear flow recovers the local Hessian description at shrinking tolerances.
\section{Numerical studies}\label{sec:experiments}

We examine three questions: whether the tail law calibrates nonlinear SGD in its joint limit, what it predicts for extreme objective values and batching, and how random sample curvature affects finite-step accuracy.
Coupled nonlinear and population-linear iterates $X,Y$ share the same innovations.
For a projected event, write
\begin{equation}\label{eq:experiment-ratios}
 R_{\rm lin}=\frac{\widehat p_Y}{\Lambda_\gamma},\qquad
 R_{\rm tr}=\frac{\widehat p_X}{\widehat p_Y},\qquad
 R_{\rm tot}=\frac{\widehat p_X}{\Lambda_\gamma}=R_{\rm lin}R_{\rm tr},\qquad
 R_{\rm const}=\frac{\Lambda_\gamma}{C_\ell L_\gamma(u)},
\end{equation}
where $\Lambda_\gamma$ is the exact-distribution single-jump sum in Corollary~\ref{cor:linear-functional-tail}. A ratio of one indicates agreement with the benchmark. Intervals are nominal 95\% pointwise Monte Carlo intervals; paired ratios retain the covariance within a coupled trajectory.

\subsection{Tail calibration in a nonconvex model}

The nonlinear theorem concerns shrinking thresholds; fixed-threshold sweeps diagnose finite-level effects. We use
\begin{equation}\label{eq:scalar-diagonal-model}
 m(x)=x+\frac{4\sin^3x}{1+x^2},\qquad
 g(x,\varepsilon)=m(x)-\sigma\varepsilon,\qquad
 \varepsilon\sim t_{1.5},\quad\sigma=0.02.
\end{equation}
Here $x^\star=0$, $H=1$, $xm(x)>0.882x^2$, and $m'(2)\simeq-0.307$.
The objective is nonconvex and the assumptions hold for every $1<p<1.5$.
For $\gamma\in\{0.01,0.003,0.001\}$, we use $600{,}000$ independent terminal pairs, $X_0=Y_0=0$, and $n=\lceil5\log(1/\gamma)/\gamma\rceil$.

\begin{figure}[!htb]
\centering
\includegraphics[width=\linewidth]{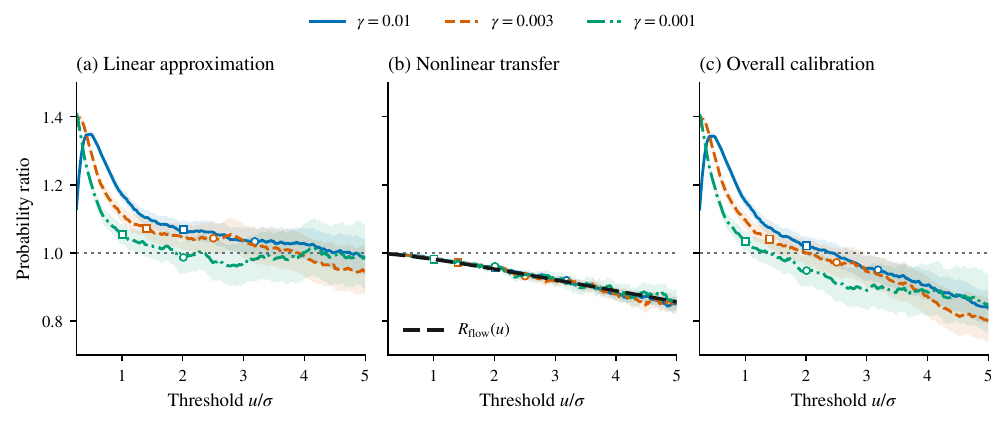}
\caption{Nonconvex tail calibration on 241 thresholds, using the same $600{,}000$ terminal pairs per stepsize. Axes are shared; shading gives pointwise intervals. Circles and squares mark $u_\gamma=0.16\gamma^a$, $a=0.2,0.3$. The black dashed curve in (b) is the deterministic flow diagnostic $R_{\rm flow}(u)$.}
\label{fig:threshold-sweep}
\end{figure}

Figure~\ref{fig:threshold-sweep}(b) compares nonlinear transfer with
$R_{\rm flow}(u)=tu^t\int_u^\infty r^{-t}/m(r)\,\dd r$, which tends to one as $u\downarrow0$ by Proposition~\ref{prop:flow-localization}.
This is a deterministic diagnostic.
It tracks the transfer deficit within $0.01$ at all six prespecified points (Table~\ref{tab:shrinking-main}); its finite-step counterpart does likewise and differs by less than $0.001$ over the displayed range (Appendix~\ref{app:flow-localization}).
The joint-limit test uses $a=0.2,0.3$, both in the theorem's window at $p=1.45$; thresholds are respectively $14.8$--$20.1$ and $9.33$--$10.1$ times $\sigma\gamma^{1/3}$.

\begin{table}[!htb]
\centering\small
\caption{All prespecified settings, with the flow diagnostic. The objective column compares $M(X_n)-M(0)>u_\gamma^2$ with Corollary~\ref{cor:risk-tail}. Brackets give pointwise intervals.}
\label{tab:shrinking-main}
\begin{tabular}{cccccc}
\toprule
$a$ & $\gamma$ & $R_{\rm flow}$ & $R_{\rm tr}$ & $R_{\rm tot}$ [95\% CI] & Objective ratio [95\% CI]\\
\midrule
0.2 & 0.01 & 0.915 & 0.920 & 0.951 [0.914, 0.988] & 0.884 [0.852, 0.918]\\
0.2 & 0.003 & 0.937 & 0.932 & 0.972 [0.931, 1.016] & 0.944 [0.907, 0.983]\\
0.2 & 0.001 & 0.953 & 0.961 & 0.948 [0.903, 0.996] & 0.907 [0.866, 0.950]\\
\midrule
0.3 & 0.01 & 0.953 & 0.956 & 1.021 [0.994, 1.048] & 0.973 [0.949, 0.998]\\
0.3 & 0.003 & 0.971 & 0.972 & 1.041 [1.013, 1.069] & 1.007 [0.982, 1.033]\\
0.3 & 0.001 & 0.981 & 0.982 & 1.034 [1.006, 1.064] & 0.993 [0.968, 1.020]\\
\bottomrule
\end{tabular}
\end{table}

Nonlinear transfer improves along both sequences; at the smallest step the total ratios are within $5.2\%$ of one.
The $a=0.2$ total ratios are nonmonotone because linear and transfer errors contribute in opposite directions.
For $M(x)=x^2/2+4\int_0^x\sin^3(s)/(1+s^2)\,\dd s$, the smallest-step objective ratios are $0.907$ and $0.993$; slower threshold shrinkage retains a larger correction.
The constant correction $R_{\rm const}-1$ is $0.025\%$--$0.25\%$ (full counts in Appendix~\ref{app:diagonal}).

\subsection{Extreme objective values and batching}

We next use $d=10$, $H=\operatorname{diag}(1,1+1/9,\ldots,2)$, $\boldsymbol a=(1,\ldots,1)^\top/\sqrt{10}$, and
\[
 m(x)=Hx+6\boldsymbol a\frac{\sin^3(\boldsymbol a^\top x)}{1+(\boldsymbol a^\top x)^2},\qquad
 M(x)=\tfrac12x^\top Hx+6\int_0^{\boldsymbol a^\top x}\frac{\sin^3s}{1+s^2}\,\dd s.
\]
This objective is nonconvex and centered-dissipative (Appendix~\ref{app:multidim-models}).
Innovations have the multivariate Student-$t_{1.5}$ law; each step averages $b$ independent samples.
Six configurations use $300{,}000$ coupled terminal pairs each, from zero, at $n=\max\{500,\lceil5\log(1/\gamma)/\gamma\rceil\}$.

The full-objective single-jump prediction is
\begin{equation}\label{eq:main-risk-benchmark}
 B_b(\epsilon)=b\sum_{j<n}\Pp\{M((\gamma/b)(I_d-\gamma H)^j\eta)-M(0)>\epsilon\}.
\end{equation}
It evaluates each linearly transported shock with the exact objective.
Appendix~\ref{app:merged-risk} separately reports the local quadratic approximation, the quadratic-event single-jump sum, and their numerical standard errors.
The matched heavy-tail reference evaluates the same objective on $Y_n$.
For comparison, we also fit a covariance-based Gaussian surrogate using $10^5$ independent innovations and repeat this fit five times.
Since the true covariance is infinite, this diagnostic tests what is lost when a finite empirical covariance is used to describe the tail.

\begin{figure}[!htb]
\centering
\includegraphics[width=\linewidth]{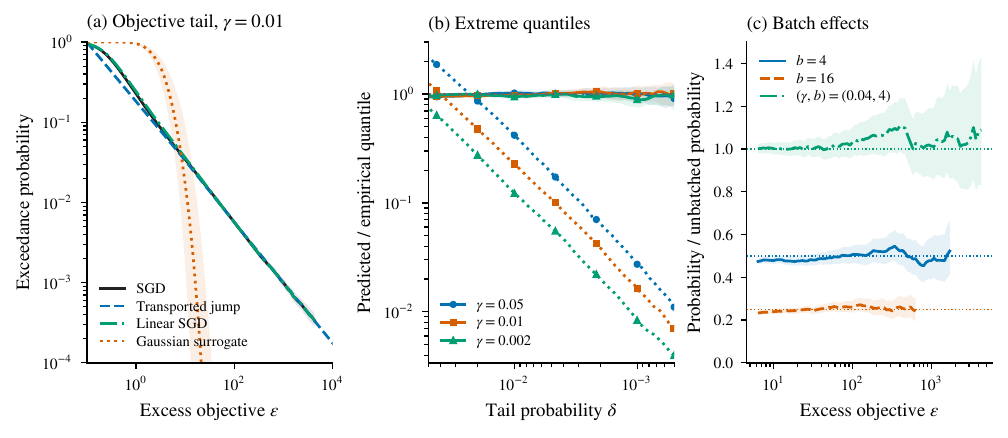}
\caption{Objective tails and batch effects. (a) Full-objective exceedances at $(\gamma,b)=(0.01,1)$; the Gaussian band spans five covariance fits. (b) Predicted/empirical quantiles for $b=1$: solid curves use~\eqref{eq:main-risk-benchmark}, dotted curves the median Gaussian fit. 
(c) Tail probabilities relative to $(0.01,1)$; horizontal references are $1/2$, $1/4$, and $1$. Empirical tail displays require at least 100 exceedances in each relevant sample; panel (c) also requires baseline probability below $0.05$.}
\label{fig:risk-batch}
\end{figure}

\begin{table}[!htb]
\centering\small
\caption{The 99.9\% objective quantile for unbatched SGD. Predictions and the matched linear reference are divided by the empirical quantile. All use the full objective.}
\label{tab:risk-main}
\begin{tabular}{ccccc}
\toprule
$\gamma$ & Empirical [95\% CI] & Linear / emp. & Jump / emp. & Gaussian / emp.\\
\midrule
0.05 & 3026.1 [2582.4, 3509.2] & 0.998 & 0.977 & 0.027\\
0.01 & 988.0 [835.7, 1145.7] & 0.995 & 1.005 & 0.016\\
0.002 & 380.5 [320.9, 448.9] & 1.000 & 0.892 & 0.008\\
\bottomrule
\end{tabular}
\end{table}

At the 99.9\% level, single-jump quantile predictions differ from the empirical estimates by $0.5\%$--$10.8\%$ and all lie within their pointwise intervals (Table~\ref{tab:risk-main}).
The matched linear quantiles differ by at most $0.5\%$.
The median Gaussian fits underestimate these quantiles by factors of about 37, 61, and 120.
Over the count-qualified thresholds in Figure~\ref{fig:risk-batch}(c), batches of four and sixteen give median probability ratios $0.487$ and $0.252$, close to $1/2$ and $1/4$.
The matched-$\gamma/b$ ratio has median $1.021$.
The comparison at matched $\gamma/b$, $(0.04,4)$ versus $(0.01,1)$, uses 2,000 and 2,303 gradient samples per trajectory, respectively.
See Appendix~\ref{app:merged-risk}.

\subsection{Random curvature and finite-step accuracy}

Random-design regression uses $Z=\sqrt d H^{1/2}U$, with $U$ uniform on the sphere, and gradient
$g(x;Z,\varepsilon)=Z(Z^\top x-\varepsilon)+6\boldsymbol a\sin^3(\boldsymbol a^\top x)/(1+(\boldsymbol a^\top x)^2)$, where $\varepsilon\sim t_{1.5}$.
The population objective is unchanged, but the sample Hessian at the target is $ZZ^\top$.
Alongside the population-Hessian benchmark, we propagate a shock through independent sample matrices $I_d-\gamma ZZ^\top$ to form the random-Hessian benchmark $\Lambda_\gamma^{\rm RH}$ (Appendix~\ref{app:random-hessian-benchmark}).
This isolates a source of finite-step error that the asymptotic population constant averages out.

\begin{table}[!htb]
\centering\small
\caption{Median projected-tail ratios over four directions and 17 thresholds in $[1,5]$. $R_{\rm tot}^{\rm RH}=\widehat p_X/\widehat\Lambda_\gamma^{\rm RH}$. The last column is the largest relative Monte Carlo standard error of $\widehat\Lambda_\gamma^{\rm RH}$ on this grid. Empirical sample sizes and intervals are in Appendix~\ref{app:multidim-experiments}.}
\label{tab:random-hessian-main}
\begin{tabular}{rrrrr}
\toprule
$d$ & $\gamma d$ & $R_{\rm tot}$ & $R_{\rm tot}^{\rm RH}$ & MC SE (\%)\\
\midrule
10 & 1 & 4.264 & 1.141 & 0.16\\
10 & 0.5 & 1.587 & 1.045 & 0.19\\
10 & 0.1 & 1.080 & 1.007 & 0.10\\
100 & 1 & 4.328 & 1.081 & 0.08\\
100 & 0.5 & 1.586 & 0.998 & 0.10\\
100 & 0.1 & 1.117 & 1.041 & 0.06\\
\bottomrule
\end{tabular}
\end{table}

At $\gamma d=1$, retaining random curvature reduces median ratios from about $4.3$ to $1.08$--$1.14$ (Table~\ref{tab:random-hessian-main}).
At $\gamma d=0.1$, both denominators are closer, and the random-Hessian ratios are $1.007$ and $1.041$.
A scalar random-curvature ablation in Appendix~\ref{app:random-curvature} gives the same separation between curvature amplification and nonlinear drift.

\section{Discussion}\label{sec:discussion}

Our results identify a regime in which local curvature and
the innovation law at the target determine rare errors of
nonlinear SGD.
The experiments show that finite-step calibration also depends
on nonlinear-flow effects at fixed thresholds and on random
sample curvature.
A stochastic fixed-tolerance extension under our assumptions,
together with quantitative control of these finite-step effects,
would strengthen the practical use of the tail approximation.
The present analysis assumes fixed dimension and centered
dissipativity around a single equilibrium.
Extending relative tail control to multiple basins and
dimension-dependent regimes is a natural direction for
further work.

\clearpage
\printbibliography
\appendix
\section{Proof of the \texorpdfstring{$L^p$}{Lp} Moment Bound}\label{app:moment-proof}

\begin{theorem}[Uniform moment bound]\label{thm:lp-moment}
Assume centered dissipativity and~\eqref{eq:global-lp-growth}.
For every $1<r\le p$, there are $C_r,c_r,\gamma_r>0$ such that
\begin{equation}\label{eq:lp-moment-bound}
 \E|X_n|^r\le(1-c_r\gamma)^n|X_0|^r+C_r\gamma^{r-1},
 \quad n\ge0,\quad 0<\gamma<\gamma_r.
\end{equation}
In particular,
\begin{equation}\label{eq:markov-tail-from-lp}
 \Pp(|X_n|>u)\le
 C_pu^{-p}\{\gamma^{p-1}+e^{-c_pn\gamma}|X_0|^p\}.
\end{equation}
\end{theorem}

This auxiliary bound requires neither regular variation nor local linearization.
We use the following elementary inequality for $1<r\le2$.

\begin{lemma}[Subquadratic increment inequality]\label{lem:subquad-increment}
For $1<r\le2$, there is $C_r<\infty$ such that, for deterministic $u,w\in\mathbb R^d$,
\begin{equation}\label{eq:subquad-increment}
|u+w|^r\le |u|^r+r|u|^{r-2}u^\top w+C_r|w|^r,
\end{equation}
with the middle term zero at $u=0$.
\end{lemma}
\begin{proof}
The gradient of $f(u)=|u|^r$ is $(r-1)$-H\"older continuous.
Integrating $\nabla f(u+sw)-\nabla f(u)$ against $w$ over $s\in[0,1]$
proves the inequality. Its conditional-expectation form also appears in
\citet[Appendix A.1]{Li2024SGD}.
\end{proof}

\begin{proof}[Proof of Theorem~\ref{thm:lp-moment}]
Fix $r\in(1,p]$. Since~\eqref{eq:global-lp-growth} holds in $L^p$, monotonicity of $L^s$ norms gives
\begin{equation}\label{eq:global-lr-growth}
    \norm{g(x,\varepsilon)}_r\le L_p(1+\abs{x-x^\star}),
    \qquad x\in\R^d.
\end{equation}
Applying Lemma~\ref{lem:subquad-increment} conditionally on $\F_{n-1}:=\sigma(\varepsilon_1,\ldots,\varepsilon_{n-1})$ with $u=x_{n-1}-x^\star$ and $w=-\gamma g(x_{n-1},\varepsilon_n)$ yields
\begin{align}
    \E[\abs{x_n-x^\star}^r\mid\F_{n-1}]
    &\le
    \abs{x_{n-1}-x^\star}^r
    -r\gamma \abs{x_{n-1}-x^\star}^{r-2}\ip{x_{n-1}-x^\star}{m(x_{n-1})}
    \notag\\
    &\quad+C_r\gamma^r\E[\abs{g(x_{n-1},\varepsilon_n)}^r\mid\F_{n-1}].
    \label{eq:moment-one-step-1}
\end{align}
By centered dissipativity~\eqref{eq:centered-dissipativity},
\begin{equation}\label{eq:dissipativity-r}
    \abs{x_{n-1}-x^\star}^{r-2}\ip{x_{n-1}-x^\star}{m(x_{n-1})}
    \ge
    \alpha \abs{x_{n-1}-x^\star}^r.
\end{equation}
By~\eqref{eq:global-lr-growth},
\begin{equation}\label{eq:noise-growth-r}
    \E[\abs{g(x_{n-1},\varepsilon_n)}^r\mid\F_{n-1}]
    \le C_r(1+\abs{x_{n-1}-x^\star}^r).
\end{equation}
Combining~\eqref{eq:moment-one-step-1}--\eqref{eq:noise-growth-r} gives
\begin{equation}\label{eq:moment-one-step-2}
    \E[\abs{x_n-x^\star}^r\mid\F_{n-1}]
    \le
    (1-r\alpha\gamma+C_r\gamma^r)\abs{x_{n-1}-x^\star}^r+C_r\gamma^r.
\end{equation}
Since $r>1$, for sufficiently small $\gamma$,
\[
    1-r\alpha\gamma+C_r\gamma^r\le 1-c_r\gamma
\]
for some $c_r>0$. Taking expectations and iterating yields
\begin{align}
    \E\abs{x_n-x^\star}^r
    &\le
    (1-c_r\gamma)^n\abs{x_0-x^\star}^r
    +C_r\gamma^r\sum_{j=0}^{n-1}(1-c_r\gamma)^j
    \notag\\
    &\le
    (1-c_r\gamma)^n\abs{x_0-x^\star}^r+C_r\gamma^{r-1}.
    \label{eq:moment-iteration-final}
\end{align}
This proves~\eqref{eq:lp-moment-bound}. The tail bound~\eqref{eq:markov-tail-from-lp} follows from~\eqref{eq:lp-moment-bound} with $r=p$ and Markov's inequality.
\end{proof}

\section{Proofs for the Linearized Recursion}\label{app:linear-proof}

\subsection{Stability of the vector linearization}

\begin{lemma}[Stability of the vector linearization]\label{lem:vector-linear-stability}
Assume Assumption~\ref{ass:centered-diss} and suppose that $m$ is differentiable at $x^\star$ with $H=\nabla m(x^\star)$. Then
\begin{equation}\label{eq:H-symmetric-positive}
        \frac{1}{2}(H+H^\top)\succeq \alpha I_d.
\end{equation}
Consequently, $H$ is positive stable. Moreover, there exists $\gamma_0>0$ such that, for all $0<\gamma<\gamma_0$, $B_\gamma:=I_d-\gamma H$ satisfies
\begin{equation}\label{eq:matrix-contraction-derived}
        \|B_\gamma^j\|_{\op}
        \le
        e^{-\alpha\gamma j/2},
        \qquad j\ge0.
\end{equation}
\end{lemma}

\begin{proof}[Proof of Lemma~\ref{lem:vector-linear-stability}]
Fix $u\in\mathbb R^d$ and $r>0$. Apply Assumption~\ref{ass:centered-diss} with $x=x^\star+ru$. Since $m(x^\star)=0$,
\[
        \left\langle u,
        \frac{m(x^\star+ru)-m(x^\star)}{r}
        \right\rangle
        \ge
        \alpha\|u\|^2.
\]
Letting $r\downarrow0$ gives $\langle u,Hu\rangle\ge\alpha\|u\|^2$. Since $\langle u,Hu\rangle=\langle u,(H+H^\top)u/2\rangle$, this proves~\eqref{eq:H-symmetric-positive}.

If $Hz=\lambda z$ for a possibly complex eigenvector $z\ne0$, then
\[
        \operatorname{Re}(\lambda)\|z\|^2
        =
        \operatorname{Re}(z^*Hz)
        =
        z^*\frac{H+H^\top}{2}z
        \ge
        \alpha\|z\|^2,
\]
so $H$ is positive stable. Moreover,
\[
        \|B_\gamma u\|^2
        =
        \|u-\gamma Hu\|^2
        \le
        (1-2\alpha\gamma+\gamma^2\|H\|_{\op}^2)\|u\|^2.
\]
For sufficiently small $\gamma$, the right-hand side is bounded by $(1-\alpha\gamma)\|u\|^2$. Hence
\[
        \|B_\gamma u\|\le (1-\alpha\gamma)^{1/2}\|u\|\le e^{-\alpha\gamma/2}\|u\|.
\]
Iterating proves~\eqref{eq:matrix-contraction-derived}.
\end{proof}

\subsection{Proof of Theorem~\ref{thm:vector-linear-tail}}

Throughout this subsection, write
\[
        Q_\gamma:=\Pp(|\eta_0|>r_\gamma),
        \qquad
        X_{\gamma,j}:=B_\gamma^j\eta_{n_\gamma-j},
        \qquad 0\le j<n_\gamma.
\]
Then $T_{\gamma,n_\gamma}=\sum_{j=0}^{n_\gamma-1}X_{\gamma,j}$. By Lemma~\ref{lem:vector-linear-stability}, there exists $c_B>0$ such that, for all sufficiently small $\gamma$,
\begin{equation}\label{eq:B-proof-contraction}
        \|B_\gamma^j\|_{\op}\le e^{-c_B\gamma j},
        \qquad j\ge0.
\end{equation}

\begin{lemma}[One-jump Riemann limit]\label{lem:one-jump-riemann}
Let $A\subset\mathbb R^d\setminus\{0\}$ be a Borel set satisfying~\eqref{eq:A-continuity}. Then
\begin{equation}\label{eq:one-jump-riemann}
        \gamma\sum_{j=0}^{n_\gamma-1}
        \frac{\Pp(r_\gamma^{-1}B_\gamma^j\eta_0\in A)}{Q_\gamma}
        \longrightarrow
        \nu_H(A).
\end{equation}
The same conclusion holds with $A$ replaced by any $\nu_H$-continuity set bounded away from the origin.
\end{lemma}

\begin{proof}
Fix $K<\infty$. Since $n_\gamma\gamma\to\infty$, the grid $\{j:\gamma j\le K\}$ is eventually contained in $\{0,\ldots,n_\gamma-1\}$. On $0\le \gamma j\le K$,
\begin{equation}\label{eq:euler-flow-uniform}
        \sup_{0\le j\le K/\gamma}\|B_\gamma^j-e^{-\gamma jH}\|_{\op}\to0.
\end{equation}
By product vague convergence and the continuous mapping $(s,z)\mapsto e^{-sH}z$,
\begin{equation}\label{eq:finite-K-exp-flow}
        \gamma\sum_{0\le j\le K/\gamma}
        \frac{\Pp(r_\gamma^{-1}e^{-\gamma jH}\eta_0\in A)}{Q_\gamma}
        \longrightarrow
        \int_0^K\mu_\eta\{z:e^{-sH}z\in A\}\,ds,
\end{equation}
where the boundary condition follows from $\nu_H(\partial A)=0$. The replacement of $e^{-\gamma jH}$ by $B_\gamma^j$ follows from~\eqref{eq:euler-flow-uniform} and the usual inner/outer enlargement argument; because $A$ is bounded away from the origin, the contribution of $|z|>R$ is controlled uniformly by Potter's bound \citep{Resnick_1987} and then sent to zero as $R\to\infty$. Hence
\begin{equation}\label{eq:finite-K-B-flow}
        \gamma\sum_{0\le j\le K/\gamma}
        \frac{\Pp(r_\gamma^{-1}B_\gamma^j\eta_0\in A)}{Q_\gamma}
        \longrightarrow
        \int_0^K\mu_\eta\{z:e^{-sH}z\in A\}\,ds.
\end{equation}
It remains to let $K\to\infty$. If $r_\gamma^{-1}B_\gamma^j\eta_0\in A$ and $a:=\operatorname{dist}(A,0)>0$, then $|\eta_0|\ge ar_\gamma/\|B_\gamma^j\|_{\op}$. Potter's bound \citep{Resnick_1987} and~\eqref{eq:B-proof-contraction} imply that, for every sufficiently small $\zeta>0$,
\[
        \frac{\Pp(r_\gamma^{-1}B_\gamma^j\eta_0\in A)}{Q_\gamma}
        \le
        C_A\|B_\gamma^j\|_{\op}^{t-\zeta}
        \le
        C_Ae^{-(t-\zeta)c_B\gamma j}.
\]
Therefore the discrete tail over $\gamma j>K$ is $O(e^{-cK})$, and the same bound controls the tail of the integral defining $\nu_H(A)$. Combining this with~\eqref{eq:finite-K-B-flow} proves~\eqref{eq:one-jump-riemann}.
\end{proof}

\begin{lemma}[Negligibility of many small jumps]\label{lem:small-jumps}
For every $\delta>0$,
\begin{equation}\label{eq:small-jumps-lemma}
        \lim_{\theta\downarrow0}\limsup_{\gamma\downarrow0}
        \frac{\gamma}{Q_\gamma}
        \Pp\left(
        \left\|
        \sum_{j=0}^{n_\gamma-1}
        X_{\gamma,j}\one_{\{\|X_{\gamma,j}\|\le\theta r_\gamma\}}
        \right\|>\delta r_\gamma
        \right)=0.
\end{equation}
\end{lemma}

\begin{proof}
Let $\widetilde X_{\gamma,j}^{(\theta)}:=X_{\gamma,j}\one_{\{\|X_{\gamma,j}\|\le\theta r_\gamma\}}$ and $\overline X_{\gamma,j}^{(\theta)}:=\widetilde X_{\gamma,j}^{(\theta)}-\E\widetilde X_{\gamma,j}^{(\theta)}$. Choose $\beta\in(t,2)$ and $0<\zeta<\beta-t$. A coordinate-wise von Bahr--Esseen inequality, equivalence of norms in finite dimension, and independence give
\[
        \E\left\|
        \sum_{j=0}^{n_\gamma-1}\overline X_{\gamma,j}^{(\theta)}
        \right\|^\beta
        \le
        C_\beta\sum_{j=0}^{n_\gamma-1}
        \E\|\widetilde X_{\gamma,j}^{(\theta)}\|^\beta.
\]
By Karamata's theorem~\citep{Resnick_1987} for the norm tail, Potter's bound~\citep{Resnick_1987}, and~\eqref{eq:B-proof-contraction},
\[
        \sum_{j=0}^{n_\gamma-1}
        \E\|\widetilde X_{\gamma,j}^{(\theta)}\|^\beta
        \le
        C r_\gamma^\beta Q_\gamma\theta^{\beta-t-\zeta}
        \sum_{j=0}^{n_\gamma-1}\|B_\gamma^j\|_{\op}^{t-\zeta}
        \le
        C r_\gamma^\beta Q_\gamma\theta^{\beta-t-\zeta}\gamma^{-1}.
\]
Markov's inequality yields
\begin{equation}\label{eq:centered-small-bound-vector}
        \Pp\left(
        \left\|
        \sum_{j=0}^{n_\gamma-1}\overline X_{\gamma,j}^{(\theta)}
        \right\|>\frac{\delta r_\gamma}{2}
        \right)
        \le
        C\delta^{-\beta}\theta^{\beta-t-\zeta}\frac{Q_\gamma}{\gamma}.
\end{equation}
Since $\E X_{\gamma,j}=B_\gamma^j\E\eta_0=0$,
\[
        \E\widetilde X_{\gamma,j}^{(\theta)}
        =-
        \E\left[X_{\gamma,j}\one_{\{\|X_{\gamma,j}\|>\theta r_\gamma\}}\right].
\]
Again by Karamata's theorem and Potter's bound~\citep{Resnick_1987},
\[
        \left\|
        \sum_{j=0}^{n_\gamma-1}\E\widetilde X_{\gamma,j}^{(\theta)}
        \right\|
        \le
        C r_\gamma\theta^{1-t-\zeta}\frac{Q_\gamma}{\gamma}
        =o(r_\gamma)
\]
for each fixed $\theta>0$, because $Q_\gamma/\gamma\to0$. Combining this with~\eqref{eq:centered-small-bound-vector} and then letting $\theta\downarrow0$ proves the lemma.
\end{proof}

\begin{proof}[Proof of Theorem~\ref{thm:vector-linear-tail}]
Let $A$ satisfy~\eqref{eq:A-continuity}, and set $a:=\operatorname{dist}(A,0)>0$. For $\varepsilon>0$, define $A^\varepsilon:=\{x:\operatorname{dist}(x,A)\le\varepsilon\}$ and $A_{-\varepsilon}:=\{x\in A:\operatorname{dist}(x,A^c)>\varepsilon\}$. Choose $\varepsilon<a/4$ so that $A^\varepsilon$ and $A_{-\varepsilon}$ are $\nu_H$-continuity sets.

For $\theta\in(0,a/4)$, let $N_{\gamma,\theta}:=\sum_{j=0}^{n_\gamma-1}\one_{\{\|X_{\gamma,j}\|>\theta r_\gamma\}}$. Lemma~\ref{lem:one-jump-riemann} applied to the norm-exceedance set gives $\E N_{\gamma,\theta}=O(Q_\gamma/\gamma)$; hence, by independence and $Q_\gamma/\gamma\to0$,
\begin{equation}\label{eq:double-big-vector}
        \Pp(N_{\gamma,\theta}\ge2)=o(Q_\gamma/\gamma).
\end{equation}
For the upper bound, if $r_\gamma^{-1}\sum_jX_{\gamma,j}\in A$, if at most one summand is larger than $\theta r_\gamma$, and if the sum of the small summands has norm at most $\varepsilon r_\gamma$, then the unique large summand, if present, must lie in $r_\gamma A^\varepsilon$. Thus
\begin{align*}
        \Pp\left(r_\gamma^{-1}\sum_jX_{\gamma,j}\in A\right)
        &\le
        \Pp\left(\bigcup_j\{r_\gamma^{-1}X_{\gamma,j}\in A^\varepsilon\}\right)
        +\Pp(N_{\gamma,\theta}\ge2)\\
        &\quad+
        \Pp\left(
        \left\|\sum_jX_{\gamma,j}\one_{\{\|X_{\gamma,j}\|\le\theta r_\gamma\}}\right\|>
        \varepsilon r_\gamma
        \right).
\end{align*}
Using the union bound, Lemma~\ref{lem:one-jump-riemann},~\eqref{eq:double-big-vector}, and Lemma~\ref{lem:small-jumps}, then letting $\theta\downarrow0$ and $\varepsilon\downarrow0$, gives
\[
        \limsup_{\gamma\downarrow0}
        \frac{\gamma}{Q_\gamma}
        \Pp\left(r_\gamma^{-1}\sum_jX_{\gamma,j}\in A\right)
        \le
        \nu_H(A).
\]
The lower bound is analogous. On the event that one summand lies in $r_\gamma A_{-\varepsilon}$, no second summand is larger than $\theta r_\gamma$, and the small-summand contribution has norm at most $\varepsilon r_\gamma$, the full sum lies in $r_\gamma A$. Bonferroni's inequality and the same estimates yield
\[
        \liminf_{\gamma\downarrow0}
        \frac{\gamma}{Q_\gamma}
        \Pp\left(r_\gamma^{-1}\sum_jX_{\gamma,j}\in A\right)
        \ge
        \nu_H(A).
\]
Therefore
\begin{equation}\label{eq:T-tail-limit-vector}
        \Pp(r_\gamma^{-1}T_{\gamma,n_\gamma}\in A)
        =
        \frac{Q_\gamma}{\gamma}\{\nu_H(A)+o(1)\}.
\end{equation}
Since $r_\gamma=u_\gamma/\gamma$, this proves the additive tail limit for $\gamma T_{\gamma,n_\gamma}/u_\gamma$, including sets with $\nu_H(A)=0$. Lemma~\ref{lem:one-jump-riemann} also gives
\[
        \sum_{j=0}^{n_\gamma-1}\Pp(r_\gamma^{-1}B_\gamma^j\eta_0\in A)
        =
        \frac{Q_\gamma}{\gamma}\{\nu_H(A)+o(1)\},
\]
which is used in the proof.

Finally, by~\eqref{eq:vector-linear-iteration}, $Y_{n_\gamma}=B_\gamma^{n_\gamma}Y_0+\gamma T_{\gamma,n_\gamma}$. The initialization condition in Assumption~\ref{ass:vector-threshold} implies $u_\gamma^{-1}B_\gamma^{n_\gamma}Y_0\to0$. A final inner/outer enlargement argument proves~\eqref{eq:vector-T-tail-asymp}. When $\nu_H(A)>0$, the two additive limits also give the stated relative single-jump equivalence.
\end{proof}

\subsection{Proof of Corollary~\ref{cor:linear-functional-tail}}

\begin{proof}
For $c>0$, let $A_{\ell,c}:=\{x:\ell^\top x>c\}$. This set is bounded away from the origin. Its boundary is the affine hyperplane $\partial A_{\ell,c}=\{x:\ell^\top x=c\}$. For each $s\ge0$,
\[
        (e^{-sH})^{-1}\partial A_{\ell,c}
        =
        \{z:\ell^\top e^{-sH}z=c\},
\]
which is an affine hyperplane not passing through the origin.

We claim that a homogeneous Radon measure on $\mathbb R^d\setminus\{0\}$ assigns zero mass to such a hyperplane. It is enough to check this on annuli $\{r\le\|z\|\le R\}$. If an affine hyperplane $L$ had positive $\mu_\eta$-mass on such an annulus, then for infinitely many distinct $\lambda_m$ sufficiently close to one, the sets $\lambda_m(L\cap\{r\le\|z\|\le R\})$ would be disjoint, have uniformly positive $\mu_\eta$-mass by homogeneity, and be contained in a slightly larger compact annulus away from zero. This contradicts Radon finiteness. Hence $\mu_\eta\{z:\ell^\top e^{-sH}z=c\}=0$ for every $s\ge0$, and therefore $\nu_H(\partial A_{\ell,c})=0$.

Taking $c=1$ and applying Theorem~\ref{thm:vector-linear-tail} to $A_\ell$ gives
\[
        \Pp(\ell^\top\gamma T_{\gamma,n_\gamma}>u_\gamma)
        \sim
        \frac{\Pp(|\eta_0|>u_\gamma/\gamma)}{\gamma}C_\ell.
\]
The deterministic shift $\ell^\top B_\gamma^{n_\gamma}Y_0$ is $o(u_\gamma)$ by Assumption~\ref{ass:vector-threshold}. Applying the same result to $(1\pm\varepsilon)u_\gamma$ and then letting $\varepsilon\downarrow0$ yields~\eqref{eq:linear-functional-Y-tail-asymp}.
\end{proof}

\subsection{Projected tails under one-dimensional regular variation}
For a general linear functional $\ell^\top Y_n$, the direction $(B_\gamma^j)^\top\ell$ changes with $j$, so the full multivariate tail measure is the natural object. In contrast, if $v$ is a left eigenvector of $H$, the projected recursion is genuinely scalar. This allows us to replace Assumption~\ref{ass:vector-rv} by a one-dimensional tail assumption on the projected innovation.

Specifically, let $v\in\R^d\setminus\{0\}$ be a fixed unit vector satisfying
\begin{equation}\label{eq:left-eigenvector}
    H^\top v=\lambda_v v,
    \qquad \lambda_v>0.
\end{equation}
Define the original-space linearized iterate $y_n:=x^\star+Y_n$. Set
\begin{equation}\label{eq:projected-linear-vars}
    y_n^{(v)}:=\ip{v}{y_n-x^\star}=\ip{v}{Y_n},
    \qquad
    \eta_n^{(v)}:=\ip{v}{\eta_n}=-\ip{v}{g(x^\star,\varepsilon_n)}.
\end{equation}
Then $\E\eta_n^{(v)}=0$, and by~\eqref{eq:left-eigenvector},

    $y_n^{(v)}=(1-\gamma\lambda_v)y_{n-1}^{(v)}+\gamma\eta_n^{(v)}.$

So the centered projected linear iterate is
\begin{equation}\label{eq:centered-linear-S}
    S_n^{(v)}:=y_n^{(v)}-\rho_{\gamma,v}^n y_0^{(v)}
    =\gamma\sum_{j=0}^{n-1}\rho_{\gamma,v}^{j}\eta_{n-j}^{(v)},
    \qquad y_0^{(v)}:=\ip{v}{x_0-x^\star},
\end{equation}
where $\rho_{\gamma,v}=1-\gamma\lambda_v$. Here we also define one-big-jump intensity
\begin{equation}\label{eq:one_big_jump_int}
    \Lambda_\gamma(u) \;:=\; \sum_{j=0}^{n_\gamma-1} P\!\left(\gamma\rho_{\gamma,v}^{\,j}\eta_0^{(v)} > u\right), \qquad u>0.
\end{equation}

\begin{assumption}[Projected tail regular variation]\label{ass:projected-tail-subq}
Let $v\in\mathbb R^d\setminus\{0\}$ be fixed and define
\begin{equation}\label{eq:projected-tail-functions}
    \bar F_v^+(u):=\Pp(\eta_0^{(v)}>u),
    \qquad
    \bar F_v^-(u):=\Pp(\eta_0^{(v)}<-u).
\end{equation}
There exist $t\in(1,2)$ and a slowly varying function $h_v$ such that
\begin{equation}\label{eq:projected-rv-tail}
    \bar F_v^+(u)\sim h_v(u)u^{-t},
    \qquad u\to\infty.
\end{equation}
Moreover, the left tail is dominated by the right-tail scale: there exist constants $C_-<\infty$ and $u_-<\infty$ such that
\begin{equation}\label{eq:left-tail-domination}
    \bar F_v^-(u)\le C_-\bar F_v^+(u),
    \qquad u\ge u_-.
\end{equation}
\end{assumption}

\begin{theorem}[Projected linear exact tail]\label{thm:linear-exact-subq}
Assume the left-eigenvector condition~\eqref{eq:left-eigenvector} and Assumption~\ref{ass:projected-tail-subq}. Let $\gamma\downarrow0$, $n=n_\gamma\to\infty$, and $u_\gamma\downarrow0$ satisfy
\begin{equation}\label{eq:projected-linear-regime}
    0<\gamma\lambda_v<1,
    \qquad
    n_\gamma\gamma\to\infty,
    \qquad
    r_\gamma:=\frac{u_\gamma}{\gamma}\to\infty,
    \qquad
    \frac{\bar F_v^+(r_\gamma)}{\gamma}\to0,
    \qquad
    \frac{\rho_{\gamma,v}^{n_\gamma}|y_0^{(v)}|}{u_\gamma}\to0.
\end{equation}
Then
\begin{equation}\label{eq:projected-linear-exact-tail}
    \Pp\left(y_{n_\gamma}^{(v)}>u_\gamma\right)
    \sim
    \Lambda_\gamma(u_\gamma)
    \sim
    \frac{1}{t\lambda_v}
    \frac{\bar F_v^+(u_\gamma/\gamma)}{\gamma}
    \sim
    \frac{1}{t\lambda_v}
    h_v(u_\gamma/\gamma)\gamma^{t-1}u_\gamma^{-t}.
\end{equation}
Equivalently,
\begin{equation}\label{eq:ratio_to_0}
    \left|\frac{\mathbb P(y_{n_\gamma}^{(v)}>u_\gamma)}{\Lambda_\gamma(u_\gamma)}-1\right|\to 0.
\end{equation}
\end{theorem}

\begin{remark}

\label{lem:one_jump_surrogate}
    Equation~\eqref{eq:ratio_to_0} identifies the leading tail probability as a sum of single-jump probabilities.

    The closed-form sharp constant $C_\ell$ is recovered by passing the Riemann-sum limit after replacing $\bar{F}_v^+$ with its regular-variation asymptotic: when an exact distribution function is available, one may evaluate the one-jump sum directly.

\end{remark}

\begin{remark}[Relation between the vector and projected assumptions]\label{rem:vector-vs-projected}
Under the stronger Assumption~\ref{ass:vector-rv}, if $A_v:=\{z:v^\top z>1\}$ satisfies $\mu_\eta(A_v)>0$, then the projected right-tail regular variation~\eqref{eq:projected-rv-tail} follows with $h_v(u)u^{-t}\sim \mu_\eta(A_v)\Pp(|\eta_0|>u)$. However, the projected theorem above only needs the one-dimensional right tail of $\eta_0^{(v)}$ and the comparable left-tail bound~\eqref{eq:left-tail-domination}.
\end{remark}

\begin{remark}[Power thresholds]\label{rem:linear-power-threshold}
For $u_\gamma=\gamma^a$ with $0<a<1$, the rare-event part of~\eqref{eq:projected-linear-regime} becomes
\begin{equation}\label{eq:linear-power-rare}
    h_v(\gamma^{a-1})\gamma^{t-1-at}\to0.
\end{equation}
The initialization part is implied, for example, by $n_\gamma\gamma/\log(1/\gamma)\to\infty$. Under these conditions,
\begin{equation}\label{eq:linear-power-tail}
    \Pp(y_{n_\gamma}^{(v)}>\gamma^a)
    \sim
    \frac{1}{t\lambda_v}h_v(\gamma^{a-1})\gamma^{t-1-at}.
\end{equation}
If $h_v(u)\to h_{v,0}\in(0,\infty)$, then~\eqref{eq:linear-power-rare} is equivalent to $a<(t-1)/t$.
\end{remark}

\subsection{Proof of Theorem~\ref{thm:linear-exact-subq}}

\begin{proof}
For notational simplicity, write $u=u_\gamma$, $n=n_\gamma$, $\rho=\rho_{\gamma,v}=1-\gamma\lambda_v$, $\lambda=\lambda_v$, and $\eta_j=\eta_j^{(v)}$. Put
\[
        Z_j:=\gamma\rho^j\eta_{n-j},\qquad 0\le j<n,
\]
so that $S_n^{(v)}=\sum_{j=0}^{n-1}Z_j$. Define the one-big-jump intensity
\[
        \Lambda_\gamma(u):=\sum_{j=0}^{n-1}\Pp(Z_j>u)
        =\sum_{j=0}^{n-1}\bar F_v^+\!\bigl(u/(\gamma\rho^j)\bigr).
\]

We first establish the asymptotic
\begin{equation}\label{eq:Lambda-projected-asymp}
        \Lambda_\gamma(u)\sim
        \frac{1}{t\lambda}\,\frac{\bar F_v^+(u/\gamma)}{\gamma}.
\end{equation}
Since the multipliers $\rho^{-j}$ range over an unbounded set as $j$ varies in $\{0,\dots,n-1\}$, neither Potter's bound \citep{Resnick_1987} alone nor the uniform convergence theorem alone yields~\eqref{eq:Lambda-projected-asymp} with the sharp constant; we therefore split the sum at a cutoff $K>0$ (assumed below to satisfy $K<n\gamma$, which holds eventually since $n\gamma\to\infty$). For $\gamma j\le K$, the values $\rho^{-j}$ lie in $[1,\rho^{-\lfloor K/\gamma\rfloor}]$, which converges to $[1,e^{\lambda K}]$ as $\gamma\downarrow 0$ and hence remains in a fixed compact subset of $[1,\infty)$ for small $\gamma$. Combining $u/\gamma\to\infty$ (from~\eqref{eq:projected-linear-regime}) with the uniform convergence theorem for slowly varying functions \citep
{Resnick_1987} applied to $h_v$,
\[
        \frac{\bar F_v^+(u/(\gamma\rho^j))}{\bar F_v^+(u/\gamma)}
        =\frac{h_v(u/(\gamma\rho^j))}{h_v(u/\gamma)}\,\rho^{tj}
        =(1+o(1))\,\rho^{tj},\qquad\gamma\downarrow 0,
\]
uniformly in $j$ with $\gamma j\le K$. Using $1-\rho^t\sim t\gamma\lambda$ and $\rho^{t\lfloor K/\gamma\rfloor}\to e^{-t\lambda K}$,
\begin{equation}\label{eq:Lambda-bounded-limit}
        \lim_{\gamma\downarrow 0}\frac{\gamma}{\bar F_v^+(u/\gamma)}\!\!
        \sum_{j:\,\gamma j\le K}\!\!\bar F_v^+(u/(\gamma\rho^j))
        =\frac{1-e^{-t\lambda K}}{t\lambda}.
\end{equation}
For $\gamma j>K$ we use Potter's bound \citep{Resnick_1987}: fix $\varepsilon\in(0,t-1)$, then there exist $A_\varepsilon<\infty$ and $X_\varepsilon<\infty$ such that for all $x\ge X_\varepsilon$ and $\mu\ge 1$,
\[
        \bar F_v^+(\mu x)/\bar F_v^+(x)\le A_\varepsilon\,\mu^{-t+\varepsilon}.
\]
For $\gamma$ small enough that $u/\gamma\ge X_\varepsilon$, taking $\mu=\rho^{-j}$ gives $\bar F_v^+(u/(\gamma\rho^j))\le A_\varepsilon\rho^{(t-\varepsilon)j}\bar F_v^+(u/\gamma)$ for every $j\ge 0$, and consequently
\begin{equation}\label{eq:Lambda-tail-limit}
        \limsup_{\gamma\downarrow 0}\frac{\gamma}{\bar F_v^+(u/\gamma)}\!\!
        \sum_{j:\,\gamma j>K}\!\!\bar F_v^+(u/(\gamma\rho^j))
        \le \frac{A_\varepsilon\,e^{-(t-\varepsilon)\lambda K}}{(t-\varepsilon)\lambda}.
\end{equation}
Combining~\eqref{eq:Lambda-bounded-limit} and~\eqref{eq:Lambda-tail-limit}, for every fixed $K>0$,
\[
        \frac{1-e^{-t\lambda K}}{t\lambda}
        \le\liminf_{\gamma\downarrow 0}\frac{\gamma\Lambda_\gamma(u)}{\bar F_v^+(u/\gamma)}
        \le\limsup_{\gamma\downarrow 0}\frac{\gamma\Lambda_\gamma(u)}{\bar F_v^+(u/\gamma)}
        \le\frac{1-e^{-t\lambda K}}{t\lambda}+\frac{A_\varepsilon e^{-(t-\varepsilon)\lambda K}}{(t-\varepsilon)\lambda}.
\]
Letting $K\to\infty$ with $\varepsilon$ fixed proves~\eqref{eq:Lambda-projected-asymp}. In particular, $\Lambda_\gamma(u)\to 0$ by~\eqref{eq:projected-linear-regime}. The same argument applied to $cu$ in place of $u$, together with $\bar F_v^+(cu/\gamma)/\bar F_v^+(u/\gamma)\to c^{-t}$ (which is direct since $cu/\gamma$ and $u/\gamma$ are comparable up to a fixed ratio), yields the shift stability
\begin{equation}\label{eq:Lambda-projected-shift}
        \Lambda_\gamma(cu)\sim c^{-t}\Lambda_\gamma(u),\qquad c>0\text{ fixed}.
\end{equation}

We next show that $\Pp(S_n^{(v)}>u)\sim\Lambda_\gamma(u)$. Fix $\delta\in(0,1/2)$ and $\theta\in(0,\delta/16)$, and define the truncated sum
\[
        W_\theta:=\sum_{j=0}^{n-1}Z_j\one_{\{|Z_j|\le\theta u\}}.
\]
We first establish moment bounds for $W_\theta$. Set $a_j:=\gamma\rho^j$ and $M_j:=\theta u/a_j$, so $Z_j=a_j\eta_{n-j}\stackrel{d}{=}a_j\eta_0$ and $\{|Z_j|>\theta u\}=\{|\eta_0|>M_j\}$. Since $M_j\ge\theta u/\gamma\to\infty$, Karamata's theorem \citep{Resnick_1987} applied to $|\eta_0|$ (whose tail is $\le(1+C_-)\bar F_v^+$ by~\eqref{eq:left-tail-domination}) gives
\[
        \E\bigl[|\eta_0|\mathbf 1_{\{|\eta_0|>M_j\}}\bigr]\le CM_j\bar F_v^+(M_j),
        \qquad
        \E\bigl[\eta_0^2\mathbf 1_{\{|\eta_0|\le M_j\}}\bigr]\le CM_j^2\bar F_v^+(M_j),
\]
where the second bound uses $t<2$. Summing the corresponding bounds for $Z_j$ and using $\sum_j\bar F_v^+(M_j)=\Lambda_\gamma(\theta u)$,
\[
        \sum_{j}\bigl|\E[Z_j\mathbf 1_{\{|Z_j|>\theta u\}}]\bigr|
        \le\sum_j a_j\E\bigl[|\eta_0|\mathbf 1_{\{|\eta_0|>M_j\}}\bigr]
        \le C\theta u\,\Lambda_\gamma(\theta u),
\]
\[
        \sum_{j}\E[Z_j^2\mathbf 1_{\{|Z_j|\le\theta u\}}]
        \le\sum_j a_j^2\E\bigl[\eta_0^2\mathbf 1_{\{|\eta_0|\le M_j\}}\bigr]
        \le C(\theta u)^2\Lambda_\gamma(\theta u).
\]
Since $\E Z_j=0$, $\E W_\theta=-\sum_j\E[Z_j\mathbf 1_{\{|Z_j|>\theta u\}}]$ and $\operatorname{Var}(W_\theta)\le\sum_j\E[Z_j^2\mathbf 1_{\{|Z_j|\le\theta u\}}]$. Combining these with the shift stability $\Lambda_\gamma(\theta u)\sim\theta^{-t}\Lambda_\gamma(u)$ from~\eqref{eq:Lambda-projected-shift}, we conclude
\begin{equation}\label{eq:scalar-trunc-moments}
        |\E W_\theta|\le C\theta^{1-t}u\Lambda_\gamma(u)=o(u),\qquad
        \operatorname{Var}(W_\theta)\le C\theta^{2-t}u^2\Lambda_\gamma(u),
\end{equation}
where the $o(u)$ uses that $\Lambda_\gamma(u)\to 0$ and $\theta$ is fixed.

For the upper bound, partition the summands according to $\mathcal L:=\{j:Z_j>\theta u\}$ and observe that on $\{S_n^{(v)}>u\}$ at least one of the following holds: some $Z_j$ exceeds $(1-\delta)u$; $|\mathcal L|\ge 2$; or $W_\theta>\delta u$ (since negative summands only reduce $S_n^{(v)}$). Therefore
\begin{equation}\label{eq:scalar-upper-start}
        \Pp(S_n^{(v)}>u)\le\Lambda_\gamma((1-\delta)u)+\Pp(|\mathcal L|\ge 2)+\Pp(W_\theta>\delta u).
\end{equation}
By independence and~\eqref{eq:Lambda-projected-shift}, $\Pp(|\mathcal L|\ge 2)\le\frac{1}{2}\Lambda_\gamma(\theta u)^2\sim\frac{1}{2}\theta^{-2t}\Lambda_\gamma(u)^2=o(\Lambda_\gamma(u))$. For the truncated sum, $|\E W_\theta|\le\delta u/2$ for $\gamma$ small enough by~\eqref{eq:scalar-trunc-moments}, and each centered summand $Z_j\one_{\{|Z_j|\le\theta u\}}-\E[\cdot]$ is bounded by $2\theta u$. Bennett--Prokhorov's inequality \citep{Prokhorov1959} therefore gives
\[
        \Pp(W_\theta>\delta u)
        \le\Pp(W_\theta-\E W_\theta>\delta u/2)
        \le\exp\!\left\{-\frac{\delta}{8\theta}\,\operatorname{arsinh}\!\left(\frac{\delta\,\theta^{t-1}}{2C\,\Lambda_\gamma(u)}\right)\right\}.
\]
Since $\operatorname{arsinh}(y)\sim\log(2y)$ as $y\to\infty$ and $\Lambda_\gamma(u)\to 0$, this bound equals $\Lambda_\gamma(u)^{\delta/(8\theta)+o(1)}$. The constraint $\theta<\delta/16$ forces the exponent to exceed $1$, so
\begin{equation}\label{eq:scalar-small-truncated}
        \Pp(W_\theta>\delta u)=o(\Lambda_\gamma(u)).
\end{equation}
Substituting~\eqref{eq:scalar-small-truncated} and~\eqref{eq:Lambda-projected-shift} into~\eqref{eq:scalar-upper-start},
\[
        \limsup_{\gamma\downarrow 0}\frac{\Pp(S_n^{(v)}>u)}{\Lambda_\gamma(u)}\le(1-\delta)^{-t}.
\]

For the lower bound, let $A_j:=\{Z_j>(1+\delta)u\}$ and $R_j:=\sum_{i\ne j}Z_i$, which is independent of $Z_j$. On $A_j\cap\{R_j>-\delta u\}$ the full sum $S_n^{(v)}=Z_j+R_j$ exceeds $u$, so by Bonferroni,
\[
        \Pp(S_n^{(v)}>u)\ge\sum_{j=0}^{n-1}\Pp(A_j)\Pp(R_j>-\delta u)-\sum_{i<j}\Pp(A_i)\Pp(A_j),
\]
and the double sum is bounded by $\frac{1}{2}\Lambda_\gamma((1+\delta)u)^2=o(\Lambda_\gamma(u))$ by independence and~\eqref{eq:Lambda-projected-shift}.

It remains to show $\sup_{0\le j<n}\Pp(R_j\le-\delta u)\to 0$. Decompose the event $\{R_j\le-\delta u\}$ according to whether some summand $Z_i$ with $i\ne j$ is itself less than $-\theta u$. The probability of at least one such large negative summand is, by the union bound and the comparable left-tail condition~\eqref{eq:left-tail-domination},
\[
        \Pp\!\left(\exists\,i\ne j:Z_i<-\theta u\right)
        \le\sum_{i\ne j}\Pp(Z_i<-\theta u)
        \le C_-\,\Lambda_\gamma(\theta u)
        =O\bigl(\theta^{-t}\Lambda_\gamma(u)\bigr)
\]
uniformly in $j$, which is $o(1)$ for fixed $\theta$ since $\Lambda_\gamma(u)\to 0$. On the complementary event $\{Z_i\ge-\theta u\text{ for all }i\ne j\}$, the sum $R_j$ coincides with the truncated sum
\[
        \widetilde W_\theta^{(j)}:=\sum_{i\ne j}Z_i\one_{\{|Z_i|\le\theta u\}}+\sum_{i\ne j}Z_i\one_{\{Z_i>\theta u\}}\,,
\]
whose first piece is a bounded centered sum (after centering) and whose second piece is non-negative. Hence on this event $R_j\ge\sum_{i\ne j}Z_i\one_{\{|Z_i|\le\theta u\}}=:W_\theta^{(j)}$, and
\[
        \Pp(R_j\le-\delta u)\le\Pp\!\left(\exists\,i\ne j:Z_i<-\theta u\right)+\Pp\bigl(W_\theta^{(j)}\le-\delta u\bigr).
\]
The stronger absolute-mean estimate
\[
        \sum_{i=0}^{n-1}
        \left|
        \E\!\left[Z_i\one_{\{|Z_i|\le\theta u\}}\right]
        \right|
        \le C\theta^{1-t}u\Lambda_\gamma(u)
\]
implies, uniformly in \(j\),
\[
        |\E W_\theta^{(j)}|
        \le C\theta^{1-t}u\Lambda_\gamma(u)=o(u).
\]
Moreover, dropping one index can only decrease the variance bound:
\[
        \operatorname{Var}(W_\theta^{(j)})
        \le C\theta^{2-t}u^2\Lambda_\gamma(u).
\]
Bennett--Prokhorov's inequality applied to $\E W_\theta^{(j)}-W_\theta^{(j)}$ then gives, exactly as in the upper-bound argument,
\[
        \Pp\bigl(W_\theta^{(j)}\le-\delta u\bigr)\le\Pp\bigl(\E W_\theta^{(j)}-W_\theta^{(j)}>\delta u/2\bigr)=o(\Lambda_\gamma(u))=o(1),
\]
uniformly in $j$, for $\theta<\delta/16$. Combining these two estimates,
\[
        \sup_{0\le j<n}\Pp(R_j\le-\delta u)\to 0,
\]
and therefore
\[
        \Pp(S_n^{(v)}>u)\ge(1-o(1))\Lambda_\gamma((1+\delta)u)-o(\Lambda_\gamma(u))\sim(1+\delta)^{-t}\Lambda_\gamma(u),
\]
giving $\liminf_{\gamma\downarrow 0}\Pp(S_n^{(v)}>u)/\Lambda_\gamma(u)\ge(1+\delta)^{-t}$. Letting $\delta\downarrow 0$ proves $\Pp(S_n^{(v)}>u)\sim\Lambda_\gamma(u)$,

and combination with~\eqref{eq:Lambda-projected-asymp} gives
\[
        \Pp(S_n^{(v)}>u)\sim\frac{1}{t\lambda}\,\frac{\bar F_v^+(u/\gamma)}{\gamma}.
\]

It remains to transfer this to $y_n^{(v)}=S_n^{(v)}+\rho^n y_0^{(v)}$. The last condition in~\eqref{eq:projected-linear-regime} gives $\rho^n|y_0^{(v)}|=o(u)$, so for any fixed $\epsilon\in(0,1)$, eventually
\[
        \Pp(S_n^{(v)}>(1+\epsilon)u)\le\Pp(y_n^{(v)}>u)\le\Pp(S_n^{(v)}>(1-\epsilon)u).
\]
Applying the asymptotic just established to $(1\pm\epsilon)u$ together with shift stability~\eqref{eq:Lambda-projected-shift}, and then letting $\epsilon\downarrow 0$, gives~\eqref{eq:projected-linear-exact-tail}; the second equivalence follows from $\bar F_v^+(u/\gamma)\sim h_v(u/\gamma)\,\gamma^t u^{-t}$.
\end{proof}

\paragraph{Relationship to the Wang--Rhee framework.}\rev{
For $u_\gamma=\gamma^a$, rescaling the linearized
recursion by $V_j=Y_j/u_\gamma$ gives
\[
V_{j+1}=(I-\gamma H)V_j+\gamma^{1-a}\eta_{j+1}.
\]
The noise exponent $1-a$ corresponds to the
general-scaling framework in Theorem~A.1 of
\citet{Rhee_2023}; its scaling restriction becomes
$1-a>1/t$, or $a<(t-1)/t$.
The two analyses therefore share this threshold-scaling restriction.

However, obtaining our last-iterate asymptotic result additionally
requires control as the observation horizon diverges,
including the contributions of old innovations and
the initial condition.
Their unclipped sample-path theorem also assumes
bounded drift and noise coefficients, which are not
assumed here. Furthermore, \citet{Rhee_2023} does not study general nonlinear SGD as we do in this paper. Their framework instead assumes that the stochastic noise can be factorized as the product of a function of the current iterate and an i.i.d. noise term. Our proof directly controls the matrix-weight limit,
remote-past contributions, and accumulated small jumps,
and removes initialization under the stated burn-in
condition. These estimates identify the
impulse-response constant $C_\ell$.
Theorem~\ref{thm:delta-q-bound} and Proposition~\ref{prop:general-threshold-transfer}
then provide a separate nonlinear transfer at the
relative rare-event scale.
Together, these estimates connect the growing observation horizon to a nonlinear final-iterate tail constant.}
\section{Central-scale moments and nonlinear comparison}\label{app:delta-proof}

\begin{lemma}[Consequences of the local and global conditions]\label{lem:local-global-consequences}
Under Assumptions~\ref{ass:centered-diss}--\ref{ass:local-quadratic}, write
$R(x)=H(x-x^\star)-m(x)$ and
$D(x,\varepsilon)=m(x)-g(x,\varepsilon)+g(x^\star,\varepsilon)$.
There are finite constants such that, for all $x$,
\begin{align}
 |m(x)|&\le C|x-x^\star|,\label{eq:m-linear-growth}\\
 \|D(x,\varepsilon)\|_p&\le C|x-x^\star|,\label{eq:D-global-bound}\\
 |R(x)|&\le C\{|x-x^\star|^2\one_{\{|x-x^\star|\le\kappa\}}
             +|x-x^\star|\one_{\{|x-x^\star|>\kappa\}}\}.\label{eq:R-global-bound}
\end{align}
Furthermore, $\|B_\gamma^j\|_{\op}\le e^{-c_Hj\gamma}$ for small $\gamma$.
\end{lemma}
\begin{proof}
Jensen's inequality, local increments, and global growth give the first two bounds as follows.
Inside the local ball, $m(x^\star)=0$ and~\eqref{eq:local-centered-lip} imply
$|m(x)|\le L_g|x-x^\star|$ and $\|D(x,\varepsilon)\|_p\le2L_g|x-x^\star|$.
Outside this ball, global growth absorbs the constant term into $C_\kappa|x-x^\star|$.
The same argument gives the linear bound on $R$ outside the ball; inside, use the local remainder.
The contraction follows from Lemma~\ref{lem:vector-linear-stability}.
\end{proof}

\begin{lemma}[Martingale moment inequality]\label{lem:vbe}
If $(d_k,\mathcal F_k)$ are $\mathbb R^d$-valued martingale differences with finite $r$th moments, $1<r\le2$, then
\begin{equation}\label{eq:vbe}
 \E\left|\sum_{k=1}^n d_k\right|^r\le C_r\sum_{k=1}^n\E|d_k|^r.
\end{equation}
\end{lemma}
\begin{proof}
Apply Lemma~\ref{lem:subquad-increment} to the partial sum through $k-1$ and $d_k$.
Conditioning on $\mathcal F_{k-1}$ removes the linear increment term.
Taking expectations and iterating proves the inequality.
This argument needs the martingale-difference property, rather than independence.
\end{proof}

\subsection{Proof of the central-scale moment bound}\label{app:central-proof}
\begin{proof}[Proof of Theorem~\ref{thm:central-moment}]
Translate $x^\star$ to zero and write $m_0(x)=m(x^\star+x)$.
The first two bounds of Lemma~\ref{lem:local-global-consequences} use only centered dissipativity and gradient regularity, so they apply without the quadratic remainder assumption.
For $h>0$, set $V_h(x)=(h^2+|x|^2)^{p/2}$.
Its Hessian satisfies
\[
 \nabla^2 V_h(x)=p(h^2+|x|^2)^{(p-2)/2}I_d
 +p(p-2)(h^2+|x|^2)^{(p-4)/2}xx^\top,
 \qquad 0\preceq\nabla^2V_h(x)\preceq p h^{p-2}I_d.
\]
The gradient is also $(p-1)$-H\"older with a constant independent of $h$.
Indeed, for $w\ne0$ put $e=w/|w|$. Along the segment $x+\theta w$, its Hessian norm is at most
$p|\langle x,e\rangle+\theta|w||^{p-2}$, whose integral over $0\le\theta\le1$ is at most
$p2^{2-p}|w|^{p-2}/(p-1)$.
Integrating the gradient difference gives the uniform Taylor bound
\begin{equation}\label{eq:regularized-taylor}
 0\le V_h(x+w)-V_h(x)-\langle\nabla V_h(x),w\rangle
 \le C_p\min\{h^{p-2}|w|^2,|w|^p\}.
\end{equation}
The scalar function $\phi_h(r)=\min(h^{p-2}r^2,r^p)$ is increasing and obeys $\phi_h(2r)\le4\phi_h(r)$.
Consequently,
\begin{equation}\label{eq:regularized-split}
 \phi_h(|u+v|)\le4\phi_h(|u|)+4\phi_h(|v|)
 \le4\phi_h(|u|)+4|v|^p.
\end{equation}

Write $Z=|\eta|$ and $\bar F(r)=\Pp(Z>r)$.
The generalized inverse in~\eqref{eq:central-scale} satisfies $\bar F(b_\gamma)\sim\gamma$; monotone regular variation with negative index ensures that relative jumps vanish, so atoms cause no difficulty.
Karamata's integral theorem \citep{BinghamGoldieTeugels1987}, with $1<p<t<2$, yields
\[
 \E[Z^2;Z\le b]\sim\frac{t}{2-t}b^2\bar F(b),\qquad
 \E[Z^p;Z>b]\sim\frac{t}{t-p}b^p\bar F(b).
\]
These formulas follow directly by integrating $\bar F$ on $(0,b)$ and $(b,\infty)$, respectively.
At $h=s_\gamma=\gamma b_\gamma$, the two branches of $\phi_h(\gamma Z)$ meet at $Z=b_\gamma$.
Thus
\begin{equation}\label{eq:central-noise-remainder}
 \E\phi_{s_\gamma}(\gamma Z)
 \sim\left(\frac{t}{2-t}+\frac{t}{t-p}\right)s_\gamma^p\bar F(b_\gamma)
 =O(\gamma s_\gamma^p).
\end{equation}

Given the past and $X_{n-1}=x$, set $F=x-\gamma m_0(x)$.
Dissipativity and $|m_0(x)|\le C|x|$ give $|F|^2\le(1-\alpha\gamma)|x|^2$ for small $\gamma$.
Concavity of $z^{p/2}$ and
$|x|^2(h^2+|x|^2)^{p/2-1}\ge V_h(x)-h^p$ imply
\begin{equation}\label{eq:central-deterministic-drift}
 V_h(F)\le(1-c_0\gamma)V_h(x)+c_0\gamma h^p,
 \qquad c_0=p\alpha/2.
\end{equation}
The next iterate is $F+\gamma(\eta_n+D(x_{n-1},\varepsilon_n))$.
Both noise terms have conditional mean zero, although they may be dependent.
Apply~\eqref{eq:regularized-taylor} at the past-measurable point $F$.
The linear Taylor term has conditional expectation zero.
Equations~\eqref{eq:regularized-split}--\eqref{eq:central-deterministic-drift} and~\eqref{eq:D-global-bound} give
\begin{align}
 \E[V_{s_\gamma}(X_n)\mid\mathcal F_{n-1}]
 &\le(1-c_0\gamma)V_{s_\gamma}(X_{n-1})
      +C\gamma s_\gamma^p+C\gamma^p|X_{n-1}|^p\notag\\
 &\le(1-c\gamma)V_{s_\gamma}(X_{n-1})+C\gamma s_\gamma^p.
 \label{eq:central-drift}
\end{align}
The last step absorbs $C\gamma^p|X_{n-1}|^p$ into the drift because $p>1$.
Iterating and using $V_{s_\gamma}(X_0)\le s_\gamma^p+|X_0|^p$ proves
\[
 \E|X_n|^p\le(1-c\gamma)^n|X_0|^p+Cs_\gamma^p.
\]
No finite second moment is used, and the initialization term vanishes exactly at $X_0=0$.
\end{proof}

\subsection{Comparison at lower moments}

\paragraph{Coupled error decomposition.}
Subtracting the recursions gives
\begin{equation}\label{eq:delta-decomposition-main}
 \Delta_n=\gamma\sum_{k=1}^nB_\gamma^{n-k}R(x_{k-1})
          +\gamma\sum_{k=1}^nB_\gamma^{n-k}D(x_{k-1},\varepsilon_k)
          =I_{1,n}+I_{2,n}.
\end{equation}
Here $R$ and $D$ are defined in Lemma~\ref{lem:local-global-consequences}.
The first term contains the population remainder; the second is a martingale convolution.

\begin{theorem}[Comparison without a radial tail assumption]\label{thm:delta-q-bound}
Under Assumptions~\ref{ass:centered-diss}--\ref{ass:local-quadratic}, for every $1\le q\le p$ there are $C_q,c_q,\gamma_q>0$ such that
\begin{equation}\label{eq:delta-q-bound}
 \E|\Delta_n|^q\le C_q\gamma^{p-1}
 +C_qe^{-c_qn\gamma}(|X_0|^p+|X_0|^q),
 \qquad n\ge1,\quad 0<\gamma<\gamma_q.
\end{equation}
\end{theorem}

\begin{proof}[Proof of Theorems~\ref{thm:central-delta} and~\ref{thm:delta-q-bound}]
Use $w_\gamma=s_\gamma^p$ under the assumptions of Theorem~\ref{thm:central-delta}, and $w_\gamma=\gamma^{p-1}$ for Theorem~\ref{thm:delta-q-bound}.
The corresponding moment theorem gives
$\E|X_k|^p\le Cw_\gamma+Ce^{-c_pk\gamma}|X_0|^p$.
For $q\in[1,p]$,~\eqref{eq:R-global-bound} implies $|R(x)|^q\le C_q|x-x^\star|^p$: locally $2q\ge p$, and outside the local ball $q\le p$.
Minkowski's inequality in~\eqref{eq:delta-decomposition-main} yields
\[
 \|I_{1,n}\|_q\le C_q\gamma\sum_{k=1}^n e^{-c_H(n-k)\gamma}
 \{w_\gamma^{1/q}+e^{-c_p(k-1)\gamma/q}|X_0|^{p/q}\}.
\]
The convolution of two decaying exponentials is bounded by $Ce^{-cn\gamma}$ for any $c$ strictly below their decay rates; a coincident-rate factor $n\gamma$ is absorbed by reducing $c$.
Therefore, including at $q=1$,
\begin{equation}\label{eq:I1-bound-final-proof}
 \E|I_{1,n}|^q\le C_qw_\gamma+C_qe^{-c_qn\gamma}|X_0|^p.
\end{equation}
Let $\mathcal F_k=\sigma(\varepsilon_1,\ldots,\varepsilon_k)$.
Fresh innovations and unbiasedness imply
\[
 \E[D(x_{k-1},\varepsilon_k)\mid\mathcal F_{k-1}]=0,\qquad
 \E[|D(x_{k-1},\varepsilon_k)|^p\mid\mathcal F_{k-1}]\le C|X_{k-1}|^p.
\]
The deterministic weights in $I_{2,n}$ preserve the martingale-difference property.
Apply Lemma~\ref{lem:vbe} at order $p$:
\begin{align}
 \E|I_{2,n}|^p
 &\le C\gamma^p\sum_{k=1}^n e^{-pc_H(n-k)\gamma}
       \{w_\gamma+e^{-c_p(k-1)\gamma}|X_0|^p\}\notag\\
 &\le C\gamma^{p-1}w_\gamma+C\gamma^{p-1}e^{-cn\gamma}|X_0|^p.
 \label{eq:I2-p-first}
\end{align}
Lyapunov's inequality now gives
\begin{equation}\label{eq:I2-final-proof}
 \E|I_{2,n}|^q\le C_q\gamma^{q(p-1)/p}w_\gamma^{q/p}
 +C_q\gamma^{q(p-1)/p}e^{-c_qn\gamma}|X_0|^q.
\end{equation}
For $w_\gamma=\gamma^{p-1}$, the first term is $C_q\gamma^{2q(p-1)/p}=o(w_\gamma)$ because $2q>p$.
For $w_\gamma=s_\gamma^p$, its ratio to $w_\gamma$ is regularly varying at zero with exponent
\[
 e(q)=\frac{q(p-1)}p+(q-p)(1-1/t)
 \ge e(1)=(p-1)(1/p+1/t-1)>0.
\]
Thus it is again $o(w_\gamma)$, including all slowly varying factors.
Combining~\eqref{eq:I1-bound-final-proof} and~\eqref{eq:I2-final-proof}, with $\gamma\le1$, proves both claims.
\end{proof}

\section{Transfer of vector and objective-value tails}\label{app:transfer-proof}

\begin{proposition}[General-threshold transfer]\label{prop:general-threshold-transfer}
Suppose the linear vector tail theorem holds at $u_\gamma$ and at fixed multiples of $u_\gamma$.
Suppose also that, for some $q\ge1$,
\[
 \E|\Delta_{n_\gamma}|^q=o\{u_\gamma^q L_\gamma(u_\gamma)\}.
\]
Then~\eqref{eq:nonlinear-vector-tail} holds for every $A$ satisfying~\eqref{eq:A-continuity}.
\end{proposition}
\begin{proof}
For $\delta>0$, define $A^\delta=\{z:\operatorname{dist}(z,A)\le\delta\}$ and
$A_{-\delta}=\{z\in A:\operatorname{dist}(z,A^c)>\delta\}$.
Take $\delta<\operatorname{dist}(A,0)/2$.
The inclusions induced by $X=Y+\Delta$ give
\[
 \Pp(Y/u\in A_{-\delta})-\Pp(|\Delta|>\delta u)
 \le\Pp(X/u\in A)
 \le\Pp(Y/u\in A^\delta)+\Pp(|\Delta|>\delta u).
\]
Markov's inequality makes the error $o(L_\gamma(u))$ for each fixed $\delta$.
Choose $\delta$ decreasing to zero through values for which the enlarged and contracted sets are continuity sets.
Such values exist because the distance level sets are disjoint and the measure is finite away from zero.
Divide by $L_\gamma(u)$, take upper and lower limits, and then let $\delta$ decrease to zero.
Continuity from above and below, together with $\nu_H(\partial A)=0$, gives the result, including the case $\nu_H(A)=0$.
\end{proof}

\begin{proof}[Proof of Theorem~\ref{thm:nonlinear-exact-tail}]
Let $u_\gamma=c\gamma^a$ with the exponents in~\eqref{eq:nonlinear-parameter-window}.
Then $r_\gamma\to\infty$ and $L_\gamma(u_\gamma)\to0$ by regular variation and Potter bounds.
First suppose $p<t$.
If $X_0=0$, both the linear initial shift and the exponential terms in~\eqref{eq:central-delta} vanish.
Otherwise, the superlogarithmic time condition makes the initial shift and each exponential error smaller than every positive power of $\gamma$.
Hence Assumption~\ref{ass:vector-threshold} holds.
Regular variation gives
\[
\begin{aligned}
\frac{s_\gamma^p}{u_\gamma L_\gamma(u_\gamma)}
&= \gamma^{p(1-1/t)-(t-1)+a(t-1)+o(1)} \\
&= \gamma^{(t-1)(a+p/t-1)+o(1)}
\longrightarrow 0,
\end{aligned}
\]
where the exponent is positive because $a>(t-p)/t$.
Theorem~\ref{thm:central-delta} with $q=1$, together with the initialization
bounds above, therefore yields
\[
\mathbb{E}|\Delta_{n_\gamma}|
=o\{u_\gamma L_\gamma(u_\gamma)\}.
\]
Apply Theorem~\ref{thm:vector-linear-tail} and Proposition~\ref{prop:general-threshold-transfer}.
If the full gradient conditions hold at the finite endpoint $p=t$, choose $\max\{1,t(1-a)\}<p'<t$ so that the strict threshold window still holds.
The gradient conditions descend to $p'$ by monotonicity of moment norms, giving the same tail conclusion.
\end{proof}

\begin{proof}[Proof of Corollary~\ref{cor:nonlinear-left-eigen-tail}]
For $H^\top v=\lambda_vv$, homogeneity gives
$C_v=\mu_\eta(A_v)\int_0^\infty e^{-t\lambda_vs}\dd s=\mu_\eta(A_v)/(t\lambda_v)$.
Regular variation also gives
$\Pp(v^\top\eta>r)\sim\mu_\eta(A_v)\Pp(|\eta|>r)$.
This proves the assertion.
\end{proof}

\begin{proposition}[Nonlinear eigenprojection with only projected regular variation]\label{prop:projected-only}
Retain Assumptions~\ref{ass:centered-diss}--\ref{ass:local-quadratic} and suppose a unit eigenvector $v$ satisfies Assumption~\ref{ass:projected-tail-subq} with index $t$.
If
\[
 2-1/t<p\le t,\qquad \frac{t-p}{t-1}<a<\frac{t-1}{t},
 \qquad u_\gamma=c\gamma^a,\quad c>0,
\]
and the initialization conditions of Theorem~\ref{thm:nonlinear-exact-tail} hold, then~\eqref{eq:nonlinear-eigen-tail} holds.
\end{proposition}
\begin{proof}
Apply Theorem~\ref{thm:linear-exact-subq} to $v^\top Y_n$.
The $q=1$ bound in Theorem~\ref{thm:delta-q-bound}, divided by
$u_\gamma\Pp(v^\top\eta>u_\gamma/\gamma)/\gamma$, has power exponent $p-t+a(t-1)>0$.
The scalar sandwich at thresholds $(1\pm\delta)u_\gamma$ proves the claim.
The lower and upper bounds on $a$ define a nonempty interval exactly when $p>2-1/t$.
\end{proof}
This result needs no radial tail control for other coordinates.
It therefore uses the moment-growth comparison, whereas the wider main-text window uses the radial regular variation in Assumption~\ref{ass:vector-rv}.

\begin{proof}[Proof of Corollary~\ref{cor:risk-tail}]
Positive definiteness of $H$ places $A_M$ between two norm-exceedance sets, so $0<C_M<\infty$.
The quadratic homogeneity implies
$\nu_H\{z:z^\top Hz/2>s\}=s^{-t/2}C_M$ for $s>0$.
Consequently the boundary of $A_M$ has zero mass.
Locally,
\[
 M(x^\star+x)-M(x^\star)=\tfrac12x^\top Hx+O(|x|^3).
\]
For a fixed $R$, the remainder divided by $u_\gamma^2$ vanishes uniformly on $|x|\le Ru_\gamma$.
Sandwich the objective exceedance event on this ball between quadratic exceedances at $1\pm\delta$.
Theorem~\ref{thm:nonlinear-exact-tail} applies to the truncated continuity sets.
The omitted region satisfies
\[
 \lim_{\gamma\downarrow0}
 \frac{\Pp(|X_{n_\gamma}|>Ru_\gamma)}{L_\gamma(u_\gamma)}
 =\nu_H\{|z|>R\}=R^{-t}\nu_H\{|z|>1\}.
\]
Let $\delta$ decrease to zero and $R$ increase to infinity.
This proves~\eqref{eq:risk-tail} without asserting uniform Taylor approximation on an unbounded set.
\end{proof}

\section{Extensions and parameter interpretation}\label{app:extensions}

\subsection{H\"older local remainders}
\begin{theorem}[Superlinear local remainder]\label{thm:holder-extension}
Replace Assumption~\ref{ass:local-quadratic} by differentiability at $x^\star$ and
$|m(x)-H(x-x^\star)|\le C|x-x^\star|^\beta$ for $|x-x^\star|\le\kappa$,
where $1<\beta\le2$.
Retain the other assumptions.
Let $q_*=\max(1,p/\beta)$.
For $1<p<t$,~\eqref{eq:central-delta} holds for $q\in[q_*,p]$.
The nonlinear vector, projection, and objective-value tail results hold under
\begin{equation}
    \frac{(t-1)(t-p)}{t(t-q_*)}<a<\frac{t-1}{t}
\end{equation}
and the same initialization conditions.
\end{theorem}
\begin{proof}
For $q\ge q_*$ the local relation $\beta q\ge p$ gives $|R(x)|^q\le C|x-x^\star|^p$.
Outside the ball, $q\le p$ and the global linear bound give the same inequality.
Thus the deterministic convolution calculation is unchanged.
The martingale $p$th-moment calculation does not use $\beta$.
This proves the error bound, including $q=q_*$ and $q=p$.
Use $q=q_*$ in Markov's inequality.
The relative-error exponent becomes $p(1-1/t)-(t-1)+a(t-q_*)$, which is positive in the stated window.
The interval is nonempty exactly when $p>q_*$, true for every $p>1$ and $\beta>1$.
Proposition~\ref{prop:general-threshold-transfer} proves the vector conclusion.
The projection is a half-space instance of the vector conclusion.
For a finite $p=t$ moment, lower the moment order slightly as in the proof of Theorem~\ref{thm:nonlinear-exact-tail}; continuity of $q_*=\max(1,p/\beta)$ preserves every strict window.
For the objective tail, the local Taylor remainder is now $O(|x|^{1+\beta})=o(|x|^2)$, sufficient for the same compact-truncation proof.
\end{proof}

\paragraph{A model outside the quadratic-remainder class.}
Let $1<\beta<2$, $c>0$, and
\[
 m(x)=x+c\frac{x|x|^{\beta-1}}{1+|x|^{\beta-1}},\qquad
 g(x,\varepsilon)=m(x)-\varepsilon,\qquad \varepsilon\sim t_t.
\]
Here $xm(x)\ge x^2$, $m'(0)=1$, and $1\le m'(x)\le1+c$.
The stochastic gradient satisfies the growth and increment assumptions for every $1<p<t$.
Its local remainder is
$-c\,\operatorname{sign}(x)|x|^\beta(1+o(1))$, which is not $O(|x|^2)$.
This example establishes a strict regularity extension; its objective is strongly convex.
Taking $M(x)=\int_0^x m(s)\dd s$ and $G(x,\varepsilon)=M(x)-\varepsilon x$ gives an integrable sample loss.

\paragraph{Available moment and threshold exponents.}
For $\beta=2$, $t=1.5$ and $p=1.35$ allow $0.1<a<1/3$.
For $t=1.5$ and $p=1.45$, $\beta=1.1$ permits $11/120<a<1/3$, while $\beta=1.2$ permits $2/35<a<1/3$.
The examples $a=.30$ and $a=.25$ have positive remainder exponents $5/132$ and $9/160$, respectively.
These sufficient exponents establish asymptotic validity, not fast finite-step convergence.
If the full gradient assumptions hold for every $p<t$, taking $p$ close enough to $t$ proves validity for every fixed $0<a<(t-1)/t$ and every $\beta>1$.
Merely knowing that the target innovation has these moments is insufficient for this last statement.

\subsection{Fixed mini-batches}
\begin{corollary}[Fixed-batch scaling]\label{cor:batch}
Let $b$ be a fixed positive integer and use the average of $b$ independent fresh stochastic gradients at each step.
Under the corresponding conditions of Theorem~\ref{thm:nonlinear-exact-tail} or~\ref{thm:holder-extension},
\[
 \Pp(X^{(b)}_{n_\gamma}/u_\gamma\in A)
 =b^{1-t}L_\gamma(u_\gamma)\{\nu_H(A)+o(1)\},
\]
where $L_\gamma$ and $\mu_\eta$ use the single-sample innovation normalization.
\end{corollary}
\begin{proof}
Finite sums of independent regularly varying vectors have the sum of their tail measures.
Rescaling the sum by $b$ therefore gives tail measure $b\cdot b^{-t}\mu_\eta$.
The batch average preserves the population gradient, centered dissipativity, and local remainder.
Minkowski's inequality preserves the required growth and increment bounds.
Apply the nonlinear tail theorem with the batch innovation and rewrite its radial normalization in terms of the original innovation.
\end{proof}
At a common admissible threshold sequence, $(\gamma,b)$ and $(\gamma/b,1)$ have the same leading tail.
Indeed, regular variation gives
$L_{\gamma/b}(u_\gamma)\sim b^{1-t}L_\gamma(u_\gamma)$.
For a shared sample budget $S_\gamma$, their effective times are both asymptotic to $S_\gamma\gamma/b$.
This comparison requires that the budget grows enough to satisfy the applicable initialization conditions.
Equivalently, $(b\gamma,b)$ matches $(\gamma,1)$.
This fixed-$b$ result does not interchange the limits in $\gamma$ and $b$.

\subsection{Extreme objective quantiles}\label{app:E3}
Under Theorem~\ref{thm:holder-extension} and its horizon conditions, suppose $\Pp(|\eta|>r)\sim c_\eta r^{-t}$ and let $\delta_\gamma=\gamma^\zeta$ with
\[
 0<\zeta<\frac{(t-1)(p-q_*)}{t-q_*}.
\]
The excess-objective quantile at probability $1-\delta_\gamma$ satisfies
\[
 Q_{1-\delta_\gamma}
 \sim\left(\frac{C_Mc_\eta\gamma^{t-1}}{\delta_\gamma}\right)^{2/t}.
\]
To prove this, apply the risk-tail result at fixed multiples of the proposed quantile and use monotonicity.
The implied exponent $a=(t-1-\zeta)/t$ lies in the valid window.
For a general slowly varying factor, the corresponding asymptotic inverse must be retained.
The result concerns extreme quantiles whose tail probabilities vanish.
\subsection{Conditional sizes and directions}\label{app:conditional-law}
Put $D_1=\{z:|z|>1\}$ and $c_H=\nu_H(D_1)>0$.
Homogeneity gives $\nu_H\{|z|>r\}=r^{-t}c_H$ and $\nu_H\{|z|=r\}=0$ for every $r>0$.
Dividing Theorem~\ref{thm:nonlinear-exact-tail} for $A\cap D_1$ by its conclusion for $D_1$ gives
\[
 \mathcal L(X_{n_\gamma}/u_\gamma\mid |X_{n_\gamma}|>u_\gamma)
 \ \Longrightarrow\ \frac{\nu_H(\,\cdot\cap D_1)}{c_H}.
\]
Indeed, this holds for every continuity set of the probability measure on the right, which characterizes weak convergence.
The continuous map $z\mapsto z/|z|$ on $D_1$ gives the conditional direction law.
In particular, for a Borel set $C\subset S^{d-1}$ with boundary of zero angular mass, its probability is
$\nu_H\{z:|z|>1,z/|z|\in C\}/c_H$.
These are conditional directional probabilities; the constants of separate eigenvector projections are marginal quantities.

\subsection{Affine population gradients at fixed tolerances}\label{app:affine-fixed}
\begin{proof}[Proof of Corollary~\ref{cor:affine-fixed}]
Because $R\equiv0$, the decomposition~\eqref{eq:delta-decomposition-main} has only its martingale term.
Theorem~\ref{thm:central-moment} and~\eqref{eq:I2-p-first} imply
\begin{equation}\label{eq:affine-error}
 \E|\Delta_n|^p\le C\gamma^{p-1}s_\gamma^p
       +C\gamma^{p-1}e^{-cn\gamma}|X_0|^p.
\end{equation}
For a fixed $u>0$, $L_\gamma(u)$ is regularly varying at zero with exponent $t-1$.
The ratio of the first term in~\eqref{eq:affine-error} to $u^p L_\gamma(u)$ has exponent
\[
 (p-1)+p(1-1/t)-(t-1)=p(2-1/t)-t>0.
\]
The strict inequality absorbs all slowly varying factors.
The initialization term vanishes at $X_0=0$; otherwise the stated superlogarithmic horizon makes it negligible.
For $u_\gamma\equiv u$, Assumption~\ref{ass:vector-threshold} holds: $u/\gamma\to\infty$, $L_\gamma(u)\to0$, and the initial linear shift is negligible.
The linear proof uses these conditions on $u/\gamma$ and $L_\gamma(u)$, with no requirement that $u$ itself decrease.
Apply Theorem~\ref{thm:vector-linear-tail}, Markov's inequality at order $p$, and Proposition~\ref{prop:general-threshold-transfer}.
If the complete gradient conditions hold at $p=t$, choose any $p'\in(t^2/(2t-1),t)$ and apply the preceding argument at order $p'$.
\end{proof}

For least squares with a bounded random design $Z$, $H=\E ZZ^\top\succ0$, and independent centered regularly varying label noise $\epsilon$, the gradient
$g(x;Z,\epsilon)=Z\{Z^\top(x-x^\star)-\epsilon\}$ has mean $H(x-x^\star)$.
The sample loss contrast $\tfrac12[(Z^\top(x-x^\star)-\epsilon)^2-\epsilon^2]$, with subtraction performed before expectation, has finite population objective $\tfrac12(x-x^\star)^\top H(x-x^\star)$ even when $\E\epsilon^2=\infty$.
Boundedness of $Z$ gives the complete gradient conditions at every $p<t$.
Conditioning on $Z$ and applying regular variation with a dominated bound gives multivariate regular variation of $Z\epsilon$, provided its tail measure is nonzero.
Thus the corollary includes random sample Hessians.
The cutoff on $p$ is a sufficient condition from the comparison bound; at $t=1.5$ it is $p>9/8$.

\subsection{Finite observation horizons}\label{app:finite-horizon}
Let $X_0=0$ and $n_\gamma\gamma\to T\in(0,\infty)$.
Under the remaining assumptions and threshold window of Theorem~\ref{thm:nonlinear-exact-tail}, define
\[
 \nu_{H,T}(A)=\int_0^T\mu_\eta\{z:e^{-sH}z\in A\}\,\dd s.
\]
For every Borel $A$ bounded away from zero with $\nu_{H,T}(\partial A)=0$,
\begin{equation}\label{eq:finite-horizon-tail}
 \Pp(X_{n_\gamma}/u_\gamma\in A)
 =L_\gamma(u_\gamma)\{\nu_{H,T}(A)+o(1)\}.
\end{equation}
To prove this, use the same triangular array as in Appendix~\ref{app:linear-proof}.
Its one-jump Riemann sum now ends at $T$, while the truncated small-jump bounds remain valid because all finite-horizon sums are bounded by the corresponding infinite sums of contracted weights.
The time-zero boundary contributes a term of order $\gamma$ after normalization and hence vanishes.
The uniform-in-$n$ bound of Theorem~\ref{thm:central-delta} has no initial term, so the identical set sandwich transfers the limit to $X$.
The case $p=t$ again follows by lowering the moment slightly.
For a unit eigenvector $v$, radial homogeneity gives
\[
 C_v(T)=\frac{\mu_\eta\{z:v^\top z>1\}}{t\lambda_v}
              (1-e^{-t\lambda_vT})=C_v(1-e^{-t\lambda_vT}).
\]

\subsection{Localization of the nonlinear flow}\label{app:flow-localization}
Assume in addition that $m$ is locally Lipschitz on $\mathbb R^d$.
Let $\varphi_s(z)$ solve $\dot\varphi_s=-m(x^\star+\varphi_s)$ with $\varphi_0=z$.
Centered dissipativity gives $|\varphi_s(z)|\le e^{-\alpha s}|z|$, so the unique solution exists for all $s\ge0$.
For $u>0$, define
\begin{equation}\label{eq:flow-constant}
 C_A^{\rm flow}(u)=\int_0^\infty\mu_\eta\{z:\varphi_s(z)\in uA\}\,\dd s.
\end{equation}
\begin{proposition}[Deterministic localization]\label{prop:flow-localization}
Under Assumptions~\ref{ass:centered-diss}, \ref{ass:local-quadratic}, and~\ref{ass:vector-rv}, and the local Lipschitz condition above,
\[
 u^t C_A^{\rm flow}(u)\longrightarrow\nu_H(A),\qquad u\downarrow0,
\]
for every Borel set $A$ with $\operatorname{dist}(A,0)>0$ and $\nu_H(\partial A)=0$.
The statement also holds for the superlinear remainder in Theorem~\ref{thm:holder-extension}.
\end{proposition}
\begin{proof}
Write $\psi_u(s,z)=\varphi_s(uz)/u$.
Homogeneity changes $u^t C_A^{\rm flow}(u)$ into
$\int_0^\infty\int\one_{\{\psi_u(s,z)\in A\}}\,\mu_\eta(\dd z)\,\dd s$.
For fixed $(s,z)$, contraction keeps the path inside the local ball for small $u$.
The rescaled equation has the form
$\dot\psi_u=-H\psi_u+\rho_u(\psi_u)$ with $|\rho_u(\psi_u)|\le Cu^{\beta-1}|z|^\beta$, where $\beta=2$ in the quadratic case.
Variation of constants therefore gives
$|\psi_u(s,z)-e^{-sH}z|\le Cu^{\beta-1}|z|^\beta s\to0$.
The exceptional points whose limits lie on $\partial A$ have product-measure mass $\nu_H(\partial A)=0$.
If $d_A=\operatorname{dist}(A,0)>0$, the indicators are dominated by
$\one_{\{|z|\ge d_Ae^{\alpha s}\}}$.
The integral of this envelope is $d_A^{-t}/(\alpha t)<\infty$.
Dominated convergence proves the assertion.
\end{proof}

For a scalar restoring drift and a positive-tail set $A=(1,\infty)$, changing the order of integration over shock size and time gives
\[
 C_{(1,\infty)}^{\rm flow}(u)
 =\mu_\eta(1,\infty)\int_u^\infty\frac{r^{-t}}{m(x^\star+r)}\,\dd r.
\]
The normalized factor is thus
$tH u^t\int_u^\infty r^{-t}/m(x^\star+r)\,\dd r$.
Figure~\ref{fig:threshold-sweep}(b) compares this factor with the saved paired endpoints; Appendix~\ref{app:scalar-experiments} gives the scalar studies.
Establishing the corresponding fixed-tolerance stochastic tail formula for general nonlinear drifts requires a separate theorem.

For the finite-step check in Figure~\ref{fig:threshold-sweep}, put $\psi_\gamma(x)=x-\gamma m(x)$ and $z_j=\psi_\gamma^{-j}(u)$.
The displayed stepsizes make $\psi_\gamma$ increasing and preserve the positive half-line.
The diagnostic is the ratio of $\sum_{j<n}\Pp\{\gamma\sigma\varepsilon>z_j\}$ to $\sum_{j<n}\Pp\{\gamma\sigma\varepsilon>u(1-\gamma)^{-j}\}$, using the exact Student survival function and all $n$ shock ages.
Across Figure~\ref{fig:threshold-sweep}'s grid it differs from the continuous factor by at most $0.00063$; at the six prespecified points its discrepancy from $R_{\rm tr}$ is at most $0.00813$.
The supplement supplies both calculations and independent quadrature and root-solver checks.

\subsection{Different innovation and state-dependent moment indices}\label{app:limited-moment-example}
Consider $g(x;\eta,\zeta)=m(x)+x\zeta-\eta$, where
$m(x)=x+4\sin^3(x)/(1+x^2)$ and independent centered Student variables $\eta$ and $\zeta$ have tail indices $1.5$ and $1.25$, respectively.
The innovation at $x^\star=0$ is $\eta$, whereas the full gradient has moments only below $1.25$ away from zero.
The regularity assumptions hold at $p=1.2$, and Theorem~\ref{thm:nonlinear-exact-tail} gives $0.2<a<1/3$.
The projected-only sufficient condition $p>2-1/t=4/3$ does not apply to this example.
At a fixed positive stepsize, any nondegenerate invariant distribution has infinite $1.25$th moment: conditional on any nonzero state, the next update contains a nonzero multiple of $\zeta$, while the additive innovation has a finite $1.25$th moment.
This observation concerns remote tails at fixed stepsize; the theorem concerns the stated joint limit.

\subsection{When the expected objective is infinite}\label{app:infinite-risk}
For additive innovations, $X_n=X_{n-1}-\gamma m(x^\star+X_{n-1})+\gamma\eta_n$.
Regular variation with $t<2$ gives $\E|\eta_n|^2=\infty$.
Conditional on the past, the other term is a fixed finite vector, so $\E(|X_n|^2\mid\mathcal F_{n-1})=\infty$ for every $n\ge1$.
Since centered dissipativity implies $M(x^\star+x)-M(x^\star)\ge\alpha|x|^2/2$, the expected excess objective is infinite.
The tail probabilities and quantiles in Corollary~\ref{cor:risk-tail} remain finite and quantify the error distribution.
\section{Shrinking-threshold and objective-value experiments}\label{app:diagonal}

\paragraph{Design and simulation unit.}
The study uses the scalar model in~\eqref{eq:scalar-diagonal-model}.
For each of the three stepsizes, we generate three batches of $200{,}000$ independent trajectories, totaling $600{,}000$ terminal pairs.
Each trajectory starts at zero and follows
\[
 X_{k+1}=X_k-\gamma m(X_k)+\gamma\sigma\varepsilon_{k+1},
 \qquad Y_{k+1}=(1-\gamma)Y_k+\gamma\sigma\varepsilon_{k+1}.
\]
The innovations are shared within a pair and independent between trajectories.
Horizon lengths are $2303$, $9682$, and $34539$ for $\gamma=0.01$, $0.003$, and $0.001$, respectively.
Both $a=0.2$ and $a=0.3$ use the same saved terminal states; there is no separate simulation per threshold.
The simulation design, sample counts, seeds, and two threshold exponents were fixed before production sampling.
No stopping rule depends on exceedance counts or observed calibration.

Figure~\ref{fig:threshold-sweep} adds a descriptive sweep of the same saved endpoints over 241 equally spaced thresholds in $[0.005,0.1]$.
At each threshold, all three exceedance counts are recomputed directly from the saved paired endpoints, without smoothing or interpolating empirical probabilities.
The common plotting range is $u/\sigma\in[0.25,5]$ and the common vertical range includes every pointwise interval.
The smallest nonlinear exceedance count across this grid is 360.

\paragraph{Structural certificates.}
With $h(x)=\sin^3(x)/(1+x^2)$, $xh(x)\ge0$ for $|x|\le\pi$ and
$h(x)/x\ge-1/\{\pi(1+\pi^2)\}$ outside that interval.
Consequently $xm(x)\ge[1-4/\{\pi(1+\pi^2)\}]x^2>0.882x^2$.
The derivative at $x=2$ is approximately $-0.307$, establishing nonconvexity.
The population remainder at zero is cubic.
Bounded $h'$ and Student moments give the global growth and local increment assumptions for every $1<p<1.5$.
Choosing $p=1.45$ makes the quadratic-remainder window $1/30<a<1/3$, containing both design exponents.

\paragraph{Deterministic denominators.}
Writing $\overline T_t$ for the Student survival function,
\[
 \Lambda_\gamma(u)=\sum_{j=0}^{n-1}
 \overline T_{1.5}\left(\frac{u}{\gamma\sigma(1-\gamma)^j}\right),
 \qquad
 A_\gamma(u)=\frac{\overline T_{1.5}(u/(\gamma\sigma))}{1.5\gamma}.
\]
Both quantities are evaluated using the exact survival function.
The second expression is the theorem's first-order formula, still retaining the exact target-innovation tail.
The reported $R_{\rm const}$ isolates the finite-sum correction.

\paragraph{Intervals and coupling.}
Let $K_X$, $K_Y$, and $K_{XY}$ count the events $X>u$, $Y>u$, and their intersection among $N$ independent pairs.
The intervals for $R_{\rm tot}$ and $R_{\rm lin}$ divide a Wilson interval by the deterministic $\Lambda_\gamma$.
For the transfer ratio, set $\widehat p_X=K_X/N$, $\widehat p_Y=K_Y/N$, and $\widehat p_{XY}=K_{XY}/N$.
The delta-method variance of its logarithm is
\[
 \widehat v=\frac1N\left(\frac1{\widehat p_X}
       +\frac1{\widehat p_Y}
       -\frac{2\widehat p_{XY}}{\widehat p_X\widehat p_Y}\right).
\]
We use $R_{\rm tr}\exp(\pm1.96\sqrt{\widehat v})$.
This keeps the within-pair covariance and uses $N$ independent units.
The intervals across different thresholds are dependent and are interpreted pointwise.

\begin{table}[htbp]
\centering
\caption{Complete prespecified shrinking-threshold results. $K_X$ is the one-sided nonlinear exceedance count among $600{,}000$ independent endpoints. The final column is the 95\% pointwise interval for $R_{\rm tot}$.}
\label{tab:diagonal}
\begin{tabular}{rrrrrrr}
\toprule
$\gamma$ & $a$ & $K_X$ & $R_{\rm lin}$ & $R_{\rm tr}$ & $R_{\rm tot}$ & 95\% CI\\
\midrule
0.01 & 0.2 & 2529 & 1.034 & 0.920 & 0.951 & [0.914, 0.988]\\
0.01 & 0.3 & 5420 & 1.069 & 0.956 & 1.021 & [0.994, 1.048]\\
0.003 & 0.2 & 2030 & 1.043 & 0.932 & 0.972 & [0.931, 1.016]\\
0.003 & 0.3 & 5192 & 1.071 & 0.972 & 1.041 & [1.013, 1.069]\\
0.001 & 0.2 & 1588 & 0.987 & 0.961 & 0.948 & [0.903, 0.996]\\
0.001 & 0.3 & 4883 & 1.054 & 0.982 & 1.034 & [1.006, 1.064]\\
\bottomrule
\end{tabular}
\end{table}

\paragraph{Exact objective event.}
The objective is even and strictly increasing on the positive half-line.
For each $u$, we solve $M(r_u)=u^2$ using adaptive quadrature and a bracketed scalar root solver, and evaluate the event as $|X_n|>r_u$.
This is an exact reformulation of the nonlinear objective event up to the reported root-solver tolerance.
With $H=1$ and symmetric scalar innovations, $\nu_H\{|z|>r\}=r^{-t}/t$.
Hence $C_M=2^{-t/2}/t$ and the denominator is
$A_{M,\gamma}=2^{-t/2}\Pp(|\sigma\varepsilon|>u/\gamma)/(t\gamma)$.
The risk intervals use the binomial Wilson construction with the actual $N$.

\begin{figure}[htbp]
\centering
\includegraphics[width=\linewidth]{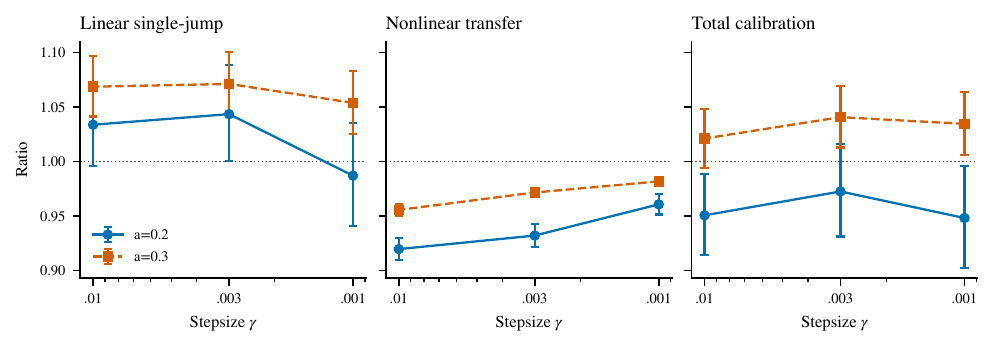}
\caption{The two prespecified shrinking-threshold sequences. Each point uses $600{,}000$ independent terminal pairs; stepsizes decrease from left to right. Error bars are pointwise probability or paired-ratio intervals.}
\label{fig:diagonal-main}
\end{figure}

\begin{figure}[htbp]
\centering
\includegraphics[width=\linewidth]{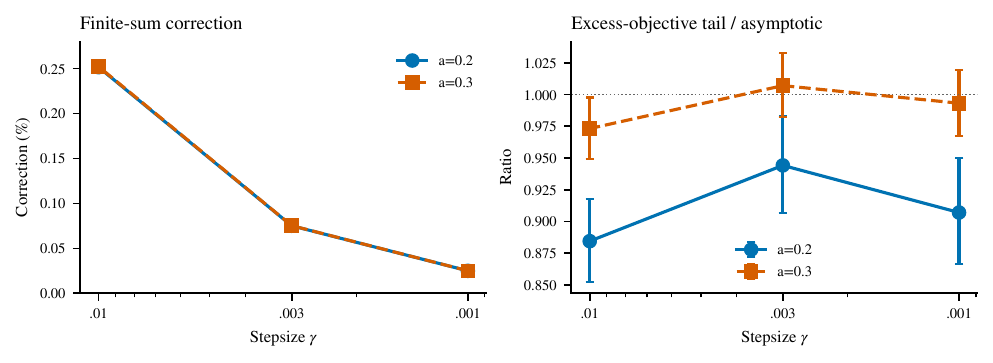}
\caption{Finite-sum correction and objective-tail calibration for both prespecified exponents. The left panel reports the percentage correction $100(R_{\rm const}-1)$; the right panel retains the full range of objective-tail ratios and pointwise intervals.}
\label{fig:diagonal-risk}
\end{figure}

Figures~\ref{fig:diagonal-main} and~\ref{fig:diagonal-risk} give the complete transfer, constant, and objective comparisons for the two prespecified sequences.

There are $1588$--$5420$ one-sided nonlinear exceedances among $600{,}000$ independent endpoints in the six settings.
All six paired transfer intervals lie below one, with point estimates moving toward one as the stepsize decreases.
For $a=0.2$, the total ratios are $0.951$, $0.972$, and $0.948$; for $a=0.3$, they are $1.021$, $1.041$, and $1.034$.
Thus the total error combines the opposite signs of the nonlinear-transfer and linear single-jump discrepancies in several settings.
The exact-CDF sum divided by the asymptotic constant ranges from $1.000250$ to $1.002522$.
Its small discrepancy cannot explain the larger observed differences from one.

The exact objective-event counts at decreasing stepsizes are $2791$, $2342$, and $1806$ for $a=0.2$, and $6128$, $5971$, and $5574$ for $a=0.3$.
The corresponding objective-tail ratios are $0.884$, $0.944$, and $0.907$, and $0.973$, $1.007$, and $0.993$.
The maximum absolute residual of the numerical equation $M(r_u)=u^2$ is $1.3\times10^{-16}$.
All counts, denominators, inverse thresholds, and intervals are supplied as machine-readable tables.
The objective result is accurate for the faster-shrinking sequence at the smallest step, while the slower-shrinking sequence retains a visible finite-threshold discrepancy.

\paragraph{Saved artifacts.}
The supplement includes every terminal pair, configuration, per-batch seed, source hash, elapsed time, aggregate count, and figure-generation script.
The generator uses NumPy's PCG64 random generator; seeds are $2026092607+100i+j$, where $i=0,1,2$ indexes the stepsize and $j=0,1,2$ the batch.
The simulator was checked against direct scalar updates and a closed-form weighted linear sum on identical deterministic innovations.
Analysis checks reconcile count identities, Wilson constructions, paired-ratio covariance, and the objective inversion.

\subsection{Random-curvature ablation}\label{app:random-curvature}

Let $A_k=1+10\zeta_k$, where $\zeta_k$ is Rademacher and independent of the Student innovation.
The model in the main text has an integrable sample loss
$G(x;\zeta,\varepsilon)=\frac12(1+10\zeta)x^2+4\Phi(x)-0.02\varepsilon x$,
where $\Phi'(x)=\sin^3x/(1+x^2)$.
Its expectation is the scalar objective.
Since $A_k$ is bounded, all growth and increment assumptions hold for every $1<p<1.5$.
The target innovation and local population Hessian are the same as in the additive study.
The state-dependent martingale term is nonzero.

The two linear controls are
\[
 Q_{k+1}=(1-\gamma A_{k+1})Q_k+\gamma\sigma\varepsilon_{k+1},
 \qquad
 Y_{k+1}=(1-\gamma)Y_k+\gamma\sigma\varepsilon_{k+1}.
\]
They share $(A_k,\varepsilon_k)$ with the full nonlinear chain.
All three start at zero.
The same three stepsizes, two threshold exponents, and horizon schedule as the additive study are used, with $N=200{,}000$ independent triples per stepsize.
Two batches of $100{,}000$ are generated with seeds $2026092611+100i+j$.
The configuration is fixed before sampling and all endpoints are retained.
The homogeneous random-curvature update satisfies
$\E(1-\gamma A)^2=1-2\gamma+101\gamma^2<1$ at all tested stepsizes.
This check concerns the update coefficient and does not assume a finite second moment of the driven process.

Define
$R_{\rm curv}=\widehat p_Q/\widehat p_Y$ and
$R_{\rm drift}=\widehat p_X/\widehat p_Q$.
Then
$R_{\rm tot}=R_{\rm lin}R_{\rm curv}R_{\rm drift}$ exactly at every displayed point.
All ratio intervals keep the paired covariance as described above.
The deterministic benchmark remains the population-Hessian exact-CDF sum.

\paragraph{A coefficient-based diagnostic.}
For the random linear control, the coefficient product of a shock of age $j$ is
$W_j=\prod_{i=1}^j(1-\gamma A_i)$.
All these factors are positive at the tested stepsizes, and independence gives
$\E W_j^t=m_t(\gamma)^j$, with
$m_t(\gamma)=\{(1-11\gamma)^t+(1+9\gamma)^t\}/2$.
The ratio of expected sums of the $t$th powers of the random and deterministic weights is therefore exactly
\[
 J_{\gamma,n}
 =\frac{1-(1-\gamma)^t}{1-m_t(\gamma)}
   \frac{1-m_t(\gamma)^n}{1-(1-\gamma)^{nt}}.
\]
For the tested horizons and $t=1.5$, this equals $1.337$, $1.081$, and $1.026$ at the three decreasing stepsizes.
It measures the amplification of a power-tail single-shock contribution by the random coefficients.
The dotted curve in the middle panel is this coefficient diagnostic, not an exact finite-step probability formula.
When the horizon factors are negligible, its small-$\gamma$ expansion is $1+25\gamma+O(\gamma^2)$.
This expansion is an algebraic property of the weights; no finite-step error rate for the nonlinear probability is inferred from it.

\begin{table}[htbp]
\centering
\caption{Complete random-curvature ablation. $K_X$ counts nonlinear exceedances among $200{,}000$ independent endpoints. The final column is the pointwise interval for $R_{\rm tot}$. The two middle ratios isolate random curvature and the smooth population remainder.}
\label{tab:random-curvature}
\begin{tabular}{rrrrrrr}
\toprule
$\gamma$ & $a$ & $K_X$ & $R_{\rm tot}$ & $R_{\rm curv}$ & $R_{\rm drift}$ & 95\% CI\\
\midrule
0.01 & 0.2 & 1206 & 1.360 & 1.428 & 0.953 & [1.286, 1.439]\\
0.01 & 0.3 & 2550 & 1.441 & 1.438 & 0.966 & [1.387, 1.498]\\
0.003 & 0.2 & 724 & 1.041 & 1.079 & 0.951 & [0.968, 1.119]\\
0.003 & 0.3 & 1930 & 1.160 & 1.101 & 0.984 & [1.110, 1.213]\\
0.001 & 0.2 & 569 & 1.019 & 1.080 & 0.963 & [0.939, 1.106]\\
0.001 & 0.3 & 1690 & 1.074 & 1.039 & 0.985 & [1.024, 1.126]\\
\bottomrule
\end{tabular}
\end{table}

\begin{figure}[htbp]
\centering
\includegraphics[width=\linewidth]{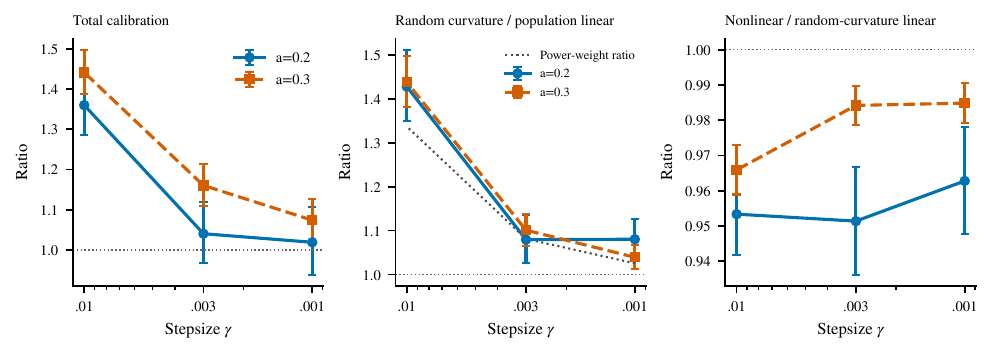}
\caption{Random-curvature ablation with both prespecified threshold exponents. Left: total calibration against the population-Hessian single-jump sum. Middle: the random-curvature linear chain relative to the population linear chain. Right: the nonlinear chain relative to the random-curvature linear chain. Error bars are pointwise and use the common-innovation covariance where appropriate.}
\label{fig:random-curvature}
\end{figure}

The $a=0.2$ random-curvature ratios at decreasing stepsizes are $1.428$, $1.079$, and $1.080$; their paired intervals are $[1.348,1.512]$, $[1.026,1.135]$, and $[1.036,1.126]$.
The last two estimates do not establish a strictly monotone finite-step trend, although both are much smaller than the largest-step amplification.
At $\gamma=0.001$, the population-remainder ratio is $0.963$ ($[0.948,0.978]$), while the population-linear single-jump ratio is $0.980$ ($[0.901,1.065]$).
Their product with the curvature ratio yields the total $1.019$ ($[0.939,1.106]$).
For $a=0.3$, the random-curvature ratios are $1.438$, $1.101$, and $1.039$, and the population-remainder ratios are $0.966$, $0.984$, and $0.985$.
At the smallest step, the total $1.074$ ($[1.024,1.126]$) remains detectably above one.
The $t$th-weight diagnostic approaches one more rapidly than some empirical ratios; it isolates a coefficient effect and does not include accumulated small jumps or nonlinear drift.
There are between $569$ and $2550$ one-sided nonlinear exceedances per setting, so the uncertainty at the rarer thresholds is visible rather than suppressed by plotting only point estimates.
\section{Objective tails, quantiles, and batch effects}\label{app:merged-risk}

\paragraph{Model and sampling units.}
In dimension $d=10$, let $H=\operatorname{diag}(1,1+1/(d-1),\ldots,2)$, $\boldsymbol a=(1,\ldots,1)^\top/\sqrt d$, and
\begin{align*}
 h(s)&=\frac{\sin^3s}{1+s^2},\qquad \Phi(s)=\int_0^s h(v)\,\dd v,\\
 m(x)&=Hx+6\boldsymbol a h(\boldsymbol a^\top x),\qquad
 M(x)=\tfrac12x^\top Hx+6\Phi(\boldsymbol a^\top x).
\end{align*}
The innovation is the multivariate Student vector $\eta=V/\sqrt{Q/t}$, with $V\sim N(0,I_d)$, $Q\sim\chi_t^2$, $t=1.5$, and independent $V,Q$.
Each step averages $b$ independent innovations.
The configurations, in seed order, are $(\gamma,b)=(0.05,1),(0.01,1),(0.04,4),(0.01,4),(0.01,16),(0.002,1)$.
Each uses $N=300{,}000$ independent terminal pairs, $X_0=Y_0=0$, and
$n=\max\{500,\lceil5\log(1/\gamma)/\gamma\rceil\}$.
The nonlinear and population-linear chains share every innovation within a pair.
Seeds are $20260926+i$ for configuration index $i=0,\ldots,5$.

\paragraph{Objective evaluation.}
The integral $\Phi$ is even and bounded.
We tabulate it on $[0,64]$ by composite Simpson integration with mesh $0.002$ and use cubic Hermite interpolation with the known derivative $h$.
For larger arguments, $\sin^3s=(3\sin s-\sin3s)/4$ gives a tail expansion after eight integrations by parts.
Adaptive weighted quadrature at 59 validation points up to $10^4$ differs by less than $10^{-12}$ in $\Phi$.
This procedure avoids applying a fixed quadrature rule across arbitrarily many oscillations.
All comparisons, including the Gaussian surrogate and the coupled linear chain, use the same full objective $M$.

\paragraph{Three distinct predictions.}
Write $\eta=R\Theta$, where $\Theta$ is uniform on $S^{d-1}$, independent of $R$, and $R^2/d\sim F(d,t)$.
Let $B_\gamma=I_d-\gamma H$.
We distinguish
\begin{align}
 A_b(\epsilon)&=b^{1-t}C_M\frac{\Pp(R>\sqrt\epsilon/\gamma)}{\gamma},\label{eq:risk-local-benchmark}\\
 \Lambda_b^{\rm quad}(\epsilon)&=b\sum_{j<n}\Pp\!\left\{\tfrac12(\gamma B_\gamma^j\eta/b)^\top H(\gamma B_\gamma^j\eta/b)>\epsilon\right\},\label{eq:risk-quadratic-benchmark}\\
 B_b(\epsilon)&=b\sum_{j<n}\Pp\{M(\gamma B_\gamma^j\eta/b)-M(0)>\epsilon\}.\label{eq:risk-exact-benchmark}
\end{align}
The first is the local quadratic tail-measure approximation.
The second uses the exact distribution of an individual gradient shock and a quadratic event.
The third evaluates the full objective at the linearly transported shock.
The factor $b$ counts the possible extreme sample within a batch; each transported sample has weight $\gamma/b$.
For $b>1$, these sums concern the individual independent samples in the unrolled recursion, rather than replacing the average innovation by a Student variable.
The theorem controls the local approximation in its stated shrinking-tolerance regime; $B_b$ is also assessed empirically at the finite objective levels used here.
The local prediction uses the same propagated constant as the stable-limit extrapolation; the theorem establishes its relative accuracy in the specified regime.

\paragraph{Numerical benchmark and its uncertainty.}
For each configuration, we draw $200{,}000$ independent pairs of shock age and spherical direction, with seed $20270925+i$.
The age is a geometric draw with success probability $1-e^{-1.2\gamma}$, capped at $n-1$; its terminal mass includes the entire capped tail.
We weight by the inverse age probability and integrate the radial Student tail analytically.
For $B_b$, the radial threshold is found by bracketed Newton iteration.
It is unique because the model's centered dissipativity makes $M(rz)$ strictly increasing for $r>0$ and $z\ne0$.
The solver's scaled residual is below $3\times10^{-11}$ over the tested grid.
Standard errors are computed across the independent age--direction draws, retaining all contributions from a draw together.
The constant $C_M$ is evaluated by angular Monte Carlo and 80-node Gauss--Legendre integration on $[0,25]$; the omitted time tail is exponentially small.
Its estimate and standard error, as well as the two finite-sum predictions, are supplied for every threshold.

\paragraph{Empirical uncertainty and quantiles.}
All tail intervals are pointwise 95\% Wilson intervals for independent terminal trajectories.
Ratios between independent configurations use the delta-method log-ratio variance
$(1/p_1-1)/N_1+(1/p_2-1)/N_2$.
Within a coupled pair, the joint exceedance count supplies the covariance correction.
Empirical quantile intervals use binomial order-statistic ranks: for level $1-\delta$, the lower and upper ranks are the $0.025$ and $0.975$ quantiles of $\operatorname{Binomial}(N,1-\delta)$, with the upper rank increased by one.
Predicted quantiles invert a monotone log--log interpolation of the corresponding survival curve.
The 181 objective levels run from $0.1$ to $10^5$; the quantile display uses 63 tail levels between $5\times10^{-4}$ and $0.05$.
The figure's quantile bands reflect empirical quantile uncertainty conditional on the independently estimated benchmark; the latter's standard errors are separately recorded.
On the count-qualified rare-event cells, the largest relative benchmark standard error is $0.117\%$.
Evaluating the survival prediction directly at each interpolated 99\% and 99.9\% quantile gives relative inversion error below $2.7\times10^{-7}$.

\paragraph{Covariance-based Gaussian surrogate.}
For each of five independent seeds, we estimate
$\widehat\Sigma=K^{-1}\sum_{k=1}^K\eta_k\eta_k^\top$, $K=10^5$, and solve the discrete Lyapunov covariance of the linear recursion,
\[
 (\widehat\Sigma_X)_{ij}=\frac{\gamma^2(\widehat\Sigma/b)_{ij}}{1-(1-\gamma H_{ii})(1-\gamma H_{jj})}.
\]
We draw two million Gaussian vectors per fitted covariance and evaluate their full objective $M$.
The figures and tables report the median across the five fits and retain their full range in the supplementary outputs.
Common standard-normal draws are reused across stepsize/batch configurations within a fit; the five covariance fits are independent.
Because the true innovation covariance is infinite, this is a deliberately misspecified covariance-based diagnostic.
The matched heavy-tail reference is the coupled linear SGD chain, evaluated using the same objective and innovations as the nonlinear chain.
Quadratic-only Gaussian quantiles are also retained to separate covariance misspecification from the choice of objective.

\begin{figure}[htb]
\centering
\includegraphics[width=\linewidth]{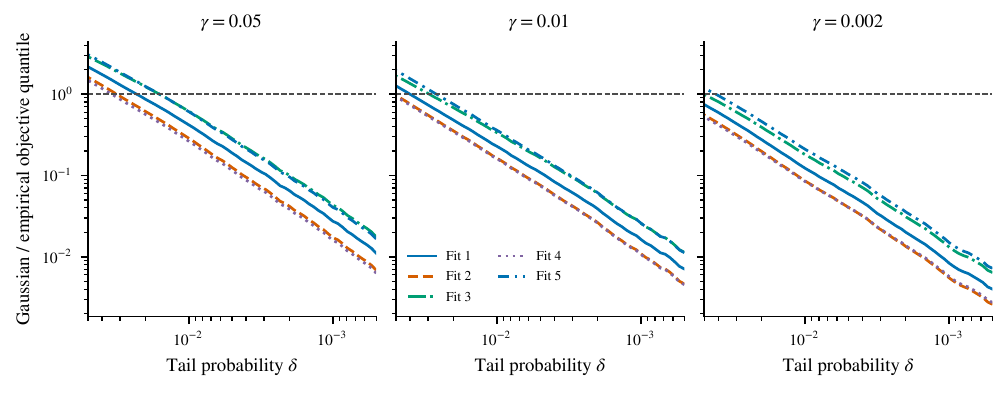}
\caption{Sensitivity of the covariance-based Gaussian diagnostic. Every curve is one of five independent covariance fits; all predictions use the full objective. The denominator is the empirical SGD quantile from $300{,}000$ terminal trajectories at the corresponding stepsize. The logarithmic vertical scale shows both fit variability and the discrepancy at increasingly extreme quantiles.}
\label{fig:gaussian-sensitivity}
\end{figure}

\begin{table}[htb]
\centering\small
\caption{All configurations at the 99\% and 99.9\% objective quantiles. ``Linear'' uses the coupled heavy-tail chain with the full objective. ``Local'', ``Quad.\ jump'', and ``Exact-loss jump'' invert~\eqref{eq:risk-local-benchmark}--\eqref{eq:risk-exact-benchmark}, respectively.}
\label{tab:all-risk-quantiles}
\begin{tabular}{cccccccc}
\toprule
$\gamma$ & $b$ & $\delta$ & Empirical [95\% CI] & Linear & Local & Quad. jump & Exact-loss jump\\
\midrule
0.05 & 1 & 0.001 & 3026.1 [2582.4, 3509.2] & 3021.5 & 2881.4 & 2953.8 & 2956.1\\
0.05 & 1 & 0.01 & 136.1 [129.3, 141.9] & 136.0 & 133.7 & 137.1 & 139.0\\
0.01 & 1 & 0.001 & 988.0 [835.7, 1145.7] & 983.1 & 985.6 & 990.8 & 993.0\\
0.01 & 1 & 0.01 & 48.7 [46.7, 51.0] & 49.0 & 45.7 & 46.0 & 47.6\\
0.04 & 4 & 0.001 & 1002.2 [858.1, 1225.1] & 1002.0 & 985.5 & 1006.1 & 1008.3\\
0.04 & 4 & 0.01 & 48.9 [46.7, 51.4] & 49.2 & 45.7 & 46.7 & 48.4\\
0.01 & 4 & 0.001 & 399.5 [356.8, 464.7] & 399.5 & 391.1 & 393.2 & 395.3\\
0.01 & 4 & 0.01 & 18.9 [18.0, 19.9] & 19.4 & 18.2 & 18.3 & 19.5\\
0.01 & 16 & 0.001 & 166.9 [146.8, 188.7] & 166.8 & 155.2 & 156.0 & 158.0\\
0.01 & 16 & 0.01 & 7.7 [7.3, 8.1] & 8.2 & 7.2 & 7.2 & 7.9\\
0.002 & 1 & 0.001 & 380.5 [320.9, 448.9] & 380.6 & 337.1 & 337.4 & 339.4\\
0.002 & 1 & 0.01 & 17.8 [17.0, 18.7] & 18.2 & 15.6 & 15.7 & 16.8\\
\bottomrule
\end{tabular}
\end{table}

\subsection{Repeated plug-in pilots}\label{app:plugin-pilot}
\begin{figure}[htb]
\centering
\includegraphics[width=\linewidth]{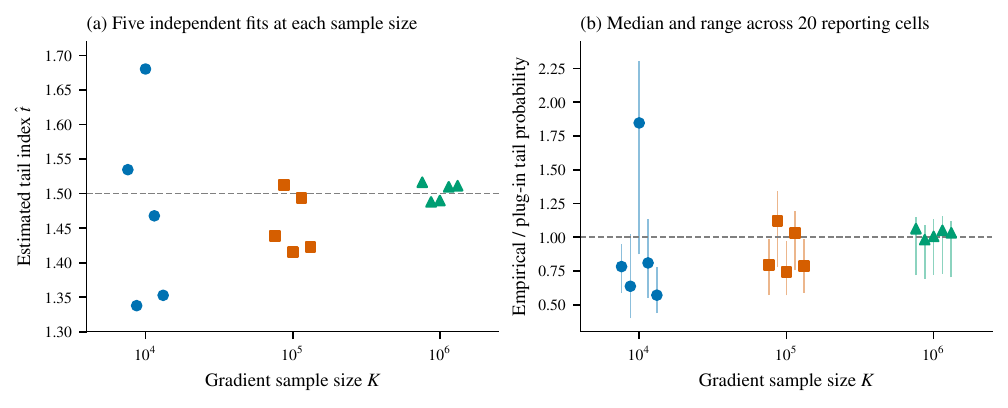}
\caption{Repeated plug-in estimation. Each marker is one independent gradient-sample fit; horizontal offsets separate the five fits at a common sample size. In (b), markers show medians and vertical segments the minimum and maximum over four directions and five thresholds. Segments describe variation across reporting cells, not confidence intervals.}
\label{fig:plugin-repeated}
\end{figure}

We also examine estimation of the projected constant from stochastic-gradient samples.
The pilot point is the coordinatewise median of 2,000 independent terminal states from $(\gamma,b)=(0.01,1)$.
At each sample size $K\in\{10^4,10^5,10^6\}$, five independent gradient samples yield a centered innovation sample, a Hill estimate from the largest 1\% of norms, and the empirical angular measure of those extremes.
A common-random-number finite difference estimates the Jacobian at the pilot point; 80-node quadrature on $[0,30]$ evaluates the estimated propagated constant.
Predictions are compared at thresholds $1,2,3,5,8$ in the four directions used in Appendix~\ref{app:multidim-experiments}.
All pilot fits and empirical exceedance counts are retained.
The additive model cancels sampling error in the finite-difference Jacobian, while the difference between curvature at the pilot point and at the target remains.
These repetitions describe sensitivity to gradient-sample size and tail estimation; they do not supply an estimator-validity or coverage theorem.
Across the five fits, the median empirical/predicted ratio over the 20 reporting cells ranges from $0.571$ to $1.847$ at $K=10^4$, from $0.745$ to $1.123$ at $10^5$, and from $0.983$ to $1.064$ at $10^6$.
The larger pilot stabilizes the median calibration, although individual directions and thresholds retain finite-level discrepancies.
Figure~\ref{fig:plugin-repeated} displays every fit and the full range across reporting cells.

\clearpage
\section{Multidimensional Projected-Tail Experiments}
\label{app:multidim-experiments}

We examined projected-tail calibration in $d=10$ and $d=100$ for two coupled nonconvex models: additive heavy-tailed gradients (Model A) and random-design regression SGD (Model B). The latter includes a random sample Hessian and state-dependent gradient noise. These experiments test the finite-step one-big-jump approximation along three eigendirections and a dense non-eigendirection.

\subsection{Population objective and stochastic-gradient models}
\label{app:multidim-models}

For each $d\in\{10,100\}$, define
\begin{equation}\label{eq:multidim-design}
 H=\operatorname{diag}\left(1+\frac{i-1}{d-1}\right)_{i=1}^d,
 \qquad \boldsymbol a=\frac{(1,\ldots,1)^\top}{\sqrt d},
 \qquad h(s)=\frac{\sin^3(s)}{1+s^2}.
\end{equation}
Both models have population objective and gradient
\begin{equation}\label{eq:multidim-objective}
 \begin{split}
 M(x)&=\tfrac12 x^\top Hx+6\Phi(\boldsymbol a^\top x),
 \qquad \Phi(s)=\int_0^s h(v)\,\mathrm dv,\\
 m(x)&=Hx+6\boldsymbol a\,h(\boldsymbol a^\top x).
 \end{split}
\end{equation}
The nonlinear correction couples all coordinates. It is $O(|x|^3)$ near zero and bounded globally, so the same linear drift $Hx$ governs both the local and large-$|x|$ behavior. The spectrum of $H$ remains in $[1,2]$ as the dimension changes.

\paragraph{Nonconvexity and centered dissipativity.}
For $|s|\le\pi$, $s h(s)\ge0$; for $|s|\ge\pi$, $h(s)/s\ge-[\pi(1+\pi^2)]^{-1}$. Hence
\begin{equation}\label{eq:multidim-dissipativity}
 x^\top m(x)\ge
 \left(1-\frac{6}{\pi(1+\pi^2)}\right)|x|^2
 >0.824|x|^2.
\end{equation}
Thus $x^\star=0$ is the unique stationary point and global minimizer. Nevertheless,
\[
 \boldsymbol a^\top\nabla^2 M(2\boldsymbol a)\boldsymbol a
 =1.5+6h'(2)\simeq-0.4604<0,
\]
so $M$ is nonconvex. These certificates hold in both dimensions.

\paragraph{Model A: additive heavy-tailed gradients.}
Independently at each step, let $V\sim N(0,I_d)$ and $Q\sim\chi^2_{1.5}$ be independent, and set $\eta=V/\sqrt{Q/1.5}$. We used
\begin{equation}\label{eq:multidim-model-a}
 X_{k+1}=X_k-\gamma m(X_k)+\gamma\eta_{k+1}.
\end{equation}
Every unit projection of $\eta$ has the univariate Student-$t_{1.5}$ law, but its coordinates are not independent. The sample loss $G_A(x;\eta)=M(x)-\eta^\top x$ has finite expectation and population gradient $m$.

\paragraph{Model B: random-design regression SGD.}
Let $U$ be uniform on the unit sphere in $\mathbb R^d$, set $Z=\sqrt d\,H^{1/2}U$, and draw $\varepsilon\sim t_{1.5}$ independently of $U$. Independent pairs $(Z,\varepsilon)$ are sampled at successive steps. We used
\begin{align}
 G_B(x;Z,\varepsilon)&=\tfrac12(Z^\top x)^2-\varepsilon Z^\top x
                  +6\Phi(\boldsymbol a^\top x),\label{eq:multidim-model-b-loss}\\
 X_{k+1}&=X_k-\gamma\left[Z_{k+1}(Z_{k+1}^\top X_k-\varepsilon_{k+1})
                  +6\boldsymbol a h(\boldsymbol a^\top X_k)\right].
                  \label{eq:multidim-model-b}
\end{align}
Since $\E[ZZ^\top]=H$ and $\E\varepsilon=0$, $\E G_B(x;Z,\varepsilon)=M(x)$. Equation~\eqref{eq:multidim-model-b-loss} is the squared-error loss with the sample-dependent constant $\varepsilon^2/2$ removed \emph{before} taking expectations. This leaves each stochastic gradient unchanged and makes the population objective finite; the uncentered expected squared loss is infinite. Model B has the random rank-one Hessian $ZZ^\top$ at the target and couples all coordinates through both the random design and the nonlinear term.

\paragraph{Relation to the assumptions.}
For every $1<p<1.5$, both innovation laws have finite $p$th moments. The bounded derivative of $h$ and the bound $|Z|^2\le2d$ give the global $L^p$ growth and local increment conditions at each fixed $d$. The cubic population remainder satisfies the local quadratic bound. Both innovations, $\eta$ in Model A and $Z\varepsilon$ in Model B, are multivariate regularly varying with index $t=1.5$. Any $1<p<1.5$ satisfies the gradient-moment conditions, and the chosen $p$ determines the shrinking-threshold window in Theorem~\ref{thm:nonlinear-exact-tail}; the endpoint $p=1.5$ is not available. These statements do not assert dimension-uniform constants. The fixed-threshold study below diagnoses finite-step behavior.

\subsection{Step sizes, projections and simulation protocol}
\label{app:multidim-protocol}

Model A used $\gamma\in\{0.1,0.05,0.01\}$ in both dimensions. Model B used the same grid at $d=10$ and $\gamma\in\{0.01,0.005,0.001\}$ at $d=100$, preserving $\gamma d\in\{1,0.5,0.1\}$. The smaller regression step sizes account for the size of the random sample Hessian, rather than a change in the population spectrum.

More precisely, the spherical fourth-moment identity gives
\[
 \E[|Z|^2ZZ^\top]
 =\frac{d}{d+2}\{\operatorname{tr}(H)H+2H^2\}.
\]
Consequently, for deterministic $v$,
\begin{equation}\label{eq:multidim-ms-update}
 \E|(I_d-\gamma ZZ^\top)v|^2
 =v^\top\left[I_d-2\gamma H+
 \gamma^2\frac{d}{d+2}\{\operatorname{tr}(H)H+2H^2\}\right]v.
\end{equation}
The homogeneous linear update is strictly contractive in Euclidean mean square for every $v\ne0$ precisely when
\begin{equation}\label{eq:multidim-stepsize}
 0<\gamma<\frac{2(d+2)}{d\{\operatorname{tr}(H)+2\lambda_{\max}(H)\}}.
\end{equation}
This upper bound is approximately $0.1263$ for $d=10$ and $0.01325$ for $d=100$.
It characterizes mean-square contraction of the homogeneous linear update and motivates the $1/d$ stepsize scale for Model B.
We assess the heavy-tailed forced recursion through tail probabilities; its finite-step errors and admissible stepsizes can depend on dimension.

We studied the four unit directions
\[
 \ell\in\{e_1,e_{d/2+1},e_d,\boldsymbol a\}.
\]
The first three are the slow, middle and fast eigendirections of $H$ (middle indices $6$ and $51$). The dense nonlinear coupling direction $\boldsymbol a$ is not an eigenvector, so it also tests the general projected formula, beyond the scalar autoregressive special case.

Each independent trajectory started at $X_0=0$ and ran for
\begin{equation}\label{eq:multidim-horizon}
 n=\max\left\{500,\left\lceil\frac{5\log(1/\gamma)}{\gamma}\right\rceil\right\}.
\end{equation}
The horizons for $\gamma=0.1,0.05,0.01,0.005,0.001$ were respectively $500,500,2303,5299,34539$. We used $N=200{,}000$ independent terminal trajectories per model--dimension--step-size setting, except Model B at $d=100$, $\gamma=0.001$, for which $N=250{,}000$. There were $2{,}450{,}000$ trajectories in total. All generated terminal states were finite, and no trajectories were discarded. We evaluated $y=0.001,0.002,\ldots,5.000$ from the same terminal projections; no additional chains were run for individual thresholds. Samples here are independent trajectories, not correlated observations taken along one chain.

\subsection{Deterministic benchmark and Monte Carlo uncertainty}
\label{app:multidim-benchmark}

Alongside each nonlinear trajectory we simulated
\[
 Y_{k+1}=B_\gamma Y_k+\gamma\eta_{k+1},
 \qquad B_\gamma=I_d-\gamma H,\qquad Y_0=0,
\]
using the same innovations as $X_k$: $\eta$ from Model A or $Z\varepsilon$ from Model B. This reference uses the \emph{population} Hessian $H$, not the random sample Hessian. The deterministic denominator was the finite-horizon one-big-jump sum
\begin{equation}\label{eq:multidim-benchmark}
 \Lambda_\gamma(\ell,y)
 =\sum_{j=0}^{n-1}\Pp\{\gamma\ell^\top B_\gamma^j\eta>y\}.
\end{equation}
It is a sum of single-jump probabilities, not the exact probability of their sum. With $w_j=(B_\gamma^j)^\top\ell$ and Student survival function $\overline T_{1.5}$, the summands are
\begin{align*}
 \text{Model A:}\quad&\overline T_{1.5}\!\left(\frac{y}{\gamma|w_j|}\right),\\
 \text{Model B:}\quad&\E_{S}\!\left[\overline T_{1.5}\!\left(
 \frac{y}{\gamma\sqrt d\,|H^{1/2}w_j|S}\right)\right],
 \qquad S^2\sim\operatorname{Beta}\!\left(\tfrac12,\tfrac{d-1}{2}\right).
\end{align*}
The regression expectation was evaluated by 128-node Gauss--Jacobi quadrature, with shape-preserving log-scale interpolation for dense-grid evaluation. Increasing the order to 256 in the $d=100$ regression checks at $y=1$ and $5$, for all displayed directions and step sizes, changed the denominator by less than $4.2\times10^{-5}$ in relative terms.

Writing $\widehat p_X$ and $\widehat p_Y$ for the projected-tail estimates, we report
\begin{equation}\label{eq:multidim-ratios}
 R_{\rm lin}=\frac{\widehat p_Y}{\Lambda_\gamma},
 \qquad R_{\rm tr}=\frac{\widehat p_X}{\widehat p_Y},
 \qquad R_{\rm tot}=\frac{\widehat p_X}{\Lambda_\gamma}.
\end{equation}
The transfer ratio in Model B combines random-Hessian and nonlinear-correction effects; it does not isolate either effect.

For $d=10$, we retained the ordinary one-sided estimates $\widehat p_X=K_X/N$, where $K_X=\sum_{i=1}^N\one_{\{\ell^\top X_n^{(i)}>y\}}$, and analogously for $Y$. The $d=100$ experiments used exact central symmetry to reduce Monte Carlo variance. Since $m$ is odd and both chains start at zero, negating all innovations in Model A, or all residuals while holding the design sequence fixed in Model B, negates the entire trajectory. Thus
\begin{equation}\label{eq:multidim-symmetry}
 p_X=\tfrac12\Pp(|\ell^\top X_n|>y),\qquad
 \widehat p_X=\frac{K_X^{\rm abs}}{2N},\qquad
 K_X^{\rm abs}\sim\operatorname{Binomial}(N,2p_X).
\end{equation}
The variance is $p_X(1-2p_X)/(2N)$, compared with $p_X(1-p_X)/N$ for the one-sided estimator. The two signs of a trajectory are \emph{not} independent trials. Symmetry is a property of these examples, not an assumption of the general tail theorem.

Shading in the figures gives nominal 95\% \emph{pointwise} Wilson score intervals for $p_X$, divided by $\Lambda_\gamma$. For $d=100$, we first computed the interval for $K_X^{\rm abs}/N$ using $N$ trials, then divided its endpoints by two. To restrict shading to cells with adequate tail counts, we displayed it only where $\Lambda_\gamma\le0.1$ and the observed nonlinear count, linear count and corresponding expected benchmark count were all at least 200. The expected count is $N\Lambda_\gamma$ for $d=10$ and $2N\Lambda_\gamma$ for $d=100$. Curves are shown on the entire threshold grid, without smoothing. These intervals quantify Monte Carlo uncertainty only: they neither include approximation error nor provide simultaneous coverage over directions and thresholds.

\subsection{Results and scope}
\label{app:multidim-results}

\begin{figure}[htbp]
\centering
\includegraphics[width=\linewidth]{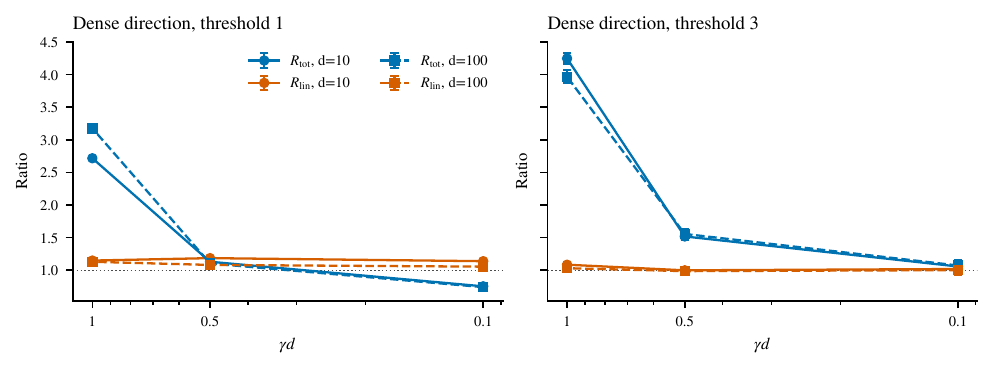}
\caption{Random-design SGD in the dense non-eigendirection, at matched $\gamma d$. Stepsizes decrease from left to right and both panels share the same vertical range. Blue curves show total calibration and orange curves the population-Hessian linear reference. Each point uses $200{,}000$ independent trajectories, except $d=100$, $\gamma=0.001$, which uses $250{,}000$.}
\label{fig:model-b-main}
\end{figure}

Table~\ref{tab:multidim-medians} uses the \emph{same deterministic reporting domain} for every row: all four directions and $y=1.000,1.001,\ldots,5.000$, comprising 16,004 direction--threshold cells. This table does not apply the count-based shading filter. The medians are descriptive summaries of correlated cells, not estimates based on 16,004 independent replications. Also, the median of $R_{\rm tot}$ need not equal the product of the other two medians, even though the pointwise identity holds.

\begin{table}[!htbp]
\centering
\caption{Median projected-tail ratios on the common domain $1\le y\le5$. $N$ is the actual number of independent trajectories. The last column is the median full width of the 95\% pointwise interval for $R_{\rm tot}$. The $d=10$ rows use one-sided counts; the $d=100$ rows use symmetry pooling.}
\label{tab:multidim-medians}
\begin{tabular}{ccrrrrrr}
\toprule
$d$ & Model & $\gamma$ & $N/10^3$ & $R_{\rm lin}$ & $R_{\rm tr}$ & $R_{\rm tot}$ & CI width\\
\midrule
10 & A & 0.100 & 200 & 1.069 & 0.992 & 1.049 & 0.085\\
10 & A & 0.050 & 200 & 1.044 & 0.992 & 1.029 & 0.101\\
10 & A & 0.010 & 200 & 1.015 & 0.993 & 0.999 & 0.150\\
10 & B & 0.100 & 200 & 1.096 & 3.923 & 4.266 & 0.158\\
10 & B & 0.050 & 200 & 1.057 & 1.533 & 1.588 & 0.117\\
10 & B & 0.010 & 200 & 1.028 & 1.051 & 1.081 & 0.147\\
\midrule
100 & A & 0.100 & 200 & 1.063 & 0.998 & 1.055 & 0.059\\
100 & A & 0.050 & 200 & 1.068 & 0.999 & 1.059 & 0.072\\
100 & A & 0.010 & 200 & 1.041 & 0.999 & 1.035 & 0.106\\
100 & B & 0.010 & 200 & 1.028 & 4.298 & 4.331 & 0.204\\
100 & B & 0.005 & 200 & 0.983 & 1.628 & 1.586 & 0.149\\
100 & B & 0.001 & 250 & 1.045 & 1.070 & 1.118 & 0.169\\
\bottomrule
\end{tabular}
\end{table}

For Model A, the median total ratios were within 6\% of one in both dimensions, and the transfer ratios were close to one. Model B showed larger finite-step discrepancies: at the matched values $\gamma d=1,0.5,0.1$, the median total ratios decreased from $4.266$ to $1.588$ to $1.081$ in $d=10$, and from $4.331$ to $1.586$ to $1.118$ in $d=100$. The similar patterns at matched $\gamma d$ support the step-size scaling motivated by~\eqref{eq:multidim-stepsize}, without establishing a dimension-uniform result.

At $d=100$ and $\gamma=0.001$ in Model B, the direction-specific median total ratios ranged from $1.078$ to $1.152$. An approximately 12\% median excess over one therefore remained at the smallest tested step size. On the common reporting domain, the median and maximum full widths of the pointwise ratio intervals were $0.169$ and $0.263$, respectively. As a sensitivity check, using ordinary one-sided counts on the same $d=100$ trajectories gave median $R_{\rm tot}=1.117$ in this setting, close to the symmetry-pooled value $1.118$.

\begin{figure}[!htbp]
\centering
\includegraphics[width=\linewidth]{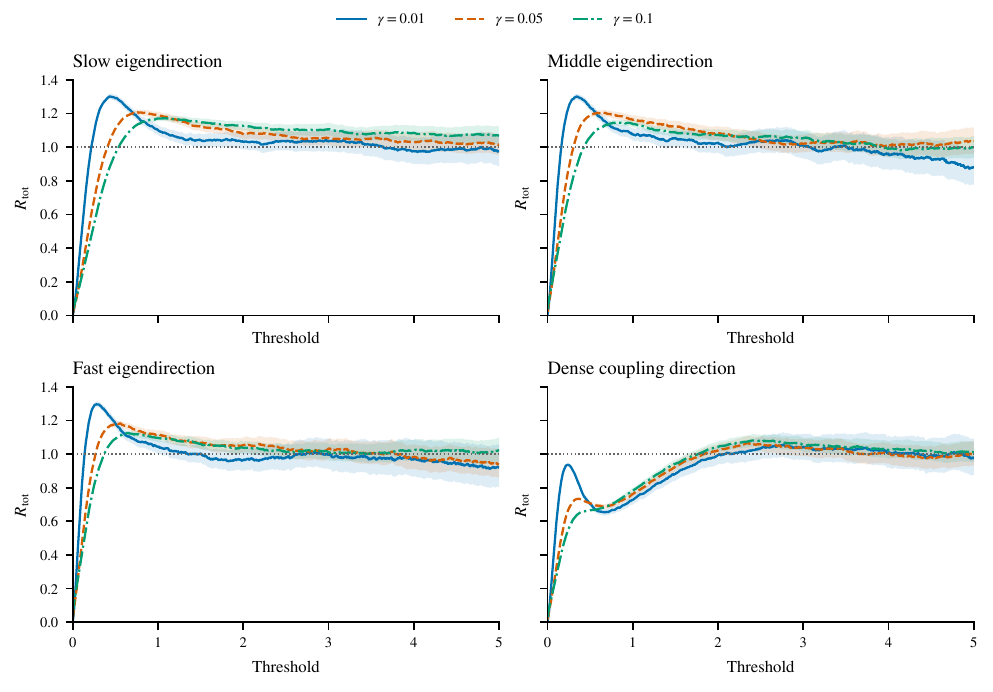}
\caption{Additive-noise SGD (Model A) in $d=10$. Curves show $R_{\rm tot}$ for $\gamma=0.1,0.05,0.01$, with $N=200{,}000$ independent trajectories per step size and ordinary one-sided tail estimates. Panels show $e_1$, $e_6$, $e_{10}$ and the dense nonlinear coupling direction $\boldsymbol a$. Shading gives 95\% pointwise Wilson ratio intervals on the count-qualified cells defined in Appendix~\ref{app:multidim-benchmark}. The dashed line marks ratio one. Axes are linear, and all four panels share a vertical scale.}
\label{fig:multidim-a10}
\end{figure}

\begin{figure}[p]
\centering
\includegraphics[width=\linewidth]{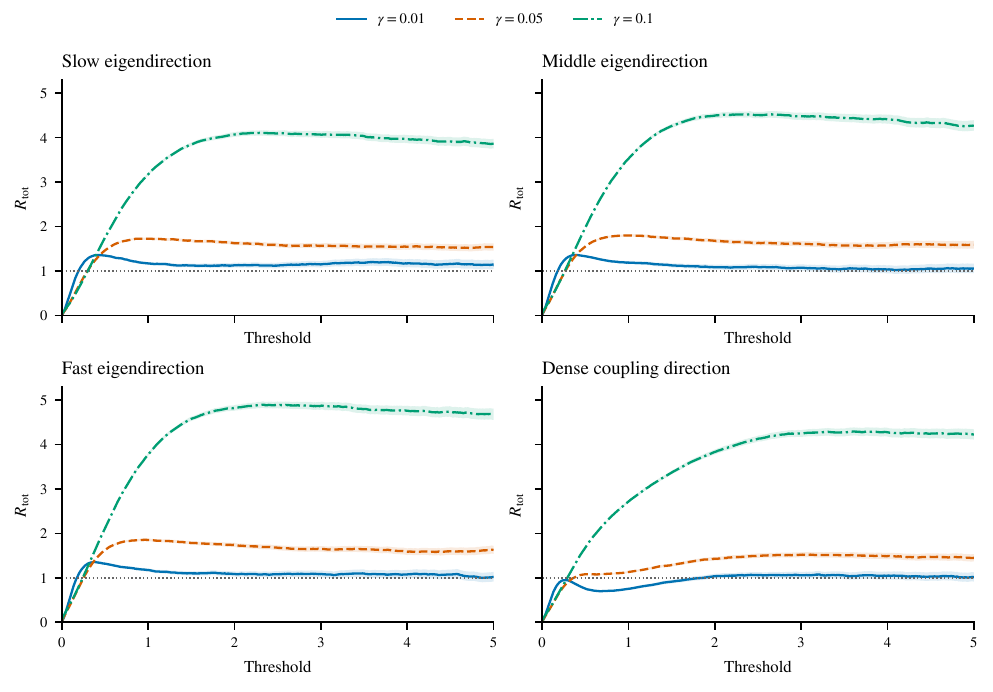}
\caption{Random-design regression SGD (Model B) in $d=10$, with $\gamma=0.1,0.05,0.01$ and $N=200{,}000$ independent trajectories per step size. The directions, one-sided estimator and pointwise interval convention are the same as in Figure~\ref{fig:multidim-a10}. The larger ratios at larger step sizes display finite-step discrepancies from the population-Hessian one-big-jump benchmark.}
\label{fig:multidim-b10}
\end{figure}

\begin{figure}[p]
\centering
\includegraphics[width=\linewidth]{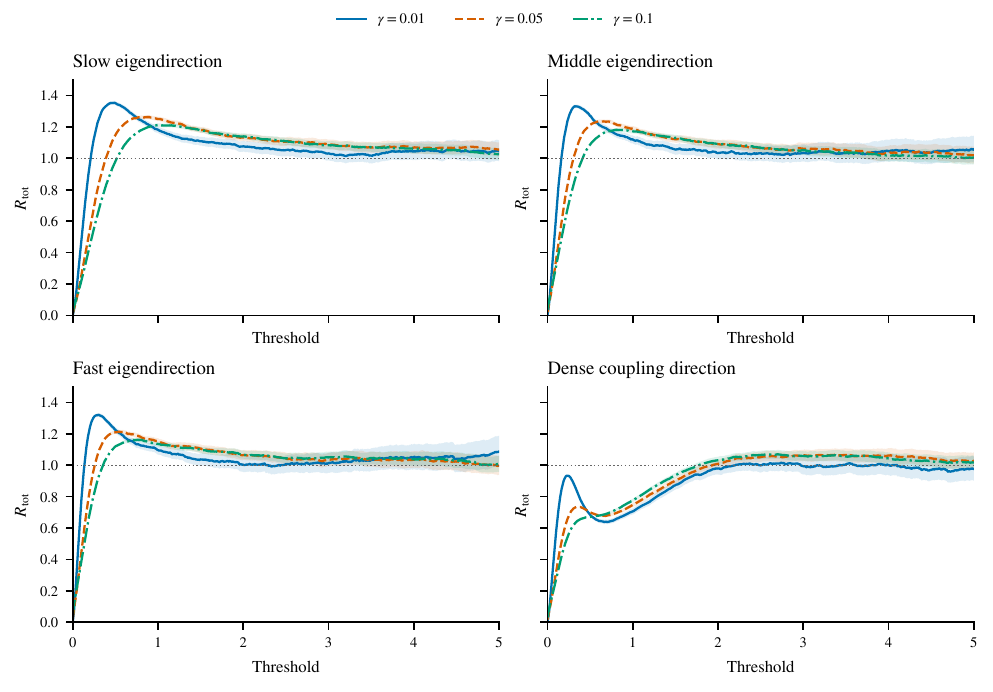}
\caption{Additive-noise SGD (Model A) in $d=100$, with $\gamma=0.1,0.05,0.01$ and $N=200{,}000$ independent trajectories per step size. Panels show $e_1$, $e_{51}$, $e_{100}$ and $\boldsymbol a$. Curves use symmetry-pooled tail estimates; shading gives 95\% pointwise Wilson ratio intervals computed with the actual $N$, as described in Appendix~\ref{app:multidim-benchmark}. Axes are linear, with a shared vertical scale across the four panels.}
\label{fig:multidim-a100}
\end{figure}

\begin{figure}[p]
\centering
\includegraphics[width=\linewidth]{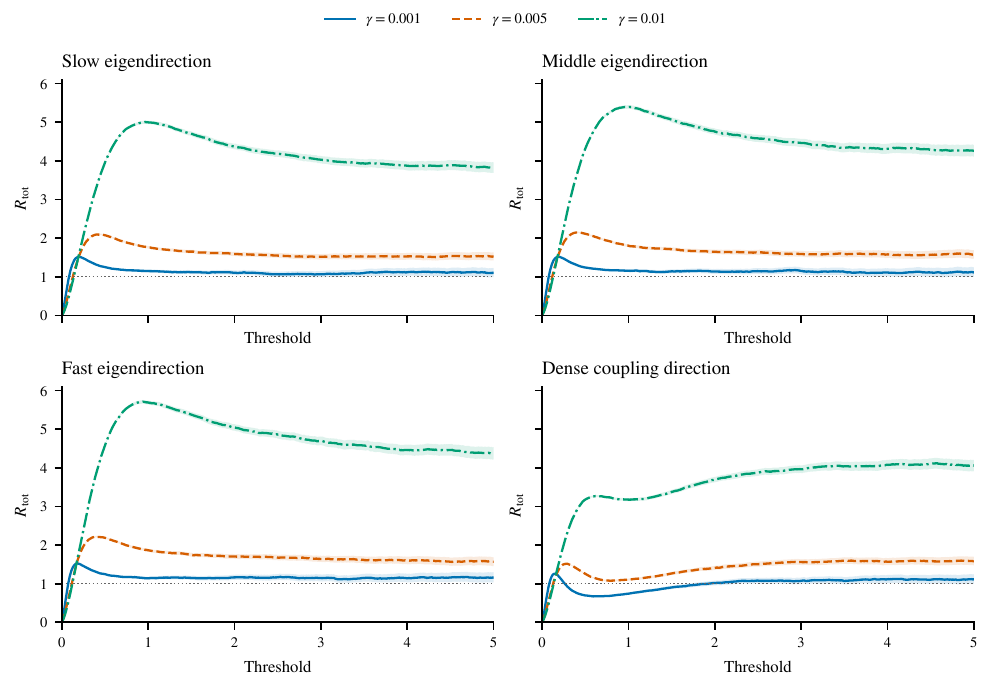}
\caption{Random-design regression SGD (Model B) in $d=100$, with $\gamma=0.01,0.005,0.001$. The first two settings use $N=200{,}000$ independent trajectories each, and the last uses $N=250{,}000$. Directions, symmetry pooling and pointwise interval conventions match Figure~\ref{fig:multidim-a100}. Decreasing the step size improves calibration, but a residual discrepancy remains at $\gamma=0.001$. The intervals measure Monte Carlo precision, not approximation error or simultaneous coverage.}
\label{fig:multidim-b100}
\end{figure}

Figures~\ref{fig:multidim-a10}--\ref{fig:multidim-b100} show the full threshold dependence. The near-origin region is not a rare-event regime, and the one-big-jump sum can exceed one there. The zero initial state and schedule~\eqref{eq:multidim-horizon} satisfy the time condition of Theorem~\ref{thm:nonlinear-exact-tail}. The fixed threshold grid serves as a finite-step diagnostic; the theorem itself uses shrinking thresholds. Together with the diagonal study, these results locate the empirical discrepancies in threshold choice, single-jump calibration, and nonlinear transfer.

\clearpage
\subsection{A random-Hessian single-jump benchmark}\label{app:random-hessian-benchmark}

For Model B, let $Z=\sqrt d H^{1/2}U$ with $U$ uniform on the sphere and scalar label noise $\epsilon\sim t_{1.5}$.
The innovation $Z\epsilon$ at the shock is independent of the subsequent sample Hessians.
Define $W_0=\ell$ and $W_{j+1}=(I_d-\gamma Z_{j+1}Z_{j+1}^\top)W_j$, with fresh designs $Z_j$.
Reversal of an iid finite product shows that $W_j$ has the distribution of the direction transported backward through $j$ sample Hessians.
The benchmark is
\begin{equation}\label{eq:random-hessian-benchmark}
 \Lambda_\gamma^{\rm RH}(\ell,u)
 =\sum_{j<n}\E\!\left[\overline T_{1.5}\!\left(\frac{u}{\gamma|W_j^\top Z|}\right)\right],
\end{equation}
where $Z$ is independent of $W_j$ and $\overline T_t$ is the Student survival function.
The population benchmark replaces $W_j$ by $B_\gamma^j\ell$.
Thus $R_{\rm tot}^{\rm RH}=\widehat p_X/\widehat\Lambda_\gamma^{\rm RH}$ separates the effect of random sample curvature from the population-flow approximation.

\paragraph{Independent sampling units.}
We use 20,000 independent design paths in $d=10$ and 8,000 in $d=100$, at all three stepsizes.
A path contributes the complete sum across ages and all four directions, so its standard error retains the dependence between ages.
We integrate over the independent shock design and Student noise analytically, reducing Monte Carlo variance.
The reported domain has four directions and 17 equally spaced thresholds on $[1,5]$, giving 68 ratios per dimension--stepsize setting.
The table's MC SE is the maximum relative standard error of the estimated denominator over this domain, expressed as a percentage.
It does not include uncertainty in the empirical SGD numerator.
Twenty independent path blocks are retained; resampling them also gives denominator-only intervals for the median ratio.

\paragraph{Analytic tail integration.}
Conditionally on $W_j$, write $v=u/(\gamma|H^{1/2}W_j|)$.
If $U_1$ is the first coordinate of a uniform sphere point, the conditional summand is
$F_d(v)=\E\overline T_t\{v/(\sqrt d|U_1|)\}$.
Expanding the Student density and integrating its tail gives
\[
 F_d(v)=\sum_{k\ge0}c_{d,k}v^{-(t+2k)},\qquad
 c_{d,k}=\frac{(-1)^k t^{t/2+k}\Gamma((t+1)/2+k)}{\sqrt\pi\Gamma(t/2)k!(t+2k)}
 \E(\sqrt d|U_1|)^{t+2k},
\]
with
\[
 \E(\sqrt d|U_1|)^r=
 \frac{d^{r/2}\Gamma(d/2)\Gamma((r+1)/2)}{\sqrt\pi\Gamma((d+r)/2)}.
\]
In this design family, $\gamma d\le1$ and $|Z|^2\le2d$ make each sample propagation matrix nonexpansive.
Consequently $|W_j|\le1$ and $td/v^2\le2td\gamma^2/u^2\le0.3$ on the reported domain.
The series is uniformly convergent and alternating there.
We keep 14 terms and bound the remainder by the first omitted term, summed over ages.
Independent adaptive spherical quadrature agrees to within $1.3\times10^{-12}$ in the validation cases.

\paragraph{Truncation and sensitivity.}
Let $h_i$ be the diagonal entries of $H$ and set
\[
 q_2=\max_i\left\{1-2\gamma h_i+
       \gamma^2\frac d{d+2}\big(\operatorname{tr}(H)h_i+2h_i^2\big)\right\}<1.
\]
The exact spherical fourth-moment formula gives
$\E|W_j|^2\le q_2^j|\ell|^2$, hence $\E|W_j|^t\le q_2^{jt/2}$.
We truncate at $J=\min\{n,\lceil40/(t|\log q_2|)\rceil\}$.
The remaining sum is bounded by
\[
 c_{d,0}(\gamma/u)^t\|H\|_{\op}^{t/2}
       \frac{q_2^{Jt/2}}{1-q_2^{t/2}},
\]
and is zero when $J=n$.
The actual cutoffs, residual bounds, and comparison with the shorter cutoff $\min\{n,\lceil16/(t\gamma)\rceil\}$ are supplied for every setting.
This controls truncation using the random-product decay rather than assuming the population contraction rate.

\begin{table}[htb]
\centering\small
\caption{Numerical accuracy of the random-Hessian denominator. Columns report the largest relative standard error, largest relative change from the shorter cutoff, and largest omitted-age bound divided by the computed denominator, over the common 68-cell domain. The residual-bound column is a fraction, whereas the two preceding columns are percentages.}
\label{tab:rh-accuracy}
\begin{tabular}{rrrrrr}
\toprule
$d$ & $\gamma d$ & $J$ & MC SE (\%) & Old-cutoff change (\%) & Residual bound\\
\midrule
10 & 1 & 444 & 0.162 & 0.22907 & 3.33e-09\\
10 & 0.5 & 400 & 0.195 & 0.00007 & 3.81e-09\\
10 & 0.1 & 1422 & 0.103 & 0.00001 & 4.07e-09\\
100 & 1 & 2303 & 0.085 & 0.28637 & 0\\
100 & 0.5 & 4237 & 0.102 & 0.00015 & 3.97e-09\\
100 & 0.1 & 14394 & 0.058 & 0.00001 & 4.08e-09\\
\bottomrule
\end{tabular}
\end{table}
The largest shorter-cutoff correction is $0.287\%$; the final omitted-age bounds are below $4.1\times10^{-9}$ relative to the calculated denominators.
The retained block summaries provide a separate check on Monte Carlo variability without treating shock ages as independent samples.
\section{Additional scalar calibration studies}\label{app:scalar-experiments}

In this section, we conduct a simulation study on three representative stochastic
recursions to examine the sharpness of the exact tail approximation in
Theorem~\ref{thm:nonlinear-exact-tail}.  Our goal is not merely to verify that
the proposed tail expression gives an upper bound, but to check whether the
empirical tail probability is well calibrated by the exact one-big-jump constant
predicted by the linearized recursion. We use $\gamma\in\{0.1,0.01,0.005,0.001\}$ and thresholds $y=0.001,0.002,\ldots,3.000$.
The denominator is the exact-CDF one-jump sum for the linearized recursion, evaluated over the finite horizon for independent terminal sampling and over the infinite horizon for stationary-chain sampling.

\paragraph{Smooth nonconvex model.}
This is~\eqref{eq:scalar-diagonal-model} with $\sigma=1$.
Each stepsize uses $10^6$ independent terminal trajectories, with horizons $1000$, $2303$, $5299$, and $34539$.
Pointwise probability intervals use the Wilson construction.

\paragraph{Threshold model.}
The drift is
\[
 m(x)=5x-2(x+\delta)\one_{\{x\le-\delta\}},
 \qquad \delta=10^{-7},
\]
with independent Student-$t_{1.5}$ additive noise.
It has slopes 3 and 5, a breakpoint at $-\delta$, and local derivative $m'(0)=5$.
The objective is strongly convex and locally quadratic near the target, while the drift has a derivative discontinuity at the breakpoint.
The local quadratic condition holds on a ball of radius $\kappa\le\delta$.
Every reported threshold is at least $10^{-3}$, far beyond that neighborhood, so this experiment probes finite-level behavior across the breakpoint.
For each stepsize, 10 independent chains run for $10^6$ steps, discard the first $10^5$ steps, and retain every tenth subsequent iterate.
The $900{,}000$ retained observations are dependent within chains.

\paragraph{Expectile estimation.}
Following \citet{NeweyPowell1987}, for $\tau\in\{0.4,0.6\}$ and $Z\sim t_{1.2}$, define
$\rho_\tau(r)=|\tau-\one_{\{r\le0\}}|r^2$ and the finite loss contrast
\[
 \theta_\tau^\star=\arg\min_\theta
 \E[\rho_\tau(Z-\theta)-\rho_\tau(Z)].
\]
Subtracting before expectation is essential because $\E Z^2=\infty$.
The stochastic gradient is
\[
 g(\theta,Z)=2[\tau\one_{\{Z>\theta\}}+(1-\tau)\one_{\{Z\le\theta\}}](\theta-Z).
\]
We use the centered recursion $X_{k+1}=X_k-\gamma g(\theta_\tau^\star+X_k,Z_{k+1})$.
The two targets are approximately $-0.831205$ and $0.831205$, and the common population derivative is $H\simeq0.907776$.
Each setting uses 100 independent chains of $10^6$ steps, discards the first $10^5$, and retains every twentieth subsequent iterate, giving $4{,}500{,}000$ observations.

\paragraph{Dependence-aware uncertainty.}
For the threshold and expectile studies, we use the Monte Carlo standard error across independent chain-level tail estimates.
We use $\widehat p\pm t_{J-1,0.975}\,\widehat{\operatorname{se}}_{\rm chain}$, clipped to $[0,1]$, with $J=10$ and $100$, respectively.
Dividing by the deterministic benchmark gives the displayed ratio intervals.
Thinning is not treated as producing independent trials.
The normal approximation across only ten chains is less reliable in the far tail; the resulting wide intervals remain visible.
These intervals quantify sampling variation across chains and do not assess residual burn-in bias.

\begin{figure}[htbp]
\centering
\includegraphics[width=\linewidth]{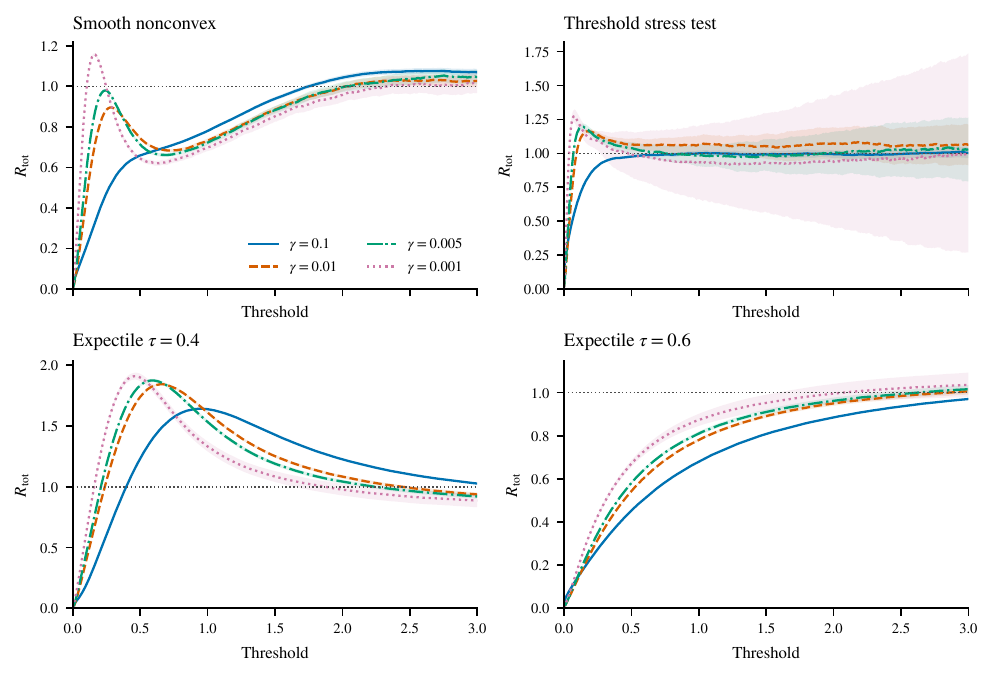}
\caption{Scalar studies with 95\% pointwise intervals. The smooth model uses independent-endpoint Wilson intervals; threshold and expectile models use independent-chain standard errors. Colors and line styles identify stepsizes. The reference line is one.}
\label{fig:legacy-scalar}
\end{figure}

The smooth model exhibits a clear threshold effect.
At $\gamma=0.001$, $R_{\rm tot}$ is $0.698$ with interval $[0.680,0.717]$ at $y=1$, but $1.016$ with interval $[0.967,1.068]$ at $y=3$.
Thus decreasing the stepsize at a fixed threshold need not yield the local-Hessian constant.
At the same stepsize and $y=3$, the $\tau=0.4$ expectile ratio is $0.887$ with interval $[0.832,0.942]$, so a systematic deficit remains.
For $\tau=0.6$, the ratio is $1.036$ with interval $[0.979,1.093]$.
The threshold model's interval is much wider, $[0.271,1.738]$, reflecting its ten independent chains.
These differences motivate separating linear and transfer errors and complement the shrinking-threshold design in the main text.

\begin{figure}[htb]
\centering
\includegraphics[width=\linewidth]{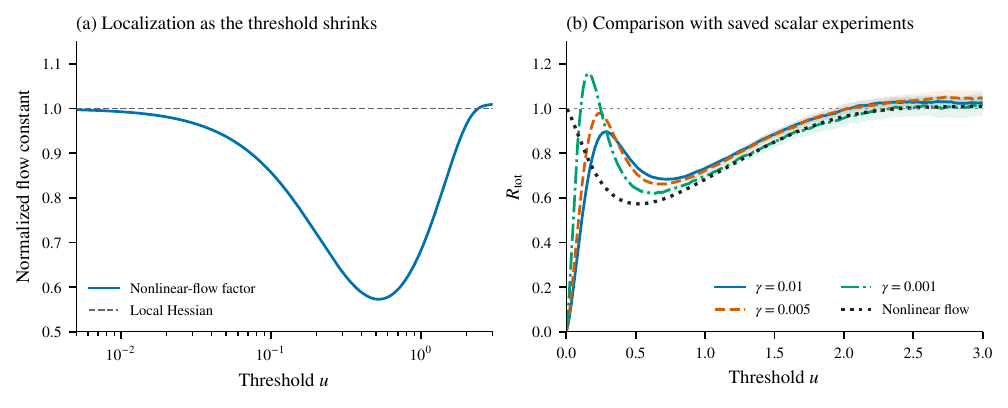}
\caption{Local curvature and nonlinear flow in the smooth model. (a) The deterministic normalized flow integral approaches one as $u\downarrow0$. (b) The same factor is overlaid on fixed-threshold SGD ratios, with independent-endpoint pointwise intervals. The selected stepsizes display the approach as $\gamma$ decreases; all four stepsizes remain in Figure~\ref{fig:legacy-scalar}.}
\label{fig:flow-localization}
\end{figure}

\paragraph{The nonlinear-flow correction.}
For the smooth scalar drift, the deterministic factor in Appendix~\ref{app:flow-localization} is
$F(u)=t u^t\int_u^\infty r^{-t}/m(r)\,\dd r$, with $t=1.5$ and $H=1$.
Figure~\ref{fig:flow-localization} shows $F(u)\to1$ at small thresholds.
Numerical integration in log coordinates and an independent physical-coordinate calculation agree within $2\times10^{-8}$ at the validation points.
The theorem proves this deterministic localization; the comparison with fixed-threshold SGD is an empirical diagnostic of the corresponding nonlinear-flow approximation.

\begin{figure}[htbp]
\centering
\includegraphics[width=\linewidth]{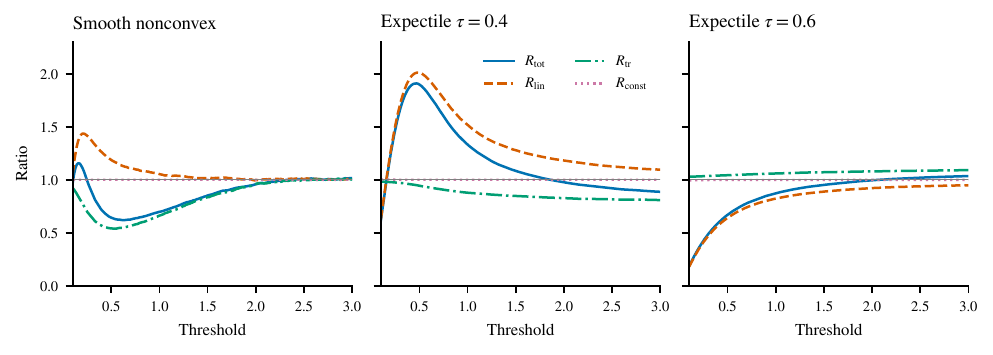}
\caption{Error decomposition of the smooth and expectile data at $\gamma=0.001$. Curves show total calibration, linear calibration, nonlinear transfer, and finite-sum versus asymptotic-constant replacement. These diagnostic curves share the same samples and are not independent estimates. Their purpose is to identify which factor contributes to a visible total discrepancy.}
\label{fig:legacy-decomposition}
\end{figure}


\end{document}